\documentclass[11pt,reqno]{amsart}  
\usepackage[dvipdfmx]{graphicx}
\usepackage{amssymb,color}                 
\allowdisplaybreaks[1]
\newtheorem{thm}{Theorem}[section]
\newtheorem{cor}[thm]{Corollary}
\newtheorem{prop}[thm]{Proposition}
\newtheorem{lem}[thm]{Lemma}
\newtheorem{assump}{\bf Assumption}

\newcommand{\sect}[1]{\section{#1}%
\setcounter{equation}{0}}
\renewcommand{\theequation}{\thesection.\arabic{equation}}
\newenvironment{prf}[1]
   {{\noindent \bf Proof of {#1}.}}{\hfill \qed}

\renewcommand{\r}{\rho}                
\newcommand{\pt}{\partial}             
\renewcommand{\th}{\theta}

\newcommand{\B}{\mathcal {B}} 
\newcommand{\D}{\mathcal{D}}
           
\newcommand{\F}{\mathcal {F}}          
\renewcommand{\S}{\mathcal {S}}

\newcommand{\Z}{\mathbb {Z}}

\newcommand{\re}{\mathbb R}

\newcommand{\Nt}{\mathbb N}
\newcommand{\N}{\nabla}

\newcommand{\Om}{\Omega}
\newcommand{\al}{\alpha}

\newcommand{\gm}{\gamma}
\newcommand{\Gm}{\Gamma}
\renewcommand{\L}{\mathcal{L}}

\newcommand{\ep}{\varepsilon}
\newcommand{\lam}{\lambda}

\newcommand{\del}{\delta}
\newcommand{\Del}{\Delta}

\newcommand{\sg}{\sigma}
\newcommand{\s}{\sigma}
\newcommand{\x}{\xi}

\renewcommand{\r}{\rho}
\renewcommand{\t}{\tau}

\newcommand{\teta}{\tilde{\eta}}

\newcommand{\dsp}{\displaystyle}

\newcommand{\dB}{\dot{B}}
\newcommand{\dF}{\dot{F}}
\newcommand{\dcB}{\dot{\mathcal{B}}}

\def\<{\langle }

\renewcommand{\>}{\rangle }

\newcommand{\bl}{\color{blue}}

\newcommand{\bk}{\color{black}}
\newcommand{\supp}{\text{ supp }}

\renewcommand{\qed}{\qquad\kern1pt   
   \vbox{\hrule height 0.6pt      
         \hbox{\vrule width 0.6pt 
               \vbox{\vskip 6pt}  
               \hskip 3pt
              \vrule width 1.3pt} 
         \hrule depth 1.3pt}     
   \kern1pt}

\newcommand{\eq}[1]{{\begin{equation}#1\end{equation}}}
\newcommand{\spl}[1]{{\begin{aligned}#1\end{aligned}}}
\newcommand{\eqn}[1]{\begin{equation*}#1\end{equation*}}

\newcommand{\algn}[1]{\begin{align*}#1 \end{align*}} 
\newcommand{\eqntag}{\addtocounter{equation}{1}\tag{\theequation}} 

\begin{document}

\title{
Maximal $L^1$-regularity 
for parabolic initial-boundary value problems with inhomogeneous data 
}
\thanks{
AMS Subject Classification: primary 35K20, secondary 42B25.\\
\quad\ \  
Keywords: parabolic equations with variable coefficients, 
maximal $L^1$-regularity, end-point estimate, 
initial-boundary value problems, 
the Dirichlet problem, 
the Neumann problem.
}

\maketitle
\vskip-1mm
 \begin{center}
{\footnotesize
{\large\sf Takayoshi Ogawa${}^{\ast}$}  and
{\large\sf Senjo Shimizu${}^{\ddag}$}
\vskip3mm
\vbox{
        Mathematical Institute/
        Research Alliance Center of Mathematical Science${}^{\ast}$\\
        Tohoku University  \\
        Sendai  980-8578, Japan  \\
        ogawa@math.tohoku.ac.jp\\
      }
\hskip5mm\noindent
\vbox{
        Graduate School of Human 
        and Environmental Studies${}^{\ddag}$,\\
        Kyoto University  \\
        Kyoto 606-8501, Japan  \\
        shimizu.senjo.5s@kyoto-u.ac.jp
      }
}
\end{center}        
\begin{abstract}
End-point maximal $L^1$-regularity for parabolic initial-boundary 
value problems is considered.  For the inhomogeneous 
Dirichlet and Neumann data, maximal $L^1$-regularity 
for initial-boundary value problems is established in time 
end-point case upon the homogeneous Besov space 
$\dB_{p,1}^s(\re^n_+)$ with $1< p< \infty$ and
$-1+1/p<s\le 0$ as well as  
optimal trace estimates.  
The main estimates obtained here are sharp in the sense 
of trace estimates and it is not available by known theory 
on the class of UMD Banach spaces.  We utilize a method of 
harmonic analysis, in particular, the almost orthogonal 
properties between the boundary potentials of the Dirichlet 
and the Neumann boundary data and the Littlewood-Paley dyadic 
decomposition of unity in the Besov and the Lizorkin-Triebel  
spaces.
\end{abstract}

\section{Introduction and Main Results}
In this paper, we are concerned with maximal $L^1$-regularity 
for initial-boundary value problems 
of parabolic equations in the half-space $\re^n_+$
with inhomogeneous boundary data. \par
 
Let $X$ be a Banach space and $A$ be a closed linear 
operator in $X$ with a densely defined domain $\mathcal D(A)$.
For an initial data $u_0\in X$ and an external force 
$f \in L^{\r}(0, T; X)$ $(1\le \r\le  \infty)$, 
let $u$ be a solution to the abstract Cauchy problem:
\begin{equation}
\left\{
\spl{
    &\frac{d}{dt}u +A u =f,\ \ t>0,\\
    & \qquad u(0)=u_0.
 }
\right.
\label{1-1}
\end{equation}
Then  $A$ has maximal $L^{\r}$-regularity if 
there exists a unique solution  $u$ of \eqref{1-1} such that  
$\dsp \frac{d}{dt}u$, $Au \in  L^\r(0, T; X)$  satisfy
the estimate:
\begin{equation} \notag
 \Big\|\frac{d}{dt} u\Big\|_{L^{\r}(0, T; X)} 
  +\|A u\|_{L^{\r}(0, T; X)}
   \le  C\Big(\|u_0\|_{(X,\D(A))_{1-\frac{1}{\r},\r}} 
              +\|f\|_{L^{\r}(0, T; X)}\Big),
\label{1-2}
\end{equation}
under the restriction $u_0\in (X,\D(A))_{1-\frac{1}{\r},\r}$, 
where $(X,\D(A))_{1-\frac{1}{\r},\r}$ denotes the real interpolation
space between $X$ and $\D(A)$, and $C$ is a positive constant 
independent of $u_0$ and $f$.  Maximal regularity for parabolic 
equations was first considered by 
Ladyzhenskaya--Solonnikov--Ural'tseva \cite{LSU}. 
Then the research of maximal regularity has developed immensely 
in these last few decades by many authors;  
\cite{BCP}, \cite{CP}, \cite{CL}, \cite{dPG}, 
\cite{deS}-\cite{Du}, 
\cite{GS}, \cite{HP97}, 
\cite{KL}, \cite{McY}, \cite{Sob}, \cite{Sol64}. 
In the general framework on Banach spaces $X$ that satisfy 
{\it the unconditional martingale differences} (called as UMD), 
the general theory of maximal regularity was well-established,  
especially by Amann \cite{Am95}, \cite{Am19}, 
Denk--Hieber--Pr\"uss \cite{DHP03}, \cite{DHP07}, 
Weis \cite{Ws} (see also \cite{KW}, \cite{KuW}, \cite{PrSi16}).

On the other hand, maximal regularity on non-UMD Banach spaces, 
for instance non-reflexive Banach space such as $L^1$ or 
$L^{\infty}$ requires independent arguments.  For example, one can 
observe some results for maximal $L^1$-regularity for the Cauchy 
problem on the homogenous Banach spaces in  Danchin \cite{Dc03}, 
\cite{Dc07}, Giga--Saal \cite{GS11},  Iwabuchi \cite{Iw15}, 
Ogawa--Shimizu \cite{OgSs10}, \cite{OgSs16} 
in various non-UMD settings. In general, maximal time $L^1$-regularity 
fails over the Lebesgue spaces in spacial variables and we need to 
introduce  restrictive function classes such as homogeneous 
and inhomogeneous Besov spaces for the spacial variables.
Hence maximal $L^1$-regularity for  the initial boundary value problems 
is not well established in the general framework.

\subsection{The Dirichlet boundary condition} 
We first recall non-endpoint maximal regularity for the initial boundary 
value problems to parabolic equations.
Let $I=(0,T)$ with $0<T\le \infty$. 
Let $u$ be a solution of the initial-boundary value problem of the 
second-order parabolic equation with variable coefficients 
and the inhomogeneous Dirichlet boundary condition in the half-space 
$\re^n_+=\{ x=(x',x_n); x'\in \re^{n-1}, x_n>0\}$: 
\eq{ \label{eqn;parabolic-D}
  \left\{
  \begin{aligned}
    &\pt_t u -\sum_{1\le i,j\le n}a_{ij}(t,x)\pt_i\pt_j u=f, 
     &\qquad &t\in I,\ &&x\in \re^n_+,\\
    &\quad  u(t,x',x_n)\big|_{x_n=0}=g(t,x'), 
     &\qquad &t\in I,\ &&x'\in \re^{n-1}, \\
    &\quad u(t,x)\big|_{t=0}=u_0(x),
      &\qquad &\quad\ &&x\in \re^n_+,  
    \end{aligned}
  \right.
}
where $\pt_t$ and $\pt_i\equiv \pt_{x_i}$ are partial derivatives with 
respect to $t$ and $x_i$,  $u=u(t,x)$ denotes the unknown function,  
$u_0=u_0(x)$, $f=f(t,x)$ and $g=g(t,x')$ 
are given initial, external force and boundary data, respectively. 
The coefficient matrix $\{a_{ij}(t,x)\}_{1\le i,j\le n}$ satisfies 
uniformly elliptic condition and have enough regularity.  Namely 
$\{a_{ij}\}$ is a real valued symmetric matrix such that 
for some constant $c>0$,
$\sum_{1\le ij\le n}a_{ij}(t,x)\xi_i\xi_j\ge c|\xi|^2$ 
for all $\xi\in\re^n$ with sufficient regularity in both $t$ and $x$.

For $1\le \r\le \infty$ 
and a Banach space $X$, we denote the Lebesgue-Bochner space 
$L^{\r}(I;X)$  and  the inhomogeneous and homogeneous 
Sobolev-Bochner spaces  as $W^{1,\r}(I;X)$, $\dot{W}^{1,\r}(I;X)$, 
respectively. 
We denote a set of all continuous bounded $X$-valued 
functions over an interval $I$ by $C_b(I;X)$ and 
on $\Om\subset \re^n$ by $C_b(\Om)$. 
We  also denote a set of all  $X$-valued 
functions which is  bounded uniformly continuous 
over an interval $I$ by $BUC(I;X)$ and 
on $\Om\subset \re^n$ by $BUC(\Om)$.

In this context, the following results were obtained by 
Weidemaier \cite{Wd05} and 
Denk--Hieber--Pr\"uss \cite{DHP07}. 

\begin{prop}[The Dirichlet boundary condition \cite{DHP07}, \cite{Wd05}]
\label{prop;D-H-P} 
Let $1<\r,p<\infty$ with $1-1/(2p)\ne1/\rho$, $I=(0,T)$ with $T<\infty$. 
Then the problem \eqref{eqn;parabolic-D} admits a unique solution 
$
   u\in W^{1,\r}\big(I;L^p(\re^n_+)\big)
   \cap L^{\r}\big(I;W^{2,p}(\re^n_+)\big)
$
if and only if 
\eq{ \label{eqn;DHP-conditions}
 \spl{
    &f\in L^{\r}\big(I;L^p(\re^n_+)\big),  
    \quad
    u_0\in  B^{2(1-1/\r)}_{p,\r}(\re^n_+), \\
    &g\in  F^{1-1/2p}_{\r,p}\big(I;L^p(\re^{n-1})\big)
     \cap L^{\r}\big(I; B^{2-1/p}_{p,p}(\re^{n-1})\big), \\
    &\text{if $1-{1}/(2p)>{1}/{\r}$, then }
     u_0(x',x_n)|_{x_n=0}=g(t,x')|_{t=0}. 
}}
Besides there exists a constant $C_T>0$ depending on $n$, $p$, 
$\r$, $\{a_{ij}\}$, $T$ such that  the solution $u$ is subject 
to the inequality: 
\vspace{-2mm}
  \eqn{
  \begin{aligned}
    &\|\pt_t u\|_{L^{\r}(I;L^p(\re^{n}_+))}
    +\|\N^2 u\|_{L^{\r}(I;L^p(\re^{n}_+))} \\
    \le& C_T\Big(
        \|u_0\|_{B^{2(1-1/2p)}_{p,\r}(\re^{n}_+)}
        +\|f\|_{L^{\r}(I;L^p(\re^{n}_+))}
        +\|g\|_{F^{1-1/2p}_{\r,p}(I;L^p(\re^{n-1}))}
        +\|g\|_{L^{\r}(I;B^{2-1/p}_{p,p}(\re^{n-1}))}
         \Big), 
    \end{aligned}
  } 
where  $|\N^2u|=(\sum_{1\le i,j\le n}|\pt_i\pt_j u|^2)^{1/2}$,
$L^{\r}(I;X)$ denotes the $\r$-th powered Lebesgue-Bochner space 
upon a Banach space $X$
and $B^{2-1/p}_{p,p}(\re^{n-1})$ and $F^{1-1/2p}_{\r,p}(I;X)$
denote the interpolation spaces of the Besov and the 
Lizorkin-Triebel  type, respectively.
\end{prop}

Weidemaier \cite{Wd94} first obtained a trace theorem for 
functions in anisotropic Sobolev spaces. 
Then he extended his result to a boundary trace of 
a solution of parabolic equations in the Bochner space 
and obtained the optimal trace estimates  (\cite{Wd02}-\cite{Wd05-2}) 
with introducing the Lizorkin-Triebel  space in the time variable. 
In the proof of the results, he employed a solution 
formula with respect to the time variable
and the proof is involved the maximal function 
for a test function. 
The results in \cite{Wd88}-\cite{Wd05-2} are obtained under 
the restriction of exponents of the Bochner spaces to space $p$ 
and time $\r$ variables as  $\frac{3}{2}< p\le \r< \infty$ of 
for the Dirichlet problem and $3<p\le  \r <\infty$ for the Neumann problem, 
respectively.
Denk--Hieber--Pr\"uss \cite{DHP07} obtained the above 
necessary and sufficient 
condition of unique existence of solutions to initial-boundary 
value problems including higher order parabolic operators subject to 
general boundary conditions in a domain $\Om$ in $\re^n$ with 
a compact boundary.     
Their ingenious idea is regarding the problem as an evolution equation 
with respect to the spatial variable $x_n$ and the boundary condition 
as an initial data.
However the proof in \cite{DHP07} is based on 
the vector valued version of Mikhlin's Fourier multiplier 
theorem, and accordingly the result is restricted 
in the cases $1<\r<\infty$. Their result is essentially 
in a time local estimate because the elliptic 
operator considered there has strictly positive spectrum 
and hence the boundary conditions is limited in
the inhomogeneous real interpolation spaces.
\par

In this paper, we show {\it time global maximal $L^1$-regularity} 
for initial-boundary value problems of the evolution equation 
of parabolic type in both the inhomogeneous Dirichlet and the 
Neumann boundary conditions. 
For $0$-Dirichlet boundary data, Danchin--Mucha \cite{DM09} 
obtained maximal $L^1$-regularity for the initial-boundary value 
problem \eqref{eqn;parabolic-D} in the half-space with 
$a_{ij}(t,x)=\del_{ij}$ and $g(t,x')\equiv 0$.
On the other hand,  maximal  $L^1$-regularity for  the case of 
non-zero boundary condition requires very different treatment 
and it is far from obvious since the time exponent is the end-point 
and it is neither clear nor straightforward if a natural extension 
from the known result Proposition \ref{prop;D-H-P} holds in a 
natural exponent.  Here we consider inhomogeneous initial-boundary 
value problems in both the Drichlet and the Neumann boundary conditions 
and show that {\it a natural extension from Proposition 
\ref{prop;D-H-P} does not hold in general} (namely the conditions 
for the boundary data in \eqref{eqn;DHP-conditions} with $\r=1$)  
and we explicitly prove it by showing the end-point time exponent 
invites the end-point interpolation exponent 
similar to the case of the initial value problem. 
\vskip1mm 
Before stating our results, we define Besov spaces and Lizorkin-Triebel  
spaces in the half-space and the half-line.
Since the global estimate requires the base space for spatial
variable $x$ in the homogeneous Besov space,  
we introduce the homogeneous Besov space over $\re^n_+$ 
(see for details Lizorkin \cite{Lz}, 
Peetre \cite{Pr75}, \cite{Pr76}, Triebel \cite{Tr73}, \cite{Tr83}).   

\vskip1mm
\noindent
{\it Definition} (The Besov and the Lizorkin-Triebel  spaces).
Let $s\in \re$,  $1\le p, \sg\le \infty$. 
Let $\{\phi_j\}_{j\in \Z}$ be the Littlewood-Paley dyadic decomposition 
of unity for $x\in \re^n$, 
namely $\widehat \phi$ is the Fourier transform of a smooth radial 
function $\phi$ with 
$\widehat \phi(\xi)\ge 0$ and ${\rm supp}\,
  \widehat\phi\subset \{\xi\in\re^n\mid 2^{-1}<|\xi|<2\}$,  
and
\vspace{-2mm}
\eq{\label{eqn;LP-decomp}
 \spl{  
    &\widehat\phi_j(\xi)=\widehat\phi(2^{-j}\xi),
     \quad \sum_{j\in\Z}\widehat\phi_j(\xi)=1
     \quad \text{for\ any}\ \xi\in\re^n\setminus\{0\}, 
     \quad j\in \Z  
   \\
    &\text{ and }\quad 
    \widehat{\phi}_{\widehat{0}}(\xi)
    +\sum_{j\ge 1}\widehat \phi_j(\xi)=1\quad 
\text{for\ any}\ \xi\in\re^n,
}}
where $\widehat{\phi}_{\widehat{0}}(\xi)
\equiv \widehat{\zeta}(|\xi|)$ with a low frequency cut-off 
$\widehat{\zeta}(r)=1$ for 
$0\le r<1$ and $\widehat{\zeta}(r)=0$ for $2<r$. 
 For $s\in\re$ and $1\le p,\sigma\le \infty$,
 $\dot{B}^s_{p,\sg}(\re^n)$ be the homogeneous Besov space with norm
\eqn{
 \|\tilde{f}\|_{\dot{B}^s_{p,\sg}}
  \equiv 
  \left\{
  \begin{aligned} 
   & \Bigl(\sum_{j\in \Z}2^{s\sg j}
          \|\phi_j*\tilde{f}\|_p^{\sg}
     \Bigr)^{1/\sg},
   & 1\le \sigma<\infty,\\
   &\, \sup_{j\in \Z} 2^{s j}\|\phi_j*\tilde{f}\|_p^{\sg}, 
   & \phantom{1\le }\sigma=\infty,
  \end{aligned} 
  \right.
}
where $f*g$ denotes the convolution between Schwartz class functions 
$f$ and $g\in \S(\re^n)$ given by 
\eqn{
 f*g(x)=c_n^{-1}\int_{\re^n}f(x-y)g(y)dy,
 \qquad c_n=(2\pi)^{-\frac{n}{2}}
}
and for $f$, $g\in \S'$ as the distribution sense,
where $\S^*$ is the tempered distributions.
${B}^s_{p,\sg}(\re^n)$ be the inhomogeneous Besov space 
with norm
$$
 \|\tilde{f}\|_{B^s_{p,\sg}}
  \equiv 
  \begin{cases} \dsp
    \Bigl(\|\phi_{\hat{0}}*\tilde{f}\|_p
    +\sum_{j\in \Z}2^{s\sg j}\|\phi_j*\tilde{f}\|_p^{\sg}\Bigr)^{1/\sg},
         & 1\le \sigma<\infty,    \\
    \dsp
    \|\phi_{\hat{0}}*\tilde{f}\|_p
      +\sup_{j\in \Z} 2^{s j}\|\phi_j*\tilde{f}\|_p^{\sg}, 
         & \phantom{1\le }\sigma=\infty.
  \end{cases} 
$$
For $s\in\re$,  $1\le p<\infty$ and
 $1\le \sigma\le  \infty$, $\dot{F}^s_{p,\sg}(\re^n)$ 
be the homogeneous Lizorkin-Triebel  space with norm
\eqn{
 \|\tilde{f}\|_{\dot{F}^s_{p,\sg}}
  \equiv 
  \left\{
  \begin{aligned} 
   & \Big\|\Bigl(\sum_{j\in \Z}2^{s\sg j}
          |\phi_j*\tilde{f}(\cdot)|^{\sg}
     \Bigr)^{1/\sg}\Big\|_p,
   & 1\le \sigma<\infty,\\
   &\, \Big\|\sup_{j\in \Z} 2^{s j}
         |\phi_j*\tilde{f}(\cdot)|\Big\|_p, 
   & \sigma=\infty.
  \end{aligned} 
  \right.
}

We define the homogeneous Besov space $\dot{B}^s_{p,\sg}(\re^n_+)$ 
as the set of all measurable functions $f$ in $\re^n_+$ satisfying 
\begin{align*}
  \|f\|_{\dot{B}^s_{p,\sg}(\re^n_+)}
   \equiv \inf \Bigg\{ \|\tilde f\|_{\dot{B}^s_{p,\sg}(\re^n)}<\infty;\  
       &\tilde f=\begin{cases}f(x',x_n) &(x_n>0)\\
                             \text{any extension}& (x_n<0)
                \end{cases}\bigg\},\\
              &\qquad\qquad\tilde{f}=\sum_{j\in\Z}\phi_j*\tilde{f} 
              \text{ in } \mathcal{S}'\bigg\}.
   \eqntag \label{eqn;Besov-halfspace}
\end{align*}
Analogously we define the inhomogeneous Besov space ${B}^s_{p,\sg}(\re^n_+)$ 
as the set of all measurable functions $f$ in $\re^n_+$ satisfying 
$$
  \|f\|_{{B}^s_{p,\sg}(\re^n_+)}
   \equiv \inf \left.\left\{ \|\tilde f\|_{{B}^s_{p,\sg}(\re^n)}<\infty;\   
       \tilde f=\begin{cases}f(x',x_n) &(x_n>0)\\
                             \text{any extension}& (x_n<0)
                \end{cases}
              \right\} \right\}.
$$ 

\noindent
{\it Definition} (The Bochner-Lizorkin-Triebel spaces).
Let $s\in \re$,  $1\le p, \sg\le \infty$ and $X(\re^n_+)$ be a Banach 
space on $\re^n_+$ with the norm $\|\cdot \|_{X}$.
Let $\{\psi_k\}_{k\in \Z}$ be the Littlewood-Paley dyadic decomposition 
of unity for $t\in \re$. 
For $s\in\re$ and $1\le p< \infty$, $\dot{F}^s_{p,\sg}(\re;X)$ be 
the Bochner-Lizorkin-Triebel  space with norm
\eqn{
 \|\tilde{f}\|_{\dot{F}^s_{p,\sg}(\re;X)}
  \equiv 
  \left\{
  \begin{aligned} 
   & \Big\|\Big(\sum_{k\in \Z}2^{s\sg k}
          \|\psi_k*\tilde{f}(t,\cdot)\|_X^{\sg}
     \Bigr)^{1/\sg}\Big\|_{L^p(\re_t)},
   & 1\le \sigma<\infty,\\
   &\,\Big\| \sup_{k\in \Z} 2^{s k}\|\psi_j*\tilde{f}(t,\cdot)\|_X
      \Big\|_{L^p(\re_t)}, 
   & \phantom{1\le }\sigma=\infty.
  \end{aligned} 
  \right.
}

Analogously above, we define the Bochner-Lizorkin-Triebel  spaces 
$\dot{F}^s_{p,\sg}(I;X)$ 
as the set of all measurable functions $f$ on $X$ satisfying 
$$
  \|f\|_{\dot{F}^s_{p,\sg}(I;X)}
   \equiv \inf 
       \left.\left\{
           \|\tilde f\|_{\dot{F}^s_{p,\sg}(\re;X)}<\infty;\ 
       \, \tilde f=\begin{cases}f(t,x) &(t\in I)\\
                             \text{any extension}& (t\in\re\setminus I)
                \end{cases}
              \right\} \right\}.
$$

\noindent
We note that all the spaces of homogeneous type are understood as 
the Banach spaces by introducing the quotient spaces identifying 
all polynomial differences.
\par
\medskip

We assume that the real valued coefficients 
$\{a_{ij}(t,x)\}_{1\le i,j\le n}$ satisfy the following conditions. 
\begin{assump}\label{assump;coefficients} For $1\le i,j\le n$ and 
$(t,x)\in I\times \re^n_+$,
\begin{enumerate}
\item 
$a_{ij}(t,x)=\del_{ij}+\widetilde{a_{ij}}+b_{ij}(t,x)$,
where $\del_{ij}$ and
$\widetilde{a_{ij}}$ denote the Kronecker delta and  
components of a positive constant matrix, respectively,
\item 
$\widetilde{a_{ij}}=\widetilde{a_{ji}}$ and 
$b_{ij}(t,x)=b_{ji}(t,x)$  for all $t>0$ and $x\in \re^n_+$,
\item 
there exists a constant $c>0$
such that for any $\xi\in \re^n$, 
$$
\sum_{1\le i,j\le n}{a_{ij}}(t,x)\xi_i\xi_j\ge c|\xi|^2, 
 \quad (t,x)\in I\times \re_+^n, 
$$
\item  
$b_{ij}\in BUC\big(\re_+;\dot{B}^{\frac{n}{q}}_{q,1}(\re_+^n)\big)$ 
for some $1\le q<\infty$, where $BUC(I;X)$ denotes a set of all 
bounded uniformly continuous functions on $I$.
\end{enumerate}
\end{assump}

It is known that for $1\le q<\infty$, 
$\dot{B}^{\frac{n}{q}}_{q,1}(\re^n)$ satisfies 
${{\mathcal S}_0 (\re^n)}\hookrightarrow%
\dot{B}^{\frac{n}{q}}_{q,1}(\re^n)\hookrightarrow C_v(\re^n)$, 
where {$\mathcal{S}_0$}  denotes the  rapidly 
decreasing smooth functions {with vanishing at the origin of its 
Fourier transform  and} 
 $C_v$ denotes all continuous functions with 
vanishing at infinity, respectively 
{\rm (cf. \cite[Proposition 2.3]{OgSs16})}. 
It also satisfies the product inequality, namely 
the following inequality holds:
$$
\|b_{ij} h\|_{\dot{B}^s_{p,1}}
  \le C \|b_{ij}\|_{\dot{B}^{\frac{n}{q}}_{q,1}}
      \|h\|_{\dot B^s_{p,1}},  
  \quad 1\le p\le \infty,\,  -\frac{n}{q}<s\le 0
$$
{\rm(cf. Adibi-Paicu \cite{AP}, \cite{BCD})}. 
 Hence to ensure the uniformly elliptic condition
 (3) in Assumption 1, we split the coefficients $a_{ij}(t,x)$
 into the constant parts $\del_{ij}+\widetilde{a_{ij}}$ and the 
 decreasing functions $b_{ij}(t,x)$ as $|x|\to\infty$. 

Then our main result for the end-point case of maximal regularity to
the problem \eqref{eqn;parabolic-D} now read as the following:

\begin{thm}[The Dirichlet boundary condition] 
\label{thm;MaxReg-Dirichlet} Let $1< p<\infty$, $-1+1/p<s\le 0$  
and assume that the coefficients $\{a_{ij}\}_{1\le i,j\le n}$ satisfy 
Assumption \ref{assump;coefficients}.\par
\noindent
{\rm(1)} Suppose that $b_{ij}(t,x)\equiv 0$ for all $1\le i,j\le n$.  
Then the problem \eqref{eqn;parabolic-D} admits a unique solution
\algn{
   &u\in \dot{W}^{1,1}\big(\re_+;\dot{B}^s_{p,1}(\re^n_+)\big), \quad
   {\Del u\in L^1\big(\re_+;\dot{B}^{s}_{p,1}(\re^n_+)\big)},
}
if and only if the external, initial and boundary data in 
\eqref{eqn;parabolic-D} satisfy 
\begin{align}
    &f\in L^1\big(\re_+;\dot{B}^s_{p,1}(\re^n_+)\big), 
     \quad u_0\in \dot{B}^{s}_{p,1}(\re^n_+),
      \label{eqn;f_cond_besov}
      \\
    &g\in \dot{F}^{1-1/2p}_{1,1}\big(\re_+;\dB^s_{p,1}(\re^{n-1})\big)
         \cap  L^1\big(\re_+; \dB^{s+2-1/p}_{p,1}(\re^{n-1})\big),
     \label{eqn;g_cond_besov}
\end{align}
respectively.
Besides the solution $u$ satisfies the following estimate 
for some constant $C_M>0$ depending only on $p$, $s$ and $n$  
  \eqn{
  \begin{aligned}
    &\|\pt_t u\|_{L^1(\re_+;\dot{B}^s_{p,1}(\re^{n}_+))}
    +\|\N^2 u\|_{L^1(\re_+;\dot{B}^s_{p,1}(\re^{n}_+))} \\
    \le& C_M\big(
        \|u_0\|_{\dot{B}^s_{p,1}(\re^{n}_+)}
        +\|f\|_{L^1(\re_+;\dot{B}^s_{p,1}(\re^{n}_+))}
        +\|g\|_{\dot{F}^{1-1/2p}_{1,1}(\re_+; \dot{B}^s_{p,1}(\re^{n-1}))}
        +\|g\|_{L^{1}(\re_+; {\dot{B}^{s+2-1/p}_{p,1}(\re^{n-1})})}
         \big).
   \end{aligned}
  }
\par
\noindent
{\rm(2)} Let $1\le q<\infty$, $-n/q<s\le 0$
and let $I=(0,T)$ for $T<\infty$. Then the problem \eqref{eqn;parabolic-D} 
admits a unique solution 
\algn{
  &u\in \dot{W}^{1,1}\big(I;\dot{B}^s_{p,1}(\re^n_+)\big),
  \quad
  {\Del u\in  L^1\big(I;\dot{B}^{s}_{p,1}(\re^n_+)\big),
  }
}
if and only if  the external force, the initial data and 
the boundary data satisfy \eqref{eqn;f_cond_besov} and 
\eqref{eqn;g_cond_besov} with replacing $\re_+$ into $I$.  
Besides the solution $u$ satisfies the following estimate 
for some constant $C_M=C_M(n,p,q,\{\widetilde{a_{ij}}\})>0$
\eq{ \label{eqn;maximal-regularity-dirichlet-2}
  \begin{aligned}
     \|\pt_t u\|_{L^1(I;\dot{B}^s_{p,1}(\re^{n}_+))}
     &+\|\N^2  u\|_{L^1(I;\dot{B}^s_{p,1}(\re^{n}_+))} 
    \le C_M 
        \int_0^T  e^{\mu (T-s)}\|f(s)\|_{\dot{B}^s_{p,1}} ds\\
        &+C_M \Bigl\{1+  \max_{1\le i,j\le n}
                   \|b_{ij}\|_{L^{\infty}(I;\dot{B}^{n/q}_{q,1}) } 
                   \big(e^{\mu T}-1\big)
            \Bigr\}\\
       &\quad\times \bigl(\|u_0\|_{\dot{B}^{s}_{p,1}} 
          +\|g\|_{\dot{F}^{1-1/2p}_{1,1}(I; \dB^s_{p,1}(\re^{n-1}))}
          +\|g\|_{L^1(I;\dot{B}^{s+2-1/p}_{p,1}(\re^{n-1}))}
         \bigr),
   \end{aligned}
  } 
where $\mu=C_M^2\log(1+C_M)$.\end{thm}

\noindent 
{\bf Remarks.}
(i)\ The solution in Theorem \ref{thm;MaxReg-Dirichlet} satisfies
$$
u\in  C_b\big([0,T); \dot B^s_{p,1}(\re^n_+)\big)
$$
with $T=\infty$ for the solution $u$ in (1) and with $T<\infty$ in (2). 
The linear evolution generated by the elliptic operator 
generates $C_0$-semigroup in $ \dot B^s_{p,1}(\re^n_+)$ and
the estimate of maximal $L^1$-regularity ensures that 
the absolute continuity in $t$ to the solution. 
\par
(ii)\ Since 
$
 1-1/(2p)<1
$ 
for all  $1< p<\infty$, the compatibility condition 
\eq{  \label{eqn;compatibility_besov}
u_0(x',x_n)|_{x_n=0}=g(t,x')|_{t=0}
}
is not necessarily required for \eqref{eqn;maximal-regularity-dirichlet-2}.
\par
(iii) If $p=\infty$, the corresponding result holds for the
homogeneous Besov space 
$
 \dcB^s_{\infty,1}(\re^n)
   \equiv \overline{C_{ 00}^{\infty}(\re^n)}^{\dB^s_{\infty,1}(\re^n)},
$
{
where $C_{00}^{\infty}(\re^n)$ denotes all compactly supported $C^{\infty}$-functions with zero at the origin of its 
Fourier transform}
and
\begin{align*}
  \|f\|_{\dcB^s_{\infty,1}(\re^n_+)}
   \equiv \inf \Bigg\{ \|\widetilde f\|_{\dcB^s_{p,1}(\re^n)}<\infty;\  
       &\tilde f=\begin{cases}f(x',x_n) &(x_n>0)\\
                             \text{any extension}& (x_n\le 0)
                \end{cases}\bigg\},\\
              &\qquad\qquad\tilde{f}=
              \sum_{j\in\Z}\phi_j*\tilde{f} 
              \text{ in } \mathcal{S}'\bigg\},
\end{align*}
instead of the Besov space $\dB^s_{\infty,1}(\re^n_+)$ 
with imposing the compatibility condition
\eqref{eqn;compatibility_besov} if  $(s,p)=(0,\infty)$.
Note that $\dcB^0_{\infty,1}(\re^n_+)\subset C_v(\re^n_+)$ 
for the endpoint case $(s,p)=(0,\infty)$.  
\par

We only show the estimate for $\re_+$ in time but a  
similar estimate for the finite time interval $I=(0,T)$ with $T<\infty$ 
is also available. In such a case, the restriction on the initial data 
$u_0$ can be relaxed into the class of inhomogeneous Besov spaces 
$B^s_{p,1}(\re^n_+)\supset \dot{B}^s_{p,1}(\re^n_+)$ and the 
constant appeared in the estimate can be estimated as 
$C_M\simeq O(\log T)$ ($T\to \infty$) 
 (see \cite{OgSs16}). 
\par

The function class for the $x$-variable in Theorem 
\ref{thm;MaxReg-Dirichlet} is restricted in 
$\dB^s_{p,1}(\re^n_+)\subsetneq \dot{W}^{s,p}(\re^n_+)$ and this 
restriction is necessary for obtaining maximal $L^1$-regularity (see \cite{OgSs16}).
Connecting this fact, the conditions on the coefficients $a_{ij}(t,x)$
are so far best available.  In general, the elliptic and parabolic type 
estimates in $L^p$-setting allow us to treat much general coefficients 
such as $a_{ij}(t,x)\in VMO(I\times\re^n)$ for the whole space case, 
where $VMO$ stands for the vanishing mean oscillation (see e.g. 
Krylov \cite{Ky}).  However since maximal regularity in $L^1$
generally fails for the Lebesgue spaces $L^p$ and we need to restrict 
the spatial function class to the Besov space $\dot{B}^0_{p,1}\subsetneq L^p$ 
even for the whole space $\re^n$ case.
Such a restriction limits the condition on the coefficients as given 
in Assumption \ref{assump;coefficients}.

On the other hand, Theorem \ref{thm;MaxReg-Dirichlet} does not cover 
the end-point spatial exponent $p=1$ nor $p=\infty$.
Instead of the above result, we show a substituting estimate holds 
in $L^p(\re^n_+)$.

\begin{thm}[The Dirichlet boundary condition] 
\label{thm;MaxLikeReg-Dirichlet}  
Assume that the coefficients 
$\{a_{ij}\}_{1\le i,j\le n}$ be constants that 
satisfy Assumption \ref{assump;coefficients}, i.e., 
$b_{ij}(t,x)\equiv 0$.
\par
\noindent
{\rm(1)} Let $1\le p< \infty$. 
If the external force, the initial data and the boundary data satisfy
\begin{align}
    &f\in L^1\big(\re_+;\dot{B}^0_{p,1}(\re^n_+)\big),
     \quad 
     u_0\in {\dot{B}^{0}_{p,1}}(\re^n_+), 
     \label{eqn;f_cond_Lp}
      \\
    &g\in  \dot{F}_{1,1}^{1-1/2p}(\re_+;L^p(\re^{n-1}))
          \cap L^1(\re_+; B^{2-1/p}_{p,1}(\re^{n-1})),
     \notag
\end{align}
then there exists a unique solution $u$ to \eqref{eqn;parabolic-D} in 
$$
   W^{1,1}(\re_+;L^p(\re^n_+))
  \cap L^1(\re_+;\dot{W}^{2,p}(\re^n_+))
$$
and which satisfies the following estimate:
\eq{\label{eqn;maximal-Lp-D}
  \begin{aligned}
    &\|\pt_t u\|_{L^1(\re_+;L^p(\re^{n}_+))}
    +\|\N^2  u\|_{L^1(\re_+;L^p(\re^{n}_+))} \\
    \le& C\big(
        \|u_0\|_{\dot{B}^0_{p,1}(\re^{n}_+)}
        +\|f\|_{L^1(\re_+;\dot{B}^0_{p,1}(\re^{n}_+))}
        +\|g\|_{\dot{F}^{1-1/2p}_{1,1}(\re_+; L^p(\re^{n-1}))}
        +\|g\|_{L^1(\re_+; B^{2-1/p}_{p,1}(\re^{n-1}))}
         \big),
  \end{aligned}
  }
where $C$ is depending only on $p$ and $n$. 
\par\noindent
{\rm(2)} For $p=\infty$, 
if the external force, the initial and 
the boundary data satisfy
\begin{align*}
    &
    f\in L^1\big(\re_+;\dcB^0_{\infty,1}(\re^n_+)\big),  
    \quad 
    u_0\in {\dcB^{0}_{\infty,1}}(\re^n_+),    
      \\
    &g\in  \dot{F}_{1,1}^{1}(\re_+;  C_v(\re^{n-1}))
          \cap L^1(\re_+; \B^{2}_{\infty,1}(\re^{n-1})),
     \notag
\end{align*}
with imposing the compatibility condition \eqref{eqn;compatibility_besov},
then there exists a unique solution $u$ to \eqref{eqn;parabolic-D} in 
$$
   W^{1,1}(\re_+;L^{\infty}(\re^n_+))
  \cap L^1(\re_+;\dot{W}^{2,\infty}(\re^n_+))
$$
and the corresponding  estimate to 
\eqref{eqn;maximal-Lp-D} holds as 
\eqn{
  \begin{aligned}
    &\|\pt_t u\|_{L^1(\re_+;L^{\infty}(\re^{n}_+))}
    +\|\N^2  u\|_{L^1(\re_+;L^{\infty}(\re^{n}_+))} \\
    \le& C\big(
        \|u_0\|_{\dcB^0_{\infty,1}(\re^{n}_+)}
        +\|f\|_{L^1(\re_+;\dcB^0_{\infty,1}(\re^{n}_+))}
        +\|g\|_{\dot{F}^{1}_{1,1}(\re_+;  C_v(\re^{n-1}))}
        +\|g\|_{L^1(\re_+; \dcB^{2}_{\infty,1}(\re^{n-1}))}
         \big).
  \end{aligned}
  }
\end{thm}

\par
 The main difference between Theorem \ref{thm;MaxReg-Dirichlet}
 and \ref{thm;MaxLikeReg-Dirichlet} is not only on the condition of the 
 end-point exponents $p=1$, $\infty$ but also the difference of 
 boundary regularity.  The estimate in Theorem \ref{thm;MaxLikeReg-Dirichlet}
 only requires the Lebesgue and the inhomogeneous Besov regularity 
 in the spatial direction and this is closer result to the known 
 estimate in \eqref{eqn;DHP-conditions} due to Denk--Hieber--Pr{\"u}ss \cite{DHP07}.
 Note that the inhomogeneous Besov space $B^0_{p,1}(\re^n_+)$ is slightly 
 wider than the one for homogeneous space $\dB^0_{p,1}(\re^n_+)$.
 However it does not stand for the strict sense of maximal regularity since 
 the regularity for the initial and external force is more than the 
 regularity for the solution itself by 
 $\dB^0_{p,1}(\re^n_+)\subsetneq L^p(\re^n_+)$ for all $1\le p\le \infty$.
 \vskip1mm
\noindent
{\bf Remarks.}
(i)\ In the second statement of Theorem \ref{thm;MaxLikeReg-Dirichlet}, 
$\dot F^{1}_{1,1}(\re_+;C_v(\re^{n-1}))$ is embedded into  the continuous 
functions $C_b([0,\infty); C_v(\re^{n-1}))$ and the compatibility condition 
\eqref{eqn;compatibility_besov} is required. 
\par
(ii)\ Shibata--Shimizu \cite{SbSz07}, \cite{SbSz08} showed a boundary 
estimate by extending the boundary  data $g$ into 
$\widetilde{g}:\re_+\times \re^n_+\to \re$  
and assume that the extended function satisfies
$$ 
  \widetilde{g}
  \in  W^{1,\r}(\re_+;L^p(\re^{n}_+))
  \cap L^{\r}(\re_+;W^{2,p}(\re^{n}_+))
$$
for $1<\r,p<\infty$.   However the estimate there does not include 
the endpoint exponent $\r=1$.  
\par
(iii)\ In the case $1<\r\le \infty$, one can 
extend Proposition \ref{prop;D-H-P} into a time global version. 
Indeed, the case $p=1$ for (1) in Theorem \ref{thm;MaxLikeReg-Dirichlet} 
is the corresponding estimate to Proposition \ref{prop;D-H-P}.

\subsection{The Neumann boundary condition}
Similar to the initial-boundary value problem with the Dirichlet 
condition, we consider the initial-boundary value problem of the 
Neumann boundary condition. We assume Assumption 
\ref{assump;coefficients} for the coefficients 
$\{a_{ij}\}_{1\le i,j\le n}$.  
\eq{ \label{eqn;parabolic-N}
   \left\{
  \begin{aligned} 
    &\pt_t u -\sum_{1\le i,j\le n}a_{ij}(t,x) \pt_i\pt_j u=f, 
     &\qquad t\in I,\ \  &x\in \re^n_+,\\
   &\quad \left.\pt_n u(t,x',x_n)\right|_{x_n=0}=g(t,x'), 
     &\qquad t\in I,\ \  &x'\in \re^{n-1}, \\
   &\quad u(t,x)\big|_{t=0}=u_0(x),
      &\qquad  &x\in \re^n_+,
   \end{aligned}     
    \right.
}
where $x=(x',x_n)\in \re^n_+$ and $\pt_n$ denotes the 
normal derivative $\pt/\pt x_n$ at any boundary point of $\re^n_+$.
Denk-Hieber-Pr\"uss \cite{DHP07} showed the following:

\begin{prop}[The Neumann boundary condition \cite{DHP07}]
\label{prop;neumann-DHP}
Let  $1<\r,p<\infty$ with  $1/2-1/(2p)\ne1/\rho$, $T<\infty$ and set $I=(0,T)$ for $T<\infty$.  
The initial-boundary value problem  \eqref{eqn;parabolic-N} has 
 a unique solution $u$ in
$W^{1,\r}(I;L^p(\re^n_+))\cap L^{\r}(I;W^{2,p}(\re^n_+))$
if and only if 
\begin{align*}
    &f\in L^{\r}(I;L^p(\re^n_+)), 
    \quad u_0\in B^{2(1-1/\r)}_{p,\r}(\re^n_+),\\
    &g\in F^{1/2-1/2p}_{\r,p}(I;L^p(\re^{n-1}))
          \cap L^{\r}(I;B^{1-1/p}_{p,p}(\re^{n-1})),
\end{align*}
if $1/2-{1}/(2p)>{1}/{\r}$, 
then assume further that the compatibility condition 
$$
(\pt_n u_0)(x',x_n)|_{x_n=0}=g(t,x')|_{t=0}.
$$
\end{prop}
\par
\indent
For the case of Neumann boundary problem \eqref{eqn;parabolic-N}, 
we obtain end-point maximal $L^1$-maximal regularity as follows:
\begin{thm}[The Neumann boundary condition]\label{thm;MaxReg-Neumann} 
Let $1< p< \infty$, $-1+1/p<s\le 0$ and 
assume that the coefficients $\{a_{ij}\}_{1\le i,j\le n}$ satisfy 
Assumption \ref{assump;coefficients}.\par
\noindent
{\rm(1)}  Suppose that $b_{ij}(t,x)\equiv 0$ for all $1\le i,j\le n$.  
Then the problem \eqref{eqn;parabolic-N} admits a unique solution 
$$
   u\in \dot{W}^{1,1}(\re_+;\dot{B}^s_{p,1}(\re^n_+)),\quad
   {
   \Del u\in L^1(\re_+;\dot{B}^{s}_{p,1}(\re^n_+)),
   }
$$
if and only if the external, initial and boundary data in 
\eqref{eqn;parabolic-N} satisfy
\begin{align}
    &f\in L^1(\re_+;\dot{B}^s_{p,1}(\re^n_+)), 
    \quad u_0\in \dB^{s}_{p,1}(\re^n_+),\label{Nbd-f}\\
    &g\in \dot{F}^{1/2-1/2p}_{1,1}(\re_+;\dB^s_{p,1}(\re^{n-1}))
     \cap L^1(\re_+;\dB^{s+1-1/p}_{p,1}(\re^{n-1})),
     \label{Nbd-g}
\end{align}  
respectively.  Moreover  end-point maximal $L^1$-regularity holds:  
\eqn{
  \begin{aligned}
   & \|\pt_t u\|_{L^1(\re_+;\dot{B}^s_{p,1}(\re^{n}_+))}
    +\|\N^2 u \|_{L^1(\re_+;\dot{B}^s_{p,1}(\re^{n}_+))} \\
    \le& C\big(
        \|u_0\|_{\dot{B}^s_{p,1}(\re^{n}_+)}
        +\|f\|_{L^1(\re_+;\dot{B}^s_{p,1}(\re^{n}_+))}
        +\|g\|_{\dot{F}^{1/2-1/2p}_{1,1}(\re_+;\dB^s_{p,1}(\re^{n-1}))}
        +\|g\|_{L^{1}(\re_+;\dB^{s+1-1/p}_{p,1}(\re^{n-1}))}
         \big), 
    \end{aligned}
  }
where $C$ is depending only on $p$, $s$ and $n$. 
\par
\noindent
{\rm(2)} Let $1\le q<\infty$ and $-n/q<s\le 0$. 
For any  $T<\infty$, let $I=(0,T)$. Then there exists 
a unique solution $u$ to \eqref{eqn;parabolic-N} 
$$
  u\in \dot{W}^{1,1}\big(I;\dot{B}^s_{p,1}(\re^n_+)\big),
  \quad
  {
  \Del u\in  L^1\big(I;\dot{B}^{s}_{p,1}(\re^n_+)\big),
  }
$$
if and only if  the external force, the initial data and 
the boundary data satisfy \eqref{Nbd-f} and \eqref{Nbd-g} with 
replacing $\re_+$ into $I$. Besides the solution $u$ satisfies 
the following estimate for some constant 
$C_M=C_M(n,p,q, \{ \widetilde{a_{ij}}\})>0$
\eq{ \label{eqn;maximal-regularity-neumann-2}
  \begin{aligned}
     \|\pt_t u\|_{L^1(I;\dot{B}^s_{p,1}(\re^{n}_+))}
     &+\|\N^2  u\|_{L^1(I;\dot{B}^s_{p,1}(\re^{n}_+))} 
    \le C_M 
        \int_0^T  e^{\mu (T-s)}\|f(s)\|_{\dot{B}^s_{p,1}} ds\\
        &+C_M \bigl\{1+  \max_{1\le i,j\le n}
                   \|b_{ij}\|_{L^{\infty}(I;\dot{B}^{n/q}_{q,1}) } 
                   \big(e^{\mu T}-1\big)
            \bigr\}\\
       &\quad\times \bigl(\|u_0\|_{\dot{B}^{s}_{p,1}} 
          +\|g\|_{\dot{F}^{1/2-1/2p}_{1,1}(I; \dB^s_{p,1}(\re^{n-1}))}
          +\|g\|_{L^1(I;\dot{B}^{s+1-1/p}_{p,1}(\re^{n-1}))}
         \bigr),
   \end{aligned}
  } 
where $\mu=C_M^2\log(1+C_M)$.  
\end{thm}

\par
\noindent
For the case of the Neumann boundary condition, we have 
an analogous compensated end-point result as in 
Theorem \ref{thm;MaxLikeReg-Dirichlet} with $L^p$ type boundary data. 
\vskip2mm
\begin{thm}[The Neumann boundary condition] \label{thm;MaxLikeReg-Neumann} 
Assume that the coefficients 
$\{a_{ij}\}_{1\le i,j\le n}$ satisfy Assumption \ref{assump;coefficients} 
and be constants,  i.e., 
$b_{ij}(t,x)\equiv 0$.
\par
\noindent 
{\rm (1)} Let $1\le p< \infty$. Assume that the external, initial and boundary data 
satisfy
\begin{align}
    &f\in L^1\big(\re_+;\dot{B}^0_{p,1}(\re^n_+)\big),  
    \quad u_0\in \dot{B}^{0}_{p,1}(\re^n_+),
    \label{eqn;external_force_cond-N}
      \\
    &g\in  \dot{F}_{1,1}^{1/2-1/2p}(\re_+;L^p(\re^{n-1}))
  \cap L^1(\re_+; B^{1-1/p}_{p,1}(\re^{n-1})),
    \label{eqn;boundary_data_cond-sharp-N}
\end{align}
then there exists a unique solution $u$ to  
\eqref{eqn;parabolic-N} in 
$$
  \dot{W}^{1,1}(\re_+;L^p(\re^n_+))
   \cap L^1(\re_+;\dot{W}^{2,p}(\re^n_+))
$$
and which satisfies the following estimate:
\algn{
    &\|\pt_t u\|_{L^1(\re_+;L^p(\re^{n}_+))}
    +\|\N^2  u\|_{L^1(\re_+;L^p(\re^{n}_+))} \\
    \le& C\big(
        \|u_0\|_{\dot{B}^0_{p,1}(\re^{n}_+)}
        +\|f\|_{L^1(\re_+;\dot{B}^0_{p,1}(\re^{n}_+))}
        +\|g\|_{\dot{F}^{1/2-1/2p}_{1,1}(\re_+; L^p(\re^{n-1}))}
        +\|g\|_{L^1(\re_+;B^{1-1/p}_{p,1}(\re^{n-1}))}
         \big),
   }
where $C$ is depending only on $p$ and $n$.  
\par\noindent
{\rm(2)}
For $p=\infty$ the corresponding result to (1) holds,  i.e.,
\begin{align*}
    &f\in L^1\big(\re_+;\dcB^0_{\infty,1}(\re^n_+)\big),  
    \quad u_0\in \dcB^{0}_{\infty,1}(\re^n_+),
      \\
    &g\in  \dot{F}_{1,1}^{1/2}(\re_+; C_v(\re^{n-1}))
  \cap L^1(\re_+; \dcB^{1}_{\infty,1}(\re^{n-1})),
\end{align*}
then there exists a unique solution $u$ to  
\eqref{eqn;parabolic-N} in 
$$
  \dot{W}^{1,1}(\re_+;L^{\infty}(\re^n_+))
   \cap L^1(\re_+;\dot{W}^{2,\infty}(\re^n_+))
$$
and which satisfies the following estimate:
\algn{
    &\|\pt_t u\|_{L^1(\re_+;L^{\infty}(\re^{n}_+))}
    +\|\N^2  u\|_{L^1(\re_+;L^{\infty}(\re^{n}_+))} \\
    \le& C\big(
        \|u_0\|_{\dcB^0_{\infty,1}(\re^{n}_+)}
        +\|f\|_{L^1(\re_+;\dcB^0_{\infty,1}(\re^{n}_+))}
        +\|g\|_{\dot{F}^{1/2}_{1,1}(\re_+;  C_v(\re^{n-1}))}
        +\|g\|_{L^1(\re_+;\B^{1}_{\infty,1}(\re^{n-1}))}
         \big),
   }
where $C$ is depending only on $n$.
\end{thm}

\noindent
{\bf Remark.} In the case $p=\infty$,
the compatibility condition \eqref{eqn;compatibility_besov} 
appeared in 
Theorem \ref{thm;MaxLikeReg-Dirichlet} (2) is redundant in 
Theorem \ref{thm;MaxLikeReg-Neumann} (2) since the regularity 
for the boundary  data \eqref{eqn;boundary_data_cond-sharp-N} 
is weaker than the case in Theorem  \ref{thm;MaxLikeReg-Dirichlet} 
and the boundary data is not the continuous function in $t$-variable.

\vskip1mm 
The rest of this paper is organized as follows. We present the 
basic formulation for the proof in particular the reduction 
to the boundary value problems of the heat equations in the next 
section.
We construct explicit solution formulas of the fundamental solutions
in Section 3.   Section 4 is devoted to prove  
for the Dirichlet condition case and  Section 5 for the Neumann 
boundary condition case.   In Section \ref{sec;almostorthogonal}, 
we devote the proof of 
the key estimate almost orthogonality 
\eqref{eqn;crucial-potential-orthogonarity-org}.  
Finally we show the optimality of the main result by showing 
the sharp boundary trace estimate in Section 7.

\vskip1mm
Throughout this paper we use the following notations. 
Let $\re_+=(0,\infty)$ and 
$\re^n_+$ denote the $n$-dimensional  Euclidean half-space; 
$\{(x',x_n);x'\in \re^{n-1}, x_n\in \re_+\}$. 
For $x\in \re^n$, $\<x\>\equiv (1+|x|^2)^{1/2}$.  
The Fourier and the inverse Fourier transforms are defined 
for any rapidly decreasing funciton $f\in \S(\re^n)$ 
with $c_n=(2\pi)^{-n/2}$ by
$$
 \widehat{f}(\xi)=\F[f](\xi)
  \equiv c_n\int_{\re^n} e^{-ix\cdot \xi} f(x)dx,
 \quad
 \F^{-1}[f](x)\equiv c_n\int_{\re^n} e^{ix\cdot \xi} f(\xi)d\xi.
$$
For $f\in \mathcal{S}'(\re^{n})$, we also denote 
$\widehat{f}(\xi')=\F[f](\xi')$.
For any functions 
$f=f(t,x',x_n)$ and $g=g(t,x',x_n)$, $f\underset{(t)}{*}g $, 
$f\underset{(t,x')}{*}g $ and $f\underset{(x_n)}{*}g$ stand for 
the convolution between $f$ and $g$ with respect to the variable 
indicated under $*$, respectively.  In the summation 
$\sum_{k\in \Z}$, the parameter $k$ runs for all integers $k\in\Z$ 
and for $\sum_{k\le j}$, $k$ runs for all integers less than or equal to  
$j\in \Z$.   
  In the norm for the Bochner 
spaces on $\dF^s_{p,\r}\big(I;X(\re^{n-1})\big)$ we use 
$$
 \|f\|_{\dF^s_{p,\r}(I;X)}
 =\|f\|_{\dF^s_{p,\r}(I;X(\re^{n-1}))}
$$
unless it may cause any confusion.
For $a\in \re^n$, we denote $B_R(a)$ as the open ball centered at $a$ 
with its radius $R>0$.  We also denote the compliment of $B_R(0)$ 
by $B_R^c$.  Various constants are simply denoted by $C$ unless 
otherwise stated.

%
%
\sect{Reduction to the heat equation and outline of the proofs}
The outline of the proof of Theorems \ref{thm;MaxReg-Dirichlet} 
and \ref{thm;MaxLikeReg-Dirichlet} is summarized as follows: 
We decompose the initial-boundary value problem 
\eqref{eqn;parabolic-D} into the following three problems 
and reduce the problem into the inhomogeneous boundary 
value problem with $0$ initial and external force:
\begin{align}
  &\left\{
   \begin{aligned}
     &\pt_t u_1 -\Del u_1=0, 
         &\qquad t\in I, x\in \re^n,\\
     &\quad u_1(t,x)\big|_{t=0}= 
      \begin{cases}  u_0(x',x_n), &x_n>0\\
                    -u_0(x',-x_n), &x_n\le 0,
      \end{cases} 
        &\qquad  x\in \re^n ,   \\
   \end{aligned} \right.\label{eqn;parabolic-D0}   \\
  &\left\{
   \begin{aligned}
    &\pt_t u_2 
       -\Del u_2=0,\quad
        & t\in I,\ x\in \re^n_+,\quad\\
    &\quad \left.u_2(t,x',x_n)\right|_{x_n=0}
        =g(t,x')- u_1(t,x',x_n)\big|_{x_n=0}, 
        \qquad\qquad\qquad\quad
    &  t\in I,\ x'\in \re^{n-1}, \\
    &\quad u_2(t,x)\big|_{t=0}=0,\quad
      & \ x\in \re^n_+,\quad\\
   \end{aligned} 
   \right. \label{eqn;parabolic-D1} 
 \\
   &\left\{ 
   \begin{aligned}
    &\pt_t u_3
      -\sum_{i,j}a_{ij}(t,x)\pt_i\pt_j u_3
       =f+\sum_{i,j}\big( \widetilde{a_{ij}}+b_{ij}(t,x)\big)\pt_i\pt_j(u_1+u_2)
      \equiv F, 
      &&t\in I,\ x\in \re^n_+,\quad\\
    &\quad \left.u_3(t,x',x_n)\right|_{x_n=0}=0, 
      &&t\in I,\ x'\in \re^{n-1}, \\
    &\quad u_3(t,x)\big|_{t=0}=0,\quad
      &&\phantom{t\in I,\ \ }x\in \re^n_+, \quad 
   \end{aligned}
   \right.
   \label{eqn;parabolic-D-h} 
\end{align}   
where  the coefficients $\widetilde{a_{ij}}$ and 
$b_{ij}(t,x)$ are defined in Assumption \ref{assump;coefficients} (1) and 
$f=f(t,x)$, $g=g(t,x')$, $u_0(x)$ are given external force,
the Dirichlet boundary data, and initial data, respectively.

Then 
$$
 u(t,x',x_n)=u_1(t,x',x_n)\big|_{x_n>0}+u_2(t,x',x_n)+u_3(t,x',x_n)
$$
is the solution of \eqref{eqn;parabolic-D}. 
\par
The external force $F$ in \eqref{eqn;parabolic-D-h} contains 
 not only the given data $f$ but the solutions $u_1$ and $u_2$
 to \eqref{eqn;parabolic-D0} and \eqref{eqn;parabolic-D1} and  
 we need to verify the regularity of $u_1$ and $u_2$ in order 
 to solve \eqref{eqn;parabolic-D-h}.
Indeed, if $u_1$ and $u_2$ have maximal $L^1$-regularity: 
For $1<p<\infty$ and $-1+1/p<s\le 0$,
\eq{\label{eqn;regularity-u1}
 u_1, u_2\in W^{1,1}\big(\re_+;\dot{B}^s_{p,1}(\re^n_+)\big),
 \quad
 {
 \Del u_1,\Del u_2\in  L^1\big(\re_+;\dot{B}^{s}_{p,1}(\re^n_+)\big),
 } 
}
then under the assumption 
$b_{ij}\in BUC\big(\re_+;\dot{B}^{\frac{n}{q}}_{q,1}(\re^{n-1}_+)\big)$ 
($1\le q<-\frac{n}{s}$) it holds that 
\begin{align*}
 \Big\|
     \sum_{i,j} &\big(\widetilde{a_{ij}}+b_{ij}(t,x)\big)\pt_i\pt_j(u_1+u_2)
 \Big\|_{L^1(\re_+;\dot{B}^s_{p,1}(\re^n_+))} \nonumber\\
 \le & C\big(1+\sup_{t>0}\|b_{ij}(t,\cdot)
                  \|_{\dot{B}^{\frac{n}{q}}_{q,1}(\re^n_+)}
            \big)
   \bigl(\|{\Del} u_1\|_{L^1(\re_+;\dot{B}^{s}_{p,1})}
        +\|{\Del} u_2\|_{L^1(\re_+;\dot{B}^{s}_{p,1})}  \bigr). 
\end{align*}
We set $g-u_1$ as $h$ and assume that it is given.  
When $(s,p)=(0,\infty)$,  the compatibility condition 
$h(t,x')=g(t,x')-u_1(t,x',0)$  on $t=0$ is required and it coincides 
with $u_2(t,x)\big|_{t=0}=0$, namely 
$$
g(t,x')\big|_{t=0}=u_0(x',x_n)\big|_{x_n=0}. 
$$
The requirement for the compatibility condition is natural for 
maximal $L^{\r}$-regularity in the cases $1<\r<\infty$. 
However such a restriction is not adopted for maximal $L^1$-regularity 
except the endpoint case $(s,p)=(0,\infty)$.

For the problem \eqref{eqn;parabolic-D-h}, we first notice that 
assuming the regularity \eqref{eqn;regularity-u1},
the external force $F$ is in $L^1(I;\dot{B}^s_{p,1})$ 
under the regularity assumption on $b_{ij}(t,x)$ and 
$\widetilde{a_{ij}}$ being constant coefficients.
Then we extend the problem  \eqref{eqn;parabolic-D-h} into the whole 
space by an appropriate extension of data and 
coefficients and maximal regularity 
\eqref{eqn;maximal-regularity-dirichlet-2} 
follows from the estimate for the Cauchy problem in $\re^n$ shown 
in \cite{OgSs16}.  
Note that the result in \cite{OgSs16} treats only 
the case $\widetilde{a_{ij}}=0$, however the analogous result 
follows for the constant positive coefficient case.

Therefore,  our main issue is to consider the initial boundary
value problem  \eqref{eqn;parabolic-D1}.  
\eq{ \label{eqn;heat-D-0}
  \left\{
  \begin{aligned}
    &\pt_t u -\Del u=0, &&t\in I,\ x\in \re^n_+,\\
    &\quad \left.u(t,x',x_n)\right|_{x_n=0}=h(t,x'), 
    &&t\in I,\ x'\in \re^{n-1}, \\
    &\quad  \left. u(t,x)\right|_{t=0}=0,
    &&\phantom{t\in I,\ \ }x\in \re^n_+ ,  
  \end{aligned}\right.
}
where the boundary function $h(t,x')$ is given by 
the function after a proper linear transformed function of
$g(t,x')- u_1(t,x',x_n)\big|_{x_n=0}$ in \eqref{eqn;parabolic-D1}. 
Once we obtain maximal $L^1$-regularity to \eqref{eqn;heat-D-0} 
with the boundary trace, then the original problem can be reduced 
into the initial value problem, 
and it can be reduced into the Cauchy problem in the whole 
Euclidian space $\re^n$.  
Note that the solution $u_1(t,x)$ that has regularity
\eqref{eqn;regularity-u1} has the boundary trace estimate 
(see Theorem \ref{thm;sharp-boundary-trace-D} in Section 7 below) and the 
condition on the boundary data $h$ is the same as the 
condition on the original data $g$.
In what follows, we rewrite $u_1$ into $u$  and consider 
the initial-boundary value problem \eqref{eqn;heat-D-0}. 

If we obtain the following theorem, then the main result 
Theorem \ref{thm;MaxReg-Dirichlet} also follows: 

\begin{thm}[Maximal $L^1$-regularity by the Dirichlet boundary data]
\label{thm;boundary-trace-D} Let $1<p<\infty$ and $-1+1/p<s\le 0$. 
there exists a unique solution 
$$
  u\in \dot{W}^{1,1}\big(\re_+;\dB^s_{p,1}(\re^n_+)\big),\quad
  {
   \Del u\in  L^1\big(\re_+;\dB^{s}_{p,1}(\re^n_+)\big)
  }
$$
to \eqref{eqn;heat-D-0}  if and only if $h$ satisfies 
\begin{align}
    &h\in \dot{F}^{1-1/2p}_{1,1}\big(\re_+;\dB^s_{p,1}(\re^{n-1})\big)
     \cap  L^1\big(\re_+; \dB^{s+2-1/p}_{p,1}(\re^{n-1})\big).
    \label{eqn;boundary_data_cond_besov-D-thm1.8}
\end{align}
Besides the solution $u$ is subject to the estimate:
  \eqn{
  \begin{aligned}
    \|\pt_t u\|_{L^1(\re_+;\dot{B}^s_{p,1}(\re^{n}_+))}
   & +\|\N^2 u\|_{L^1(\re_+;\dot{B}^s_{p,1}(\re^{n}_+))} \\
    \le& C\big(
        \|h\|_{\dot{F}^{1-1/2p}_{1,1}(\re_+; \dot{B}^s_{p,1}(\re^{n-1}))}
        +\|h\|_{L^{1}(\re_+; \dot{B}^{s+2-1/p}_{p,1}(\re^{n-1}))}
         \big), 
   \end{aligned}
  }
where $C$ is depending only on $p$, $s$ and $n$.

When $p=\infty$, the analogous result holds under arranging 
the function classes as in  Theorem \ref{thm;MaxLikeReg-Dirichlet}  
with the compatibility condition 
\eq{\label{eqn;compatibility_besov-0}
  h(t,x')|_{t=0}=0.
} 
  
\end{thm}

\vskip3mm

Similarly Theorem  \ref{thm;MaxLikeReg-Dirichlet} can be reduced into the 
following:
\begin{thm}[$L^p$-estimate by the Dirichlet boundary data]
\label{thm;Lp-boundary-D}  
{\rm (1)} Let $1\le p< \infty$. If $h$ satisfies
\begin{align*}
    &h\in \dot{F}^{1-1/2p}_{1,1}\big(\re_+;L^{p}(\re^{n-1})\big)
          \cap  L^1\big(\re_+; B^{2-1/p}_{p,1}(\re^{n-1})\big),
\end{align*}
then there exists a unique solution $u$ 
to \eqref{eqn;heat-D-0} which satisfies the following estimate:
\begin{align*}
         \|\pt_t u\|_{L^1(\re_+;L^{p}(\re^{n}_+))}
      & +\|\N^2  u\|_{L^1(\re_+;L^{p}(\re^{n}_+))}  \\
   \le& C\big(
          \|h\|_{\dot{F}^{1-1/2p}_{1,1}(\re_+; L^{p}(\re^{n-1}))}
         +\|h\|_{L^{1}(\re_+; B^{2-1/p}_{p,1}(\re^{n-1}))}
         \big),  
\end{align*}
where $C$ is depending only on $p$ and $n$. 
\par\noindent
{\rm(2)} For $p=\infty$ the corresponding result to (1) holds 
under the analogous arrangement for the function classes 
as in Theorem \ref{thm;MaxLikeReg-Dirichlet} with imposing the compatibility condition \eqref{eqn;compatibility_besov-0}.
\end{thm}

For the case of the Neumann boundary condition, we decompose the problem 
\eqref{eqn;parabolic-N} into the following problems 
and reduce the boundary condition into the case of heat equation.
{\allowdisplaybreaks
\begin{align*}  
  &\left\{
    \begin{aligned}
     &\pt_t u_1 - \Del u_1=0, 
     &\qquad t\in I,\ x\in \re^n, \\
     &\quad u_1(t,x)\big|_{t=0}=\begin{cases} u_0(x',x_n), &x_n>0, \\
                                   u_0(x',-x_n), &x_n\le 0,
                     \end{cases} 
                             &\qquad  x\in \re^n,
   \end{aligned} 
   \right.\\ 
  &\left\{
   \begin{aligned}
     &\pt_t u_2 - \Del u_2=0, \quad 
     & t\in I,\  x\in \re^n_+,\quad\\
     &\quad\quad  \pt_n u_2 (t,x',x_n)\big|_{x_n=0}
                 =g(t,x')- \pt_n u_1(t,x',x_n)\big|_{x_n=0}, \quad
     &\qquad   t\in I,\ x'\in \re^{n-1}, \\
     &\quad\quad  u_2(t,x)\big|_{t=0}=0,\quad 
      & x\in \re^n_+,\quad
   \end{aligned} 
   \right. \\ 
  &\left\{ 
   \begin{aligned}
   &\pt_t u_3
      -\sum_{1\le i,j\le n} a_{ij}(t,x)\pt_i\pt_j u_3 
      \equiv F, 
      &\quad t\in I,\ x\in \re^n_+,\quad\\
    &\quad \quad \pt_n u_3(t,x',x_n)\big|_{x_n=0}=0, 
      &\quad  t\in I,\ x'\in \re^{n-1}, \\
    &\quad \quad u_3(t,x)\big|_{t=0}=0,
      &\quad  x\in \re^n_+,\quad  
   \end{aligned}
   \right.
\end{align*}
}
where $F$ is similarly defined as in \eqref{eqn;parabolic-D-h} and 
$f=f(t,x)$, $g=g(t,x')$ are given external data and the Neumann 
boundary condition, respectively and $u_0(x)$ is the initial data. 
Therefore the problem can be reduced by setting 
$h(t,x')=g(t,x')- \pt_n u_1(t,x',x_n)\big|_{x_n=0}$ 
into the following problem:
\eq{ \label{eqn;heat-N-0}
  \left\{
  \begin{aligned}
    &\pt_t u -\Del u=0, &\qquad t\in I,\ x\in \re^n_+,\quad\\
    &\quad \pt_n u(t,x',x_n)\big|_{x_n=0}=h(t,x'), 
    &\qquad  t\in I,\ x'\in \re^{n-1}, \\
    &\qquad u(t,x)\big|_{t=0}=0,
      &\qquad  x\in \re^n_+ .  \quad
  \end{aligned}
  \right.
}
Then the following result yields our main result for the Neumann 
problem Theorem \ref{thm;MaxReg-Neumann}. 

\begin{thm}[Maximal $L^1$-regularity by the Neumann boundary data]
\label{thm;boundary-trace-N}
Let $1<p<\infty$ and $-1+1/p<s\le 0$.
There exists a unique solution 
\begin{equation*}
  u \in \dot{W}^{1,1}(\re_+;\dot{B}^s_{p,1}(\re^{n}_+)),\quad
  {
  \Del u\in  L^1(\re_+;\dot{B}^{s}_{p,1}(\re^{n}_+))
  }
\end{equation*}
to \eqref{eqn;heat-N-0} if and only if 
\begin{equation}
    h\in \dot{F}^{1/2-1/2p}_{1,1}\big(\re_+;\dB^s_{p,1}(\re^{n-1})\big)
           \cap  L^1\big(\re_+; \dB^{s+1-1/p}_{p,1}(\re^{n-1})\big).
    \label{eqn;boundary_data_cond_besov-N}
\end{equation}
Besides it holds the estimate:
  \eqn{
  \begin{aligned}
    \|\pt_t u\|_{L^1(\re_+;\dot{B}^s_{p,1}(\re^{n}_+))}
    &+\|\N^2 u\|_{L^1(\re_+;\dot{B}^s_{p,1}(\re^{n}_+))} \\
    \le& C\big(
    \|h\|_{\dot{F}^{1/2-1/2p}_{1,1}(\re_+; \dot{B}^s_{p,1}(\re^{n-1}))}
        +\|h\|_{L^{1}(\re_+; \dB^{s+1-1/p}_{p,1}(\re^{n-1}))}
         \big),
  \end{aligned}
  }
where $C$ is depending only on $p$, $s$ and $n$.  

When $p=\infty$, the analogous result holds under 
arranging the function classes as in Theorem \ref{thm;MaxLikeReg-Neumann}.  
\end{thm}

For the $L^p$-estimate in Theorem  \ref{thm;MaxLikeReg-Neumann} 
is reduced into the following:

\begin{thm}[$L^p$-estimate by the Neumann boundary data]
\label{thm;Naumann-BC-Lp-trace}\par\noindent
{\rm (1)} Let $1\le p<\infty$. If the boundary data $h$ satisfies
\begin{align*}
    &h\in \dot{F}^{1/2-1/2p}_{1,1}\big(\re_+;L^p(\re^{n-1})\big)
        \cap  L^1\big(\re_+; B^{1-1/p}_{p,1}(\re^{n-1})\big),
\end{align*}
then there exists a unique solution $u$ to \eqref{eqn;heat-N-0} 
which fulfills the following estimate:
\begin{align*}
    \|\pt_t u\|_{L^1(\re_+;L^{p}(\re^{n}_+))}
   & +\|\N^2 u\|_{L^1(\re_+;L^{p}(\re^{n}_+))}   \\
    \le& C\big(
        \|h\|_{\dot{F}^{1/2-1/2p}_{1,1}(\re_+; \dot{B}^0_{p,1}(\re^{n-1}))}
        +\|h\|_{L^{1}(\re_+; B^{1-1/p}_{p,1}(\re^{n-1}))}
         \big),
\end{align*}
where $C$ is depending only on $p$ and $n$. 
\par\noindent
{\rm(2)} For $p=\infty$ the corresponding result to {\rm(1)} holds 
under the analogous arrangement for the function classes 
as in Theorem \ref{thm;MaxLikeReg-Neumann}.
\end{thm}

Therefore, in order to obtain main results Theorems 
\ref{thm;MaxReg-Dirichlet}--\ref{thm;MaxLikeReg-Dirichlet} and 
Theorems \ref{thm;MaxReg-Neumann}--\ref{thm;MaxLikeReg-Neumann},
it is enough to show Theorems \ref{thm;boundary-trace-D}--\ref{thm;Lp-boundary-D} 
and Theorems \ref{thm;boundary-trace-N}--\ref{thm;Naumann-BC-Lp-trace}, 
respectively.

To show Theorem \ref{thm;boundary-trace-D}, we first 
apply the Laplace transform with respect to $t$,
the partial Fourier transform with respect to $x'$ and  
we obtain the solution formula of \eqref{eqn;heat-D-0}
as
\eqn{ 
  \Del u(t,x',x_n)
  =
   \int_{\re_+}\int_{\re^n}\Psi_D(t-s,x'-y',\eta) h(s,y')dy'ds
  \Big|_{\eta=x_n}
}
by using the boundary potential term:
 \eq{ \label{eqn;Dirichlet-potential-0} 
 \Psi_D(t,x',\eta)
   = \frac{c_{n-1}}{2\pi i}
     \int_{\Gm}\int_{\re^{n-1}} 
           e^{\lam t +ix'\cdot \xi'}
           \lam e^{-\sqrt{\lam+|\xi'|^2}\eta}\,d\xi'd\lam,
 }
where $c_{n-1}=(2\pi)^{-(n-1)/2}$ and $\Gm$ is a pass 
parallel to the imaginary axis.
We  then extend the boundary data $h(t,x')$ into $t<0$ by the 
zero extension and it enable us to formulate the above 
formula by the Fourier transform.
The key idea to derive the boundary estimate is applying
an {\it almost orthogonal estimate} between the boundary 
potential $\psi_D$ in \eqref{eqn;Dirichlet-potential-0}  
and the Littlewood-Paley dyadic decompositions   
of unity in both space $x'$ and time $t$ variables;
$\{\phi_j(x')\}_{j\in \Z}$, $\{\psi_k(t)\}_{k\in \Z}$ .
From \eqref{eqn;Dirichlet-potential-0}, 
let $\eta>0$ as a spectral parameter 
(such as  $\eta\simeq 2^{-\ell}$ with $\ell\in \Z$) and 
consider 
\eqn{
 \Psi_{D}(t,x',\eta)\eta^2
 =\frac{c_{n-1}}{2\pi i}
     \int_{\Gm}\int_{\re^{n-1}} 
           e^{\lam t +ix'\cdot \xi'}
           \lam \eta^2
           e^{-\sqrt{\lam+|\xi'|^2}\eta}
           \,d\xi'd\lam
  }
and the almost orthogonality between the boundary 
potential  $\Psi_{D}(t,x',\eta)\eta^2$ and 
the Littlewood-Paley dyadic decompositions can be shown.  
Indeed, by setting 
\eq{ \label{eqn;Dirichlet-potential-1} 
  \Psi_{D,k,j}(t,x',\eta)
  \equiv 
     \int_{\re}\int_{\re^{n-1}}
     \Psi_D(t-s,x'-y',\eta)\psi_k(s)\phi_j(y')dy'ds
}
for $\eta>0$, the almost orthogonal property 
is presented in two-way 
estimates separated by the time-like estimate and space-like
regions.  If $\eta \in [2^{-\ell}, 2^{-\ell+1})$ 
for $\ell\in \Z$, then 
\eq{ \label{eqn;crucial-potential-orthogonarity-org}
  \|\Psi_{D,k,j}(t,\cdot,\eta) \eta^2\|_{L^1(\re^{n-1}_{x'})}
  \le
  \left\{
   \begin{aligned}
      &C_n 2^{k-2\ell}
        \big(1+2^{(n+2)(k-2\ell)}\big)
       e^{-2^{\frac12(k-2\ell)}}{2^k}{\<2^kt\>^{-2}},
        &{k\ge 2j},\\
      &C_n 2^{2j-2\ell}\big(1+2^{(n+2)(2j-2\ell)}\big)
       e^{-2^{j-\ell}}{2^k}{\<2^kt\>^{-2}},
       &  {k< 2j},
   \end{aligned}
  \right.
  }
where $\<t\>=(1+|t|^2)^{1/2}$  which are essentially shown 
in Lemma \ref{lem;pt-orthogonal-1} in Section 6 below.
Then our main strategy to show maximal $L^1$-regularity for the boundary 
term \eqref{eqn;heat-D-0} relies on the above estimate 
\eqref{eqn;crucial-potential-orthogonarity-org} 
and the proof can be complete after exchanging the solution potential 
$\Psi_{D}$ into the Littlewood-Paley dyadic decompositions $\phi_j$ 
and $\psi_k$. 

In order to complete such procedure, the above type estimate 
\eqref{eqn;crucial-potential-orthogonarity-org} plays a key role.
The only difference between the case of Dirichlet boundary condition 
and the Neumann boundary condition is to the regularity 
of the boundary data.  This is because the solution can be 
realized by the potential term 
 \eqn{ 
 \Psi_N(t,x',\eta)
   = -\frac{c_{n-1}}{2\pi i}\int_{\Gm}\int_{\re^{n-1}} 
           e^{\lam t +ix'\cdot \xi'}
            \frac{\lam}
                 {\sqrt{\lam+|\xi'|^2}}
                e^{-\sqrt{\lam+|\xi'|^2}\eta}\,d\xi'd\lam
 }
instead of \eqref{eqn;Dirichlet-potential-0}.

\sect{Boundary potentials}
\subsection{The Dirichlet boundary potential and the compatibility 
condition}

In this subsection, we derive the exact solution formula of 
\eqref{eqn;heat-D-0}.  Let $h=h(t,x')$ be the boundary 
data extended into $t<0$ by the zero extension.
We apply the Laplace transform  $\L$ in time and 
the Fourier transform in  $\re^{n-1}$-dimensional 
spatial variables to the equation \eqref{eqn;heat-D-0}.
Noting $\hat{u}(0,\xi',x_n)=0$, we obtain that  
\eqn{ 
  \left\{
  \begin{aligned} 
  & (\lam+|\xi'|^2-\pt_n^2)\widehat{\L  u}(\lam,\xi',x_n)=0, 
   \\
  & \widehat{\L u}(\lam,\xi',0)=\widehat{\L h}(\lam,\xi').
  \end{aligned}
  \right.
}
Then it follows that 
\eqn{
    \widehat{\L u}(\lam,\xi',x_n)=
   \widehat{\L h}(\lam,\xi')e^{-\sqrt{\lam+|\xi'|^2}x_n}.
}
Hence  the solution is expressed by 
\eq{ \label{eqn;potential-D1}
 u(t,x)=
   \frac{c_{n-1}}{2\pi i}\int_{\Gm}
      e^{\lam t} \int_{\re^{n-1}}e^{ix'\cdot \xi'}
     \widehat{\L h}(\lam,\xi')e^{-\sqrt{\lam+|\xi'|^2}x_n} d\xi'd\lam,
}
where $c_{n-1}=(2\pi)^{-\frac{n-1}{2}}$ and $\Gm$ denotes an integral 
path in holomorphic domain parallel to the imaginary 
axis ${\rm Re} \lam>0$. 
The solution \eqref{eqn;potential-D1} satisfies the heat equation
and the boundary condition $u(t,x',0)=h(t,x')$. 
From \eqref{eqn;potential-D1}, the Green function $G_D(t,x)$ of 
the initial-boundary value problem \eqref{eqn;heat-D-0} is identified 
as 
\eq{\label{eqn;Green-D}
 G_D(t,x',x_n)= \frac{c_{n-1}}{2\pi i}\int_{\Gm}\int_{\re^{n-1}}
                  e^{\lam t}  e^{ix'\cdot \xi'}
                  e^{-\sqrt{\lam+|\xi'|^2}x_n} d\xi'd\lam.
}

Introducing the Dirichlet boundary potential by
\eq{\label{eqn;potential-D5}
 \Psi_D(t,x',\eta)=
   \frac{c_{n-1}}{2\pi i}
     \int_{\Gm}\int_{\re^{n-1}} 
           e^{\lam t}e^{ix'\cdot \xi'}
                \lam 
                e^{-\sqrt{\lam+|\xi'|^2}\eta}
                d\xi'd\lam,  
}
where we set $\eta=x_n$, we decompose this boundary potential 
\eqref{eqn;potential-D5} by a combination of two families of 
the Littlewood-Paley dyadic decomposition of unity. 
If we put an extra-parameter $\eta\simeq 2^{-\ell}$ 
to $\Psi_D$ with $\ell\in \Z$ which is 
a substitution of the boundary parameter $x_n$, 
then we find that $\Psi_D\eta^2$ can be expressed by 
the time Littlewood-Paley decomposition  
$\{\psi_k(t)\}_k$ and the space Littlewood-Paley decomposition 
$\{\phi_j(x')\}_j$ and the relations between the parameters 
$\ell$ and $(k,j)$ are explicitly estimated.
Such estimates stand for {\it the almost orthogonality} between 
$\big\{\Psi_D(t,x',\eta)\eta^2\big|_{\eta=2^{-\ell}}\big\}_{\ell\in \Z}$ and 
$\{\psi_k(t),\phi_j(x') \}_{k,j\in \Z}$.
Here we notice that from 
\eqref{eqn;potential-D1}-\eqref{eqn;potential-D5}, 
the potential $\Psi_D$ represents  the solution operated by 
the Laplace operator, namely 
\eq{\label{eqn;green-fn-D}
  \Psi_D(t,x',\eta)  \equiv  \Del G_D(t,x',\eta).
}
Then it follows that the potential $\Psi_D$ represents 
\eqn{
   \Del u(t,x',\eta)  =\int_{\re^{n-1}}
      \int_{\re}   \Psi_D(t-s,x'-y',\eta) h(s,y')  ds
       dy'.
 }               

If we consider the case when $1-\frac{1}{2p}< \frac{1}{\r}$, 
then the class $F^{1-\frac{1}{2p}}_{\r,p}(\re_+;L^q(\re^{n-1}))$ 
of the boundary data with $p>1$ is embedded into the class 
$C_b(I;L^q(\re^{n-1}))$ and the data have to satisfy continuity 
at $t=0$ as in Proposition \ref{prop;D-H-P}. 
On the other hand, if we consider maximal $L^1$ regularity,  
the class $\dot{F}^{1-\frac1{2p}}_{1,1}(\re_+;\dB^s_{p,1}(\re^{n-1}))$ 
is not embedded into  $C_b(\re_+;\dB^s_{p,1}(\re^{n-1}))$  
because $1-\frac{1}{2p}<1$ under $p<\infty$, it does not necessarily require 
continuity at $t=0$. This shows that it is not necessary to require 
the compatibility condition  $g(t,x')\big|_{t=0}=u_0(x',x_n)\big|_{x_n=0}$
point-wisely when $p<\infty$. 
In order to obtain maximal $L^p$ regularity for $1<p<\infty$, 
the boundary data is extended to the  zero extension for $t<0$. 
\par
From \eqref{eqn;potential-D5}, we change the integral path into
$\Gm=\gm+\Gm_{\ep}\to \Gm_{\ep}$ with 
$$
  \Gm_{\ep}=L_{\ep}\cup C_{\ep}, 
$$
where 
\eq{\label{eqn;path}
L_{\ep}=\big\{\lam=i\t; \t\in (-\infty,\ep)\cup(\ep,\infty)\big\},\quad
C_{\ep}= \big\{ \lam=\ep e^{i\th};\ep>0, 
               \th:-\frac{\pi}{2}\to \frac{\pi}{2}\big\}.
 }
Then by setting $\lam\in \gm+\Gm_{\ep}$,
 \begin{align}
   \Psi_D(t,x',\eta)
   =&\lim_{\gm\to 0}c_{n+1}i^{-1}
      \int_{\gm+L_{\ep}}\int_{\re^{n-1}} 
                e^{\lam t +ix'\cdot \xi'}
                \lam 
                e^{-\sqrt{\lam+|\xi'|^2}\eta}
                d\xi'd\lam \nonumber\\
   & + \lim_{\gm\to0} c_{n+1}i^{-1}  \int_{\gm+C_{\ep}}\int_{\re^{n-1}} 
                e^{\lam t +ix'\cdot \xi'}
                \lam  
                e^{-\sqrt{\lam+|\xi'|^2}\eta}
                d\xi'd\lam \nonumber\\
  \equiv & I_{L_{\ep}}+I_{C_{\ep}},\label{eqn;potential-D7} 
 \end{align}
where $c_{n+1}=(2\pi)^{-(n+1)/2}$.
%
%
  \begin{center}
 \begin{picture}(300,220)(-50,0)
    \thicklines
    \put(0,120){\vector(1,0){215}}  
    \put(100,40){\vector(0,1){185}}  
    \thicklines 
    {\bl
     \multiput(130,40)(0,2){36}{\line(1,0){1.2} }
     \multiput(130,130)(0,2){42}{\line(1,0){1.2} }
     \multiput(103,40)(0,2){36}{\line(1,0){1.2} }
     \multiput(103,130)(0,2){42}{\line(1,0){1.2} }
    }
   \thicklines  
   {\bl
   \put(100,120){\oval(17,17)[r]}
   }
   {\bl
   \put(124,120){\oval(17,17)[r]}
   }
 \put(120,60){\vector(-1,0){12}}
 \put(120,88){\vector(-1,0){12}}
 \put(120,140){\vector(-1,0){12}}
 \put(120,180){\vector(-1,0){12}}
    \put(138,105){$\gm+\Gm_{\ep}$}
    \put(105,105){$\ep$}
    \put(130,200){as $\gm\to 0$}
 \large
    \put(70,210){Im}   
    \put(200,100){Re}
\end{picture}
\vskip -12mm
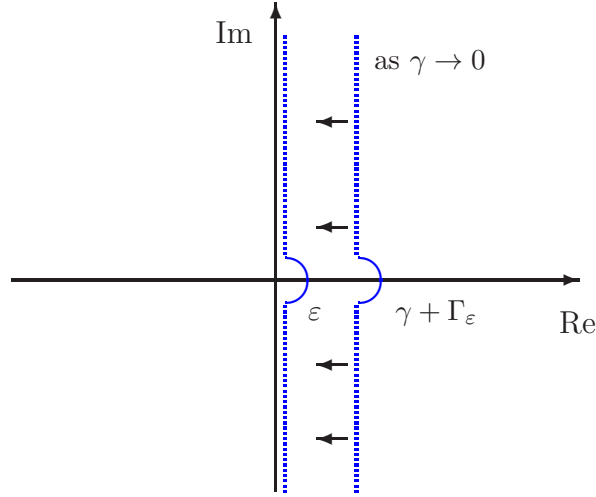
Fig. 1: The integral path $\gm+\Gm_{\ep}$ and $\Gm_{\ep}$
\end{center}
 
The first term of the right hand side of  \eqref{eqn;potential-D7}
is converging even $\gm\to 0$ 
if we observe the real part of the exponential integrant is 
\eq{\label{eqn;potential-D9}
   \big|\lam  e^{-\sqrt{i\t+|\xi'|^2}\eta}\big|
    =  |\t|  \exp\big(-(\t^2+|\xi'|^4)^{\frac14}\eta
               \cos\big(\frac{1}2\tan^{-1}\frac{\t}{|\xi'|^2} \big)
                \big)
}
 from 
 \eqn{
 \begin{aligned}
  i\t e^{-\sqrt{i\t+|\xi'|^2}\eta}
    = & i\t \exp\big(-(\t^2+|\xi'|^4)^{\frac14}\eta
             e^{\frac{i}2\tan^{-1}\frac{\t}{|\xi'|^2} }\big) \\
    = & i\t \exp\big(-(\t^2+|\xi'|^4)^{\frac14}\eta
               \big\{\cos\big(\frac{1}2\tan^{-1}\frac{\t}{|\xi'|^2} \big)
               +i\sin\big(\frac{1}2\tan^{-1}\frac{\t}{|\xi'|^2} \big) 
               \big\}\big).
  \end{aligned}
 }
 We emphasize that for 
 $-\frac{\pi}{4}<\frac12\tan^{-1}\frac{\t}{|\xi'|^2}<\frac{\pi}{4}$ 
 the exponent appeared in \eqref{eqn;potential-D9} is negative 
 definite unless   $(\t,\xi')=(0,0)$ and hence  
 $\lim_{\ep\to 0}I_{L_{\ep}}$ converges. i.e.,
 {\allowdisplaybreaks
  \begin{align}
  \lim_{\ep\to 0}I_{L_{\ep}}
    =&\lim_{\ep\to 0}\lim_{\gm\to 0}c_{n+1}i^{-1}e^{\gm t}
      \int_{L_{\ep}}\int_{\re^{n-1}} 
                e^{i\t t +ix'\cdot \xi'}
                (\gm+i\t) 
                e^{-\sqrt{\gm+i\t+|\xi'|^2}\eta}
                d\xi' i d\t
 \nonumber\\
     =&\lim_{\ep\to 0}c_{n+1}
       \int_{L_{\ep}}\int_{\re^{n-1}} 
                e^{i\t t +ix'\cdot \xi'}
                i\t 
                e^{-\sqrt{i\t+|\xi'|^2}\eta}
                d\xi'd\t
\nonumber\\
     =&c_{n+1}
       \int_{\re\setminus\{0\}}\int_{\re^{n-1}} 
                e^{i\t t +ix'\cdot \xi'}
                 i \t 
                e^{-\sqrt{i\t+|\xi'|^2}\eta}
                d\xi'd\t.\label{eqn;potential-D11}
 \end{align}
 }
 For the convergence of the second term $I_{C_{\ep}}$ of 
 the right hand side  of \eqref{eqn;potential-D7}, 
 we see  by noting 
 \eq{\label{eqn;potetial-D13}
 \left|\tan^{-1}\big(\frac{\ep \sin \th}{\ep \cos \th+|\xi'|^2}\big)\right|
 \le  \left|\tan^{-1}\big(\frac{ \sin \th}{ \cos \th}\big)\right|
  = |\th|, 
  \qquad  -\frac{\pi}{2}\le \th\le \frac{\pi}{2}
 }
 that for $\eta>0$
{\allowdisplaybreaks
  \begin{align*}
  \lim_{\ep\to 0}|I_{C_{\ep}}|
     \le & 
       \lim_{\ep\to0}  \Big|c_{n+1}i^{-1}\int_{C_{\ep}}\int_{\re^{n-1}} 
           e^{\lam t +ix'\cdot \xi'}
                 \lam 
                e^{-\sqrt{\lam+|\xi'|^2}\eta}
                d\xi'd\lam \Big| \\
     \le & \lim_{\ep\to 0} c_{n+1}
         \int_{\re^{n-1}}
          \Big|\int_{-\frac{\pi}{2}}^{\frac{\pi}{2}}           
                e^{\ep te^{i\th} }
                \ep e^{i\th}  
                e^{-\sqrt{\ep e^{i\th}+|\xi'|^2}\eta}
                \ep i e^{i\th}d\th
          \Big| d\xi'  \\
     = & \lim_{\ep\to 0} c_{n+1}
         \int_{\re^{n-1}}
            \int_{-\frac{\pi}{2}}^{\frac{\pi}{2}}           
                e^{\ep t\cos \th }
                \ep^2  \eta^2 \\
           &\hskip1.3cm
            \times 
            \exp\Big(-\eta\big((\ep \cos \th+|\xi'|^2)^2+\ep^2\sin^2 \th\big)^{1/4}
                      \cos\big(\frac12\tan^{-1}
                      \big(\frac{\ep \sin \th}{\ep \cos \th+|\xi'|^2}\big)
                \Big)
              d\th  d\xi'  \\
     \le & \lim_{\ep\to 0} c_{n+1}
         \int_{\re^{n-1}}
              \int_{-\frac{\pi}{2}}^{\frac{\pi}{2}}           
                e^{\ep t\cos \th}
                 \ep^2   
                 \exp\big(-\eta\big((\ep \cos \th+|\xi'|^2)^2\big)^{1/4}
                     \cos\big(\frac{\th}2 \big)
                     \big)
                d\th  d\xi'  \\
    \le & \lim_{\ep\to 0} c_{n+1}\ep^2  
             \int_{\re^{n-1}}  
                 \exp\Big(-\frac{\sqrt{2}}{2} \eta|\xi'|
                     \Big)d\xi' 
              \int_{-\frac{\pi}{2}}^{\frac{\pi}{2}}           
                e^{\ep t\cos \th } 
                d\th  \\
       \le & c_{n+1}\lim_{\ep\to 0} \ep^2  
              \Big(\int_{\re^{n-1}}  
                 \exp\big(-\frac{\sqrt{2}}{2} \eta|\xi'|
                     \big)d\xi'
              \Big) 
              \Big(\int_{-1}^{1}           
                \frac{e^{\ep t|\zeta|}}{\sqrt{1-\zeta^2}}             
                d\zeta
              \Big)
        =0. 
     \eqntag \label{eqn;potential-D15}
  \end{align*}
 }
Therefore, by passing $\ep\to 0$ in \eqref{eqn;potential-D7}, 
we obtain the following formula from \eqref{eqn;potential-D11} 
and \eqref{eqn;potential-D15}:
 \eq{ \label{eqn;potential-D}
 \Psi_D(t,x',\eta)
   = c_{n+1}\int_{\re}\int_{\re^{n-1}} 
           e^{it\t +ix'\cdot \xi'}
                i\t  
                e^{-\sqrt{i\t+|\xi'|^2}\eta}
                d\xi'd\t.
 }
\noindent
{\bf Remark.} We note that from \eqref{eqn;potential-D1} and similar 
way in \eqref{eqn;potential-D5}-\eqref{eqn;potential-D15}, 
the Green's function \eqref{eqn;Green-D} is expressed by the Fourier 
inverse transform as 
\eqn{
 G_D(t,x',\eta)
     = c_{n+1}\int_{\re}\int_{\re^{n-1}} 
                e^{it\t +ix'\cdot \xi'}
                e^{-\sqrt{i\t+|\xi'|^2}\eta}
                d\xi'd\t.
 }
Therefore the potential function $\Psi_D$ given by \eqref{eqn;potential-D}
is understood as the Green's function operated by the Laplacian;
$$
  \Psi_D(t,x',\eta)=\Del_{(x',\eta)} G_D(t,x',\eta).
$$ 
Note that the $0$-initial condition for $G_D(t,x',\eta)\big|_{t=0}=0$ 
is fulfilled by the Cauchy integral theorem on the same complex path 
\eqref{eqn;path} avoiding the branch cut at the negative real-line 
and passing the limit $\gm\to 0$ and $\ep\to 0$ along the analogous  
estimates \eqref{eqn;potential-D11} and \eqref{eqn;potential-D15}.

\subsection{The Neumann boundary potential and the compatibility condition}
Following the method of the case of the Dirichlet boundary condition, 
we consider the initial-boundary value problem \eqref{eqn;heat-N-0}.
Then we deduce the boundary potential (the Green's function) that 
yields the solution to \eqref{eqn;heat-N-0}.
We apply the Laplace transform in time and the Fourier transform 
in  $\re^{n-1}$-dimensional spatial variables and noting 
$\widehat{\L u}(0,\xi',x_n)=0$ to have 
\eqn{
  \left\{
  \begin{aligned}
  &(\lam+|\xi'|^2-\pt_n^2)\widehat{\L  u}=0, \\
  &\pt_n\widehat{\L u}(\lam,\xi',0)=\widehat{\L h}(\lam,\xi').
  \end{aligned}
  \right.
}
Then we write it explicitly, 
\eq{ \label{eqn;potential-N3}
 u(t,x)=
   -\frac{c_{n-1}}{2\pi i}\int_{\Gm}
      e^{\lam t} \int_{\re^{n-1}}e^{ix'\cdot \xi'}
     \frac{\widehat{\L h}(\lam,\xi')}
         {\sqrt{\lam+|\xi'|^2}}
       e^{-\sqrt{\lam+|\xi'|^2}x_n} d\xi'd\lam,
}
where $c_{n-1}=(2\pi)^{-\frac{n-1}{2}}$ and  $\Gm$ is 
a proper path on the analytic region.  
From \eqref{eqn;potential-N3}, the Green's function $G_N(t,x)$ of 
the initial-boundary value problem \eqref{eqn;heat-N-0} is given by 
\eq{\label{eqn;Green-N}
 G_N(t,x',x_n)=  -\frac{c_{n-1}}{2\pi i}\int_{\Gm}\int_{\re^{n-1}}
                  e^{\lam t}  e^{ix'\cdot \xi'}
                  \frac{1}{\sqrt{\lam+|\xi'|^2}}
                  e^{-\sqrt{\lam+|\xi'|^2}x_n} d\xi'd\lam.
}
One can choose as 
the parallel line to the imaginary axis in $Re \lam>0$. 
Let
\eq{ \label{eqn;potential-N5}
 \Psi_N(t,x',\eta)=
   -\frac{c_{n-1}}{2\pi i}
     \int_{\Gm}\int_{\re^{n-1}} 
           e^{\lam t}e^{ix'\cdot \xi'}
           \frac{\lam}{\sqrt{\lam+|\xi'|^2}}
           e^{-\sqrt{\lam+|\xi'|^2}\eta}
           d\xi'd\lam,
}
where $\Gm$ is a proper integral path basically parallel to 
the imaginary axis.    
For the case of maximal $L^1$-regularity, the 
boundary regularity is 
$\dot{F}^{\frac12-\frac{1}{2p}}_{1,p}\big(I;L^p(\re^{n-1})\big)$
and the continuity in time direction does not hold, hence 
the compatibility condition is redundant.
Along the Dirichlet boundary case before, the boundary condition can be 
prolonged to $t<0$ by  zero extension and hence the solution is 
understood by  zero extension.
From \eqref{eqn;potential-N5}, we pass the integral 
path \eqref{eqn;path} into
 $\Gm=\gm+\Gm_{\ep}\to \Gm_{\ep}$ with 
 $$
  \Gm_{\ep}=L_{\ep}\cup C_{\ep} 
 $$
by $\gm\to 0$.
 Then
 \eq{\label{eqn;potential-N7}
 \begin{aligned}
   \Psi_N(t,x',\eta)
   =&\lim_{\gm\to 0}c_{n+1}i 
      \int_{\gm+L_{\ep}}\int_{\re^{n-1}} 
                e^{\lam t +ix'\cdot \xi'}
                \frac{\lam }{\sqrt{\lam+|\xi'|^2}}
                e^{-\sqrt{\lam+|\xi'|^2}\eta}
                d\xi'd\lam \\
   & + \lim_{\gm\to0} c_{n+1}i  \int_{\gm+C_{\ep}}\int_{\re^{n-1}} 
                e^{\lam t +ix'\cdot \xi'}
                \frac{\lam}{\sqrt{\lam+|\xi'|^2}}
                e^{-\sqrt{\lam+|\xi'|^2}\eta}
                d\xi'd\lam \\
  \equiv & II_{L_{\ep}}+II_{C_{\ep}}.
 \end{aligned}
 } 
The first term of the right hand side of  \eqref{eqn;potential-N7}
is converging even $\gm\to 0$ 
if we observe the real part of the exponential integrant 
is  for the case $\gm=0$ that \
\eq{ \label{eqn;potential-N9}
  \frac{i\t }{\sqrt{i\t+|\xi'|^2}}
     =\frac{i\t 
            e^{\frac{i}{2}\tan^{-1}\frac{-\t}{|\xi'|^2}}
           }
          {(\t^2+|\xi'|^4)^{\frac14}},
 }
and \eqref{eqn;potential-D5}.  Then we obtain that 
\eq{\label{eqn;potential-N11}
   \left|\frac{i\tau }{\sqrt{i\t+|\xi'|^2}} 
      e^{-\sqrt{i\t+|\xi'|^2}\eta}\right|
    =   \frac{|\t|}
             { (\t^2+|\xi'|^4)^{\frac14} }
         \exp\big(-(\t^2+|\xi'|^4)^{\frac14}\eta
               \cos\big(\frac{1}2\tan^{-1}\frac{\t}{|\xi'|^2} \big)
              \big).
 }
 Under the condition on the argument 
 $-\frac{\pi}{4}<\frac12\tan^{-1}\frac{\t}{|\xi'|^2}<\frac{\pi}{4}$
 the real part of the exponential function in \eqref{eqn;potential-N9} 
 is negative unless $(\t,\xi)=(0,0)$.  We observe that for 
 $-\frac{\pi}{4}<\frac12\tan^{-1}\frac{\t}{|\xi'|^2}<\frac{\pi}{4}$
\eqref{eqn;potential-N11} is integrable unless   $(\t,\xi')=(0,0)$ and 
hence $II_{L_{\ep}}$ converges. i.e.,
\eq{  \label{eqn;potential-N13}
  \lim_{\ep\to 0}II_{L_{\ep}}
     =-c_{n+1}
       \int_{\re\setminus\{0\}}\int_{\re^{n-1}} 
                e^{i\t t +ix'\cdot \xi'}
                \frac{i\t}{\sqrt{i\t+|\xi'|^2}}
                e^{-\sqrt{i\t+|\xi'|^2}\eta}
                d\xi'd\t.
}  
 For the convergence of the second term $I_{C_{\ep}}$ of the right hand side 
 of \eqref{eqn;potential-N7}, we see  by noting \eqref{eqn;potetial-D13} 
 and $\eta>0$ that 
{\allowdisplaybreaks 
 \begin{align}
  \big|\lim_{\ep\to 0}I_{C_{\ep}}\big|
     \le & \lim_{\ep\to 0} 
          c_{n+1}\int_{\re^{n-1}}
          \Big|\int_{-\frac{\pi}{2}}^{\frac{\pi}{2}}           
                e^{\ep e^{i\th} }
                \frac{\ep e^{i\th} }{\sqrt{\ep e^{i\th}+|\xi'|^2}}
                      e^{-\sqrt{\ep e^{i\th}+|\xi'|^2}\eta}
                \ep i e^{i\th}d\th
          \Big| d\xi'  \nonumber\\
    \le & \lim_{\ep\to 0} c_{n+1}
         \int_{\re^{n-1}}
              \int_{-\frac{\pi}{2}}^{\frac{\pi}{2}}  
                e^{\ep \cos \th t}         
                \frac{\ep^2 }%
                {\big((\ep \cos \th+|\xi'|^2)^2+\ep^2\sin^2 \th\big)^{1/4}
                 }\nonumber\\
           &\qquad
            \times 
            \exp\Big(-\eta\big((\ep \cos \th+|\xi'|^2)^2+\ep^2\sin^2 \th\big)^{1/4}
                     \cos\big(\frac12\tan^{-1}
                     \big(\frac{\ep \sin \th}{\ep \cos \th+|\xi'|^2}\big)
                     \Big)
                d\th  d\xi'  \nonumber\\
     \le & \lim_{\ep\to 0} c_{n+1}
           \int_{\re^{n-1}}
              \int_{-\frac{\pi}{2}}^{\frac{\pi}{2}}           
                e^{\ep \cos \th t}
                   \frac{ \ep^2   }
                        {\big(\ep \cos \th+|\xi'|^2\big)^{1/2}
                         }
                 \exp\big(-\eta|\xi'|\cos\big(\frac{\th}2 \big)
                     \big)
                d\th  d\xi'  \nonumber\\
       \le & c_{n+1} \lim_{\ep\to 0} \ep^{\frac32}  
              \Big(\int_{\re^{n-1}}  
                   \exp\big(-\frac{\sqrt{2}}{2} \eta|\xi'|
                     \big)d\xi'
              \Big) 
              \Big(\int_{-1}^{1}           
                \frac{e^{\ep t|\zeta|}}{\sqrt{|\zeta|} \sqrt{1-\zeta^2}}             
                d\zeta
              \Big)=0.  \label{eqn;potential-N15}
\end{align}
}
Then it follows from \eqref{eqn;potential-N13} and \eqref{eqn;potential-N15} 
that 
\eq{ \label{eqn;potential-N}
 \Psi_N(t,x',\eta)
   = -c_{n+1}\int_{\re}\int_{\re^{n-1}} 
           e^{it\t +ix'\cdot \xi'}
           \frac{i\t}{\sqrt{i\t+|\xi'|^2}}
           e^{-\sqrt{i\t+|\xi'|^2}\eta}
           d\xi'd\t.
 }
We notice that the potential of the solution operated by 
the Laplace operator is given by 
\eq{\label{eqn;green-fn-N}
\Psi_N(t,x',\eta)\equiv \Del G_N(t,x',\eta),
}
where
\eq{\label{eqn;Green-N}
 G_N(t,x',x_n)=  -c_{n+1}\int_{\re}\int_{\re^{n-1}}
                   e^{it\tau+ix'\cdot \xi'}
                  \frac{1}{\sqrt{i\t+|\xi'|^2}}
                  e^{-\sqrt{\lam+|\xi'|^2}\eta} d\xi'd\t.
}
 We then regard this boundary potential as 
a role of Littlewood-Paley dyadic decomposition of unity
and the main argument for the proof consists on exchanging the boundary  
potential into the standard Littlewood-Paley decomposition 
$\{\psi_k\}_{k\in \Z}$ and $\{\phi_j\}_{j\in\Z}$.

%
%
\sect{The Dirichlet boundary condition}
\label{sec;4}
%
%
\subsection{The Besov spaces on the half-spaces}\par
First we recall the summary for the Besov spaces over the half-space 
on the Euclidean space $\re^n_+$.

\vskip2mm
\noindent
{\it Definition.}  Let $1\le p<\infty$ and $1\le \sg <\infty$ with $s\ge 0$.
Let 
\begin{align*}
   &\overset{\circ\quad}{B^s_{p,\sg}}(\re^n_+)
       =\overline{C_0^{\infty}(\re^n_+) }^{\dB^s_{p,\sg}(\re^n_+)}, 
\\
   &\overset{\odot}{B^s}_{p,\sg}(\re^n_+)
      =\overline{
        \{f\in \dB^s_{p,\sg}(\re^n); 
            \supp f \subset \re^n_+\}}^{\dB^s_{p,\sg}(\re^n)}, 
\end{align*}
 by the Besov norm $\dB^s_{p,\sg}(\re^n_+)$
 (see  Triebel \cite{Tr78} Section 2.9.3).
 It is shown that the above defined space coincides 
the space $\dB^s_{p,\sg}(\re^n_+)$ defined by the restriction 
in \eqref{eqn;Besov-halfspace}.
Namely, the following proposition is shown by 
Triebel \cite{Tr78} and Danchin--Mucha \cite{DM09}.

\begin{prop}[\cite{DM09}, \cite{Tr78}]
\label{prop;DM-09} Let $1<p<\infty$.\par\noindent
{\rm (1)} For $0<s$, $1\le \sg<\infty$, 
\eqn{
   \dB^{-s}_{p',\sg'}(\re^n_+)
   \simeq \big(\overset{\circ\quad}{B^s_{p,\sg}}(\re^n_+)\big)^*.
}
{\rm (2)}  For $-\infty<s\le \frac{1}{p}$ and 
for $1<\sg<\infty$,
\eqn{
 \overset{\circ\quad}{B^s_{p,\sg}}(\re^n_+)\simeq \dB^s_{p,\sg}(\re^n_+). 
}
{\rm (3)} For  $-\infty<s< \frac{1}{p}$ and  $\sg=1$,
\eqn{
 \overset{\circ\quad}{B^s_{p,1}}(\re^n_+)\simeq \dB^s_{p,1}(\re^n_+). 
}
{\rm (4)} For $-1+\frac{1}{p}<s<\frac{1}{p}$ and $1\le \sg <\infty$,
\eqn{
  \overset{\odot\quad}{B^s_{p,\sg}}(\re^n_+) 
  \simeq \dB^s_{p,\sg}(\re^n_+).
  }
\end{prop}

\vskip2mm
We consider the restriction operator $R_0$ 
by multiplying a cut-off function 
\eqn{
 \chi_{\re^n_+}(x)=
  \begin{cases}
   1,    & \text{ in } \re^n_+, \\
   0,    & \text{ in } \overline{\re^n_-}.
  \end{cases}
}
i.e., for $f\in \dB^s_{p,\sg}(\re^n)$, 
set $R_0f=\chi_{\re^n_+}(x)f(x)$ in $\dB^s_{p,\sg}(\re^n)$
if $s>0$ and it is understood as a distribution sense.
Let the extension operator $E_0$ from $\overset{\odot\quad}{B^s_{p,\sg}}(\re^n_+)$ given by the zero-extention, i.e.,
 for any $f\in \overset{\odot\quad}{B^s_{p,\sg}}(\re^n_+)$, set 
$$
 E_0f=
 \begin{cases}
   f(x), & \text{ in } \re^n_+, \\
   0,    & \text{ in } \overline{\re^n_-}.
 \end{cases}
$$

One can find that those operators are basic tool to 
recognize the homogeneous Besov spaces. Using 
{Proposition \ref{prop;DM-09}},
the following statement is a variant introduced by 
Triebel \cite[p.228]{Tr78}{.}  
\begin{prop}\label{prop;retraction-coretraction}
Let $1\le p<\infty$, $1\le \sg<\infty$ and $-1+\frac{1}{p}<s<\frac{1}{p}$.
It holds that  
\begin{align}
    &R_0:\dB^s_{p,\sg}(\re^n)
        \to\dB^s_{p,\sg}(\re^n_+),  \label{eqn;retraction}\\ 
    &E_0:\dB^s_{p,\sg}(\re^n_+)
        \to \dB^s_{p,\sg}(\re^n), \label{eqn;co-retraction}
\end{align}
are linear bounded operators.  Besides
it holds that
\eq{
    R_0E_0=Id: \dB^s_{p,\sg}(\re^n_+)\to \dB^s_{p,\sg}(\re^n_+),
      }
where $Id$ denotes the identity operator.  
Namely $E_0$ and $R_0$ are a retraction and a co-retraction, respectively.
\end{prop}
\vskip2mm
The proof of Proposition \ref{prop;retraction-coretraction} 
is along the same line of the proof in \cite{Tr78}.  Note that the 
spaces are homogeneous Besov spaces and then the arrangement 
appears in Proposition 3 in Danchin--Mucha \cite{DM09} is required.

\vskip2mm
\begin{prf}{Proposition \ref{prop;retraction-coretraction}}
It is clear that both 
\eq{\label{eqn;retraction-co-retraction-2}
  \spl{
    &R_0:\dB^s_{p,\sg}(\re^n)
        \to\dB^s_{p,\sg}(\re^n_+), \\ 
    &E_0:\dB^s_{p,\sg}(\re^n_+)
        \to \dB^s_{p,\sg}(\re^n)
  }     
  }
are linear operators and 
$$
R_0E_0=Id:\dB^s_{p,\sg}(\re^n_+)
        \to \dB^s_{p,\sg}(\re^n_+),
$$
then it sufficient to show that both operators are  
bounded.  Then they are retraction and co-retraction, respectively.

To see the first operator \eqref{eqn;retraction} is bounded,
let $f\in \dB^s_{p,\sg}(\re^n)$ 
and  we show that 
$$
 \|R_0 f\|_{\dB^s_{p,\sg}(\re^n_+)}
 =    \inf \|\widetilde{R_0f}\|_{\dB^s_{p,\sg}(\re^n)}
 \le \| \chi_{\re^n_+}f\|_{\dB^s_{p,\sg}(\re^n)}
 \le \| f\|_{\dB^s_{p,\sg}(\re^n)}
$$
is valid under the restriction $s>0$.
Let $-1/p'<s<0$ and {$f\in \dB^{s}_{p,\sg}(\re^n)$}.
Suppose that $\phi \in C_0^{\infty}(\re^n)$, Then 
\eq{
 \spl{
  \big|\<\chi_{\re^n_+}f,\phi\>\big|
  =  & \big|\<f, \chi_{\re^n_+}\phi\>\big|
  \le \|f\|_{\dB^{s}_{p,\sg}(\re^n)} 
      \|\chi_{\re^n_+}\phi\|_{\dB^{-s}_{p',\sg'}(\re^n)} \\
 \le &\|f\|_{\dB^{s}_{p,\sg}(\re^n)} 
      \|\phi\|_{\dB^{-s}_{p',\sg'}(\re^n)}, 
}}
since $0<-s< 1/p'$ and the last inequality follow from the pointwise sense.
Thus from the definition of the norm in $\dB^s_{p,\sg}(\re^n_+)$,
it holds similarly to the above that 
\eqn{
    \|R_0 f \|_{\dB^{s}_{p,\sg}(\re^n_+)}
   \le \| \chi_{\re^n_+} f \|_{\dB^{s}_{p,\sg}(\re^n)} 
   = \sup_{\phi\in\dB^{-s}_{p',\sg'}(\re^n)\setminus\{0\} }
    \frac{  \big|\<\chi_{\re^n_+} f,\phi\>\big|  }
         {  \|\phi\|_{\dB^{-s}_{p',\sg'}(\re^n)} }
   \le \|f\|_{\dB^{s}_{p,\sg}(\re^n)}. 
  }

To see the second bound \eqref{eqn;co-retraction},
we introduce a quotient space $\dB^s_{p,\sg}(\re^n)/\sim$,
where we identify all $f\in \dB^{s}_{p,\sg}(\re^n)$ 
that coincides on $\re^n_+$. Then  the restriction operator $R_0$ is 
one to one mapping from $\dB^s_{p,\sg}(\re^n)/\sim$ onto $\dB^s_{p,\sg}(\re^n_+)$, and then
the extension operator $E_0$ is an inverse operator of $R_0$.
Thus the open mapping theorem implies the required boundedness directly.
\end{prf}

\vskip2mm
In what follows, we restrict ourselves to the regularity 
range of the Besov spaces $\dB^s_{p,\sg}(\re^n_+)$ 
in $-1+1/p<s<1/p$ for $1<p<\infty$.
Hence we freely use the above mentioned results. As a consequence, 
any component in $\dB^s_{p,\sg}(\re^n_+)$ under such restriction 
on $s$ and $p$ can be extended into the one over whole space 
$\re^n$ and 
conversely the component $f\in \dB^{s}_{p,\sg}(\re^n)$ is 
restricted into the one over the half-space $\re^n_+$.
We frequently use those facts without noticing for every case
below.

\subsection{The L-P decomposition with separation of variables}\par
In order to split the variables $x'\in \re^{n-1}$ and $x_n\in \re_+$, 
we introduce an $x'$-parallel decomposition and an $x_n$-parallel
decomposition by Littlewood-Paley type. In what follows $\eta\in \re_+$
denotes a parameter for $x_n$-axis in $\re^n_+$. 
We introduce $\{\overline{\Phi_m}\}_{m\in \Z}$ as a Littlewood-Paley dyadic 
frequency decomposition of unity in separated variables 
$(\x',\x_n)\in \re^{n-1}\times \re$. 
%
%
\begin{center}
  \begin{picture}(300,220)(-50,0)
    
    \thicklines
    \put(-10,80){\vector(1,0){220}}  
    \put(100,50){\vector(0,1){165}}  
    \thicklines
    {
     \multiput(0,80)(0,2){50}{\line(1,0){1.2} }
     \multiput(200,84)(0,2){50}{\line(1,0){1.2} }
     \multiput(55,80)(0,2){31}{\line(1,0){1.2} }
     \multiput(145,80)(0,2){31}{\line(1,0){1.2} }
     \multiput(0,180)(2,0){100}{\line(0,1){1.2} }
     \multiput(55,140)(2,0){46}{\line(0,1){1.2} }
%
    }
    {
    \multiput(145,160)(2,0){28}{\line(0,1){1.2} }
    \multiput(145,150)(2,0){28}{\line(0,1){1.2} }
    \multiput(145,140)(2,0){28}{\line(0,1){1.2} }
    \multiput(145,130)(2,0){28}{\line(0,1){1.2} }
    \multiput(145,120)(2,0){28}{\line(0,1){1.2} }
    \multiput(145,110)(2,0){28}{\line(0,1){1.2} }
    \multiput(145,100)(2,0){28}{\line(0,1){1.2} }
    \multiput(145,90)(2,0){28}{\line(0,1){1.2} }
    \multiput(145,80)(2,0){28}{\line(0,1){1.2} }
    \multiput(1,160)(2,0){28}{\line(0,1){1.2} }
    \multiput(1,150)(2,0){28}{\line(0,1){1.2} }
    \multiput(1,140)(2,0){28}{\line(0,1){1.2} }
    \multiput(1,130)(2,0){28}{\line(0,1){1.2} }
    \multiput(1,120)(2,0){28}{\line(0,1){1.2} }
    \multiput(1,110)(2,0){28}{\line(0,1){1.2} }
    \multiput(1,100)(2,0){28}{\line(0,1){1.2} }
    \multiput(1,90)(2,0){28}{\line(0,1){1.2} }
    \multiput(1,80)(2,0){28}{\line(0,1){1.2} }
    }
    {
    \multiput(0,140)(0,2){22}{\line(1,0){1.2} }
    \multiput(10,140)(0,2){22}{\line(1,0){1.2} }
    \multiput(20,140)(0,2){22}{\line(1,0){1.2} }
    \multiput(30,140)(0,2){22}{\line(1,0){1.2} }
    \multiput(40,140)(0,2){22}{\line(1,0){1.2} }
    \multiput(50,140)(0,2){22}{\line(1,0){1.2} }
    \multiput(60,140)(0,2){22}{\line(1,0){1.2} }
    \multiput(70,140)(0,2){22}{\line(1,0){1.2} }
    \multiput(80,140)(0,2){22}{\line(1,0){1.2} }
    \multiput(90,140)(0,2){22}{\line(1,0){1.2} }
    \multiput(100,140)(0,2){22}{\line(1,0){1.2} }
    \multiput(110,140)(0,2){22}{\line(1,0){1.2} }
    \multiput(120,140)(0,2){22}{\line(1,0){1.2} }
    \multiput(130,140)(0,2){22}{\line(1,0){1.2} }
    \multiput(140,140)(0,2){22}{\line(1,0){1.2} }
    \multiput(150,140)(0,2){22}{\line(1,0){1.2} }
    \multiput(160,140)(0,2){22}{\line(1,0){1.2} }
    \multiput(170,140)(0,2){22}{\line(1,0){1.2} }
    \multiput(180,140)(0,2){22}{\line(1,0){1.2} }
    \multiput(190,140)(0,2){22}{\line(1,0){1.2} }
    }
{
    \put(75,185){$2^{m+1}$}
    \put(77,160){$2^{m}$}
    \put(77,128){$2^{m-1}$}
    \put(89,65){$0$}
    \put(175,65){$2^m$}
    \put(142,65){$2^{m-1}$}
    \put(200,65){$2^{m+1}$}}
\large 
    \put(80,207){$\xi_n$}   
    \put(235,65){$|\xi'|$}
\end{picture}\\
\vspace{-10mm}
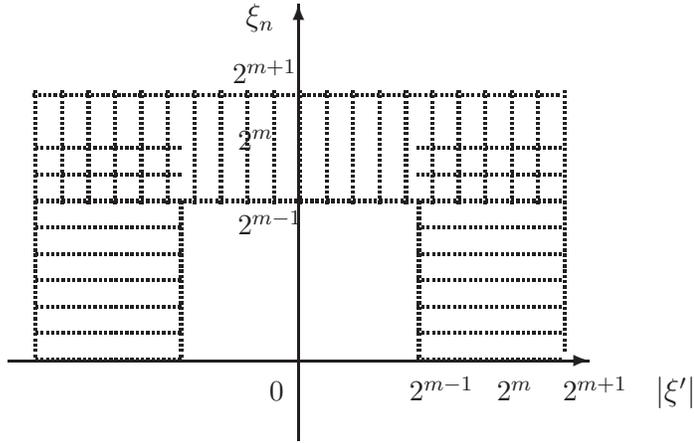
Fig 2:  The support of Littlewood-Paley decomposition $\{\overline{\Phi_m}\}_{m\in \Z}$
\end{center}

\vskip3mm
\noindent{\it Definition} (The Littlewood-Paley decomposition 
of separated variables).
For  $m\in \Z$, let
\eq{\label{eqn;zeta} 
 \begin{aligned}
 &\widehat{\zeta_m}(\xi_n)
   =\left\{
    \begin{aligned}
       1,&\qquad 0\le |\xi_n|\le 2^{m}, \\
       \text{smooth}&, 2^{m}<|\xi_n|<2^{m+1},\\
       0,&\qquad  2^{m+1}\le |\xi_n|, 
    \end{aligned}
    \right. \\
 &\widehat{\zeta_m}(\xi_n) 
    =\widehat{\zeta_{m-1}}(\xi_n)+\widehat{\phi_m}(\xi_n)
 \end{aligned}
}
(one can choose 
$\widehat{\zeta_m}(r)
 =\sum_{\ell\le m-1}\widehat{\phi_{\ell}}(r)+\widehat{\phi_{-\infty}}(r)$
with a correction distribution $\widehat{\phi_{-\infty}}(r)$ at $r=0$)
and set
\eq{\label{eqn;direct-sum-L-P}
  \widehat{\overline{\Phi_m}}(\xi)
  \equiv\widehat{\phi_m}(|\xi'|)\otimes \widehat{\zeta_{m-1}}(\xi_n)
   +\widehat{\zeta_{m}}(|\xi'|)\otimes \widehat{\phi_{m}}(\xi_n).
}
Then it is obvious from Fig. 2  (restricted on the upper half 
region in $\re^n$) that  
\eq{ \label{eqn;Sum_Phi_is_1}
 \sum_{m\in \Z}\widehat{\overline{\Phi_m}}(\xi)\equiv 1,
 \quad \xi=(\xi',\xi_n) \in  \re^n\setminus \{0\}.
}
Indeed, from \eqref{eqn;zeta} and \eqref{eqn;direct-sum-L-P},
\algn{
   &\sum_{m\in \Z}\widehat{\overline{\Phi_m}}(\xi)\\
  =&  \sum_{m\in\Z} \widehat{\phi_m}(|\xi'|)
                    \otimes \sum_{-\infty\le\ell\le m-1}\widehat{\phi_{\ell}}(\xi_n)
    + \sum_{m\in\Z} \Big(\sum_{-\infty\le\ell\le m}\widehat{\phi_{\ell}}(|\xi'|)\Big)
                    \otimes \widehat{\phi_{m}}(\xi_n) 
 \\
  =&  \sum_{m\in\Z} \widehat{\phi_m}(|\xi'|)
                    \otimes \sum_{-\infty\le\ell\le m-1}\widehat{\phi_{\ell}}(\xi_n)
    + \sum_{m\in\Z} \sum_{\ell\le m} \widehat{\phi_{\ell}}(|\xi'|)
                    \otimes\widehat{\phi_{m}}(\xi_n)
    + \sum_{m\in\Z} \widehat{\phi_{-\infty}}(|\xi'|)\otimes \widehat{\phi_{m}}(\xi_n) 
  \\
  =&  \sum_{m\in\Z} \widehat{\phi_m}(|\xi'|)
                    \otimes \sum_{-\infty\le\ell\le m-1}\widehat{\phi_{\ell}}(\xi_n)
    + \sum_{\ell\in\Z} \widehat{\phi_{\ell}}(|\xi'|)
                    \otimes \sum_{m\ge \ell}\widehat{\phi_{m}}(\xi_n)
    + \sum_{m\in\Z}  \widehat{\phi_{-\infty}}(|\xi'|)
                    \otimes \widehat{\phi_{m}}(\xi_n) 
 \\
  =&  \sum_{m\in\Z} \widehat{\phi_m}(|\xi'|)\otimes 
                   \Big(\sum_{-\infty\le\ell\le m-1}\widehat{\phi_{\ell}}(\xi_n)
                   + \sum_{\ell\ge m}\widehat{\phi_{\ell}}(\xi_n)\Big) 
      + \widehat{\phi_{-\infty}}(|\xi'|)\otimes \sum_{m\in\Z} \widehat{\phi_{m}}(\xi_n) 
 \\
  =&  \sum_{m\in\Z} \widehat{\phi_m}(|\xi'|)\otimes 
      \sum_{\ell\in\Z\cup\{-\infty\}} \widehat{\phi_{\ell}}(\xi_n)
        + \widehat{\phi_{-\infty}}(|\xi'|)
          \otimes \sum_{\ell\in\Z\cup\{-\infty\}} \widehat{\phi_{m}}(\xi_n) 
        - \widehat{\phi_{-\infty}}(|\xi'|)\otimes \widehat{\phi_{-\infty}}(\xi_n)
 \\
   =& \Big( \sum_{m\in\Z\cup\{-\infty\}} \widehat{\phi_m}(|\xi'|)\otimes 1\Big) 
      - \widehat{\phi_{-\infty}}(|\xi'|)\otimes \widehat{\phi_{-\infty}}(\xi_n) \\
   =& 1- \widehat{\phi_{-\infty}}(|\xi'|)\otimes \widehat{\phi_{-\infty}}(\xi_n).
}

\vskip2mm
\noindent
{\it Definition} (Varieties of the Littlewood-Paley dyadic decompositions).
Let $(\t,\xi',\xi_n)\in  \re \times \re^{n-1}\times \re$ be  Fourier adjoint variables corresponding to 
$(t,x',\eta)\in \re_+\times \re^{n-1}\times \re_+$.
\begin{enumerate}
\item[$\bullet$] $\{\Phi_m(x)\}_{m\in \Z}$:
      the standard (annulus type) Littlewood-Paley dyadic decomposition by \hfill\break $ x=(x',\eta)\in \re^n_+$.  
\item[$\bullet$]  $\{ \overline{\Phi_m}(x) \}_{m\in \Z}$:
      the Littlewood-Paley dyadic decomposition  over $ x=(x',\eta)\in \re^n_+$
      given by \eqref{eqn;direct-sum-L-P}.
\item[$\bullet$]  $\{\psi_k(\tilde{t})\}_{k\in \Z}$:
     the Littlewood-Paley dyadic decompositions  in  $\tilde t\in \re$.
\item[$\bullet$]  $\{\phi_j(x')\}_{j\in \Z}$ and $\{\phi_j(\tilde\eta)\}_{j\in \Z}$:
      the standard (annulus type) Littlewood-Paley dyadic decompositions in $x'\in \re^{n-1}$ 
      and  $\tilde\eta\in \re$, 
      respectively.
\item[$\bullet$]  $\{\zeta_m(x')\}_{m\in \Z}$ and $\{\zeta_m(\tilde\eta)\}_{m\in \Z}$:
      the lower frequency smooth cut-off 
      given by \eqref{eqn;zeta}, respectively.
\item[$\bullet$] Let $\widetilde{\phi_j}=\phi_{j-1}+\phi_j+\phi_{j+1}$ be the 
 Littlewood-Paley dyadic decompositions with its $j$-neighborhood.
\item[$\bullet$]  Since all the above defined 
decompositions are even functions,
 we identify $\tilde{t}\in \re$ and $\tilde\eta\in \re $ with 
 $|\tilde t|=t>0$ and $|\tilde\eta|=\eta>0$, respectively. 
\end{enumerate}
\vskip3mm

Then it is easy to see that the norm of the Besov spaces on $\re^n_+$
defined by $\{\Phi_m\}_m$ is equivalent to the one from the 
Littlewood-Paley decomposition of direct sum type, $\{\overline{\Phi_m}\}_m$:
\begin{align}
\|\Del u(t)\|_{\dB^s_{p,1}(\re^n_+)} 
    =& \sum_{m\in \Z} 2^{sm}
       \|\Phi_m\underset{(x)}{*}  \sum_{|m-k|\le 1} 
          \overline{\Phi_{k}}\underset{(x)}{*} \Del u(t)\|_{L^p(\re^n_+)} \nonumber\\
 \le & 3C\sum_{m\in \Z} 2^{sm} 
       \|\overline{\Phi_m}\underset{(x)}{*}\Del u(t)\|_{L^p(\re^n_+)}.   
\label{eqn;potential-besov-1}        
\end{align}

\subsection{Separation on the Dirichlet potential}\par
In order to show the sufficiency part of Theorem \ref{thm;boundary-trace-D},
we first decompose the solution potential $\Psi_D$ defined in 
\eqref{eqn;potential-D}.
From the solution formula with the Green function 
\eqref{eqn;green-fn-D} and \eqref{eqn;potential-D}, 
we see that 
{\allowdisplaybreaks
 \algn{  
  \overline{\Phi_m}\underset{(x',\eta)}{*}
     &\big(\Psi_D(t,x',\eta)\big) \nonumber\\
     =& \Big(\phi_m(|x'|)\otimes \zeta_{m-1}(\eta)\Big)
          \underset{(x',\eta)}{*}\Psi_D(t,x',\eta) 
       +\Big(\zeta_m(|x'|)\otimes \phi_{m}(\eta)\Big)
          \underset{(x',\eta)}{*}\Psi_D(t,x',\eta) \\
     =&  \zeta_{m-1}(\eta)\underset{(\eta)}{*}
         \Big(\phi_m(|x'|)\underset{(x')}{*}\Psi_D(t,x',\eta) \Big)
       + \phi_{m}(\eta) \underset{(\eta)}{*}
         \Big(\zeta_m(|x'|)\underset{(x')}{*}\Psi_D(t,x',\eta)\Big) \\
     =&  \zeta_{m-1}(\eta)\underset{(\eta)}{*}
                \frac{c_{n+1}}{c_{n-1}}
                \int_{\re}\int_{\re^{n-1}}
                e^{it\t +ix'\cdot \xi'}  
                \widehat{\phi_m}(|\xi'|)
                \, i\t\, e^{-\sqrt{i\t+|\xi'|^2}\eta}
                d\xi'd\t    
        \nonumber\\
    &+   \phi_{m}(\eta)\underset{(\eta)}{*}   
              \frac{c_{n+1}}{c_{n-1}}\int_{\re}\int_{\re^{n-1}}
                e^{it\t +ix'\cdot \xi'}
                \widehat{\zeta_{m}}(|\xi'|)
                \, i\t\, e^{-\sqrt{i\t+|\xi'|^2}\eta}
                d\xi'd\t.
 \eqntag\label{eqn;potential-besov-3} 
 }
 }
%
%
To estimate the solution by the Besov-norm, we involve the 
Littlewood-Paley dyadic decomposition of the direct sum type
\eqref{eqn;direct-sum-L-P}. 
 Besides $\zeta_m(\eta)*$ can be treated by  the 
Hausdorff-Young inequality of $\eta$-variable, 
i.e., the $L^p(\re^{n}_+)$ norm of the first term of 
the right hand side of \eqref{eqn;potential-besov-3}
is estimated as follows: 
\eqn{
 \spl{
 \Big\|\zeta_{m-1}\underset{(\eta)}{*} &
       \Big(\phi_m(x')\underset{(x')}{*}\Psi_D(t,x',\eta) \Big)
 \Big\|_{L^p(\re_{+,\eta};L^p(\re^{n-1}_{x'}))} \\
 \le & \|\zeta_{m-1}\|_{L^1(\re_{+,\eta})}
       \big\|\phi_m(x')\underset{(x')}{*}\Psi_D(t,x',\eta) 
       \big\|_{L^p(\re_{+,\eta};L^p(\re^{n-1}_{x'}))} \\
  =  &C\big\|\phi_m(x')\underset{(x')}{*}\Psi_D(t,x',\eta) 
       \big\|_{L^p(\re_{+,\eta};L^p(\re^{n-1}_{x'}))}.
}}
Applying the Hausdorff-Young inequality to 
 the first term of the right hand side of 
\eqref{eqn;potential-besov-3} and 
restrict the range of the summation of $j$ by 
$\zeta_m(x')\underset{(x')}{*}$, we have from
\eqref{eqn;potential-besov-1} that 
{\allowdisplaybreaks 
\algn{ 
  \int_0^{\infty} &\|\Del u(t)\|_{\dB^s_{p,1}(\re^{n}_+)}dt  
  \\
  \le &
     C\Bigl\|\sum_{m\in \Z} 2^{sm}
            \Bigl\|
             \overline{\Phi_m}
             \underset{(x',\eta)}{*}
             \Psi_D(t,x',\eta)
             \underset{(t,x')}{*}h(t,x')
            \Bigr\|_{L^p(\re^n_+)}
        \Bigr\|_{L^1_t(\re_+)}
   \\
     \le &C\Bigl\|\sum_{m\in \Z}2^{sm}
         \Bigl(\int_{\re_+} 
         \Bigl\| \Psi_D(t,x',\eta)
                \underset{(t,x')}{*} 
                \big(\phi_m(x')\underset{(x')}{*}h(t,x') \big) 
         \Bigr\|_{L^p(\re^{n-1}_{x'})}^p 
         d\eta \Bigr)^{1/p}\Bigr\|_{L^1_t(\re_+)}
 \\
        &+C\bigg\|\sum_{m\in \Z}2^{sm}
         \Bigl(\int_{\re_+}
         \Bigl\|
               \phi_{m}(\eta)\underset{(\eta)}{*}
               \Psi_D(t,x',\eta) 
               \underset{(t,x')}{*} \\
         &\hskip4cm \times 
               \sum_{k\in \Z} \sum_{j\le m+1}
               \psi_k(t)\underset{(t)}{*}\phi_j(y')
                        \underset{(y')}{*} h(t,y')
         \Big\|_{L^p(\re^{n-1}_{x'})}^p  d\eta \Bigr)^{1/p}
         \bigg\|_{L^1_t(\re_+)}
 \\
\equiv &\|P_1^D\|_{L^1_t(\re^+)}+\|P_2^D\|_{L^1_t(\re^+)},
\eqntag\label{eqn;potential-besov-9}
}
}
where the first term of the right hand side of \eqref{eqn;potential-besov-9}  
includes  $\phi_m(x')$,
once the outer decomposition $\sum_{m\in \Z}$ is fixed then 
the inner decomposition $\{\phi_j(x')*\}_{j\in \Z}$
is restricted into only $|j-m|\le 1$ and the summation for $j$ disappears.

%
%
\subsection{ Time-space splitting argument}
We separate the estimate of \eqref{eqn;potential-besov-9}
into two regions ; one is  time-dominated area  and 
the other is  space-dominated area.

\vskip4mm
\noindent
$\bullet$  The relation between each variables:

\begin{center}
 \begin{picture}(300,200)(-50,0)
    \thicklines
    \put(0,80){\vector(1,0){215}}  
    \put(10,-20){\vector(0,1){215}}  
    \thinlines
    \put(10,80){\line(1,4){4}} 
    \put(13,90){\line(1,3){4}} 
    \put(13,90){\line(1,2){4}} 
    \put(17,100){\line(1,1){18}} 
    \put(35,118){\line(4,3){20}} 
    \put(50,130){\line(2,1){40}} 
    \put(89,150){\line(3,1){40}} 
    \put(130,163){\line(4,1){50}} 
    \put(10,80){\line(1,-3){4}} 
    \put(13,70){\line(1,-2){5}} 
    \put(17,60){\line(1,-1){18}} 
    \put(35,42){\line(4,-3){20}} 
    \put(50,30){\line(2,-1){40}} 
    \put(89,10){\line(3,-1){40}} 
    \put(130,-3){\line(4,-1){50}} 
    \thicklines
    { 
     \multiput(10,72)(0,2){7}{\line(1,0){1.2} }
     \multiput(12,68)(0,2){10}{\line(1,0){1.2} }
     \multiput(15,65)(0,2){15}{\line(1,0){1.2} }
     \multiput(20,58)(0,2){24}{\line(1,0){1.2} }
     \multiput(25,50)(0,2){30}{\line(1,0){1.2} }
     \multiput(30,48)(0,2){35}{\line(1,0){1.2} }
     \multiput(40,40)(0,2){42}{\line(1,0){1.2} }
     \multiput(50,33)(0,2){49}{\line(1,0){1.2} }
     \multiput(60,30)(0,2){53}{\line(1,0){0.4} }
     \multiput(70,20)(0,2){60}{\line(1,0){1.2} }
     \multiput(80,10)(0,2){67}{\line(1,0){0.4} }
     \multiput(90,8)(0,2){71}{\line(1,0){0.4} }
     \multiput(100,5)(0,2){75}{\line(1,0){1.2} }
     \multiput(110,0)(0,2){80}{\line(1,0){0.4} }
     \multiput(120,0)(0,2){82}{\line(1,0){0.4} }
     \multiput(130,0)(0,2){83}{\line(1,0){0.4} }
     \multiput(140,-3)(0,2){84}{\line(1,0){1.2} }
    }
    {
     \multiput(8,170)(2,0){73}{\line(0,1){0.3} }
     \multiput(8,160)(2,0){55}{\line(0,1){0.3} }
     \multiput(8,150)(2,0){41}{\line(0,1){1.2} }
     \multiput(8,140)(2,0){30}{\line(0,1){0.3} }
     \multiput(8,130)(2,0){21}{\line(0,1){1.2} }
     \multiput(8,120)(2,0){14}{\line(0,1){1.2} }
     \multiput(8,110)(2,0){8}{\line(0,1){1.2} }
     \multiput(8,50)(2,0){8}{\line(0,1){1.2} }
     \multiput(8,40)(2,0){15}{\line(0,1){1.2} }
     \multiput(8,30)(2,0){20}{\line(0,1){1.2} }
     \multiput(8,20)(2,0){30}{\line(0,1){0.6} }
     \multiput(8,10)(2,0){38}{\line(0,1){1.2} }
     \multiput(8,0)(2,0){56}{\line(0,1){0.6} }
     \multiput(8,-10)(2,0){70}{\line(0,1){0.6} }
    }
   \thinlines
    \put(-5,65){$0$}
 \large
    \put(10,180){$|\x'|\simeq 2^j, j\in \Z$}   
    \put(150,60){$\t\simeq 2^k, k\in \Z$}
\end{picture}
\vskip5mm
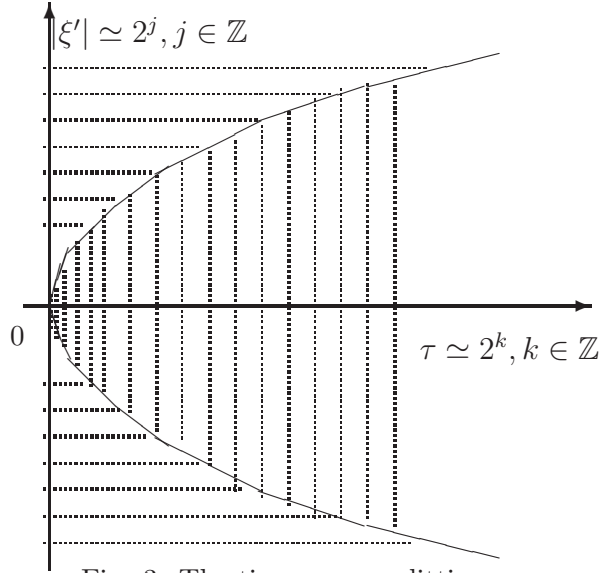
Fig. 3: The time-space splitting.
\end{center}

In order to prove Theorem \ref{thm;boundary-trace-D}, it is enough 
to prove Lemma \ref{lem;P1-P2-estimate}. 
\begin{lem}\label{lem;P1-P2-estimate}
Let $1<p\le \infty$.
The term $P_1^D$ defined in 
\eqref{eqn;potential-besov-9} is estimated as follows:
\eq{ \label{eqn;P1-est}
 \|P_1^D\|_{L^1_t(\re_+)}
 \le C\Big(
       \big\|h \big\|_{ \dF^{1-1/2p}_{1,1}(\re_+;\dB^s_{p,1}(\re^{n-1})) }
      +\big\|h \big\|_{ L^1(\re_+;\dB^{s+2-1/p}_{p,1}(\re^{n-1}))  }
      \Big).
}
Simultaneously the term $P_2^D$ defined in 
\eqref{eqn;potential-besov-9} is estimated as follows:
 \eq{ \label{eqn;P2-est}
 \|P_2^D\|_{L^1_t(\re_+)}
 \le C\Big(
       \big\|h \big\|_{ \dF^{1-1/2p}_{1,1}(\re_+;\dB^s_{p,1}(\re^{n-1}))}
      +\big\|h \big\|_{ L^1 (\re_+;\dB^{s+2-1/p}_{p,1}(\re^{n-1})) }
      \Big).
}
\end{lem}

\begin{prf}{Lemma \ref{lem;P1-P2-estimate}}
We split the data $h$ into the  time-dominated region and the 
space-dominated region. 
\begin{align}
 h(t,x')=&
  \sum_{k\in \Z}\sum_{j\in \Z}
   \psi_k(t)\underset{(t)}{*}
   \phi_j(x')\underset{(x')}{*}h(t,x') \nonumber\\
   =& \sum_{k\in\Z} \sum_{2j\le k}
      \psi_k(t)\underset{(t)}{*}
      \phi_j(x')\underset{(x')}{*}h(t,x') 
    + \sum_{k\in\Z} \sum_{2j>k}   
      \psi_k(t)\underset{(t)}{*}
      \phi_j(x')\underset{(x')}{*}h(t,x').
 \label{eqn;space-time-splitting}     
 \end{align}
Letting $h_m(t,x')\equiv  \widetilde{\phi_m}\underset{(x')}{*}h(t,x')$ 
($m\in \Z$), we proceed 
{\allowdisplaybreaks
\algn{ 
 P_1^D
 =&\; C\sum_{m\in \Z} 2^{sm}
       \bigg( \int_{\re_+}
        \big\|
         \sum_{k\in\Z} \sum_{|j-m|\le 1, 2j<k} 
         \Psi_D(t,x',\eta)
         \underset{(t,x')}{*}        
         \psi_k(t)\underset{(t)}{*}
         \phi_j(x')\underset{(x')}{*}  \\
   &\hskip7.5cm \times
         \widetilde{\phi_m}(x')\underset{(x')}{*} 
           h(t,x')  
        \big\|_{L^p(\re^{n-1}_{x'})}^p
         d\eta
      \bigg)^{1/p} \\
 &+ C\sum_{m\in \Z} 2^{sm}
        \bigg(\int_{\re_+}
         \big\| 
          \sum_{j\in \Z}\sum_{|j-m|\le 1,k<2j}
          \Psi_D(t,x',\eta)\underset{(t,x')}{*}
          \psi_k(t)\underset{(t)}{*}
          \phi_j(x')\underset{(x')}{*}
      \\
    &\hskip7.5cm \times
         \widetilde{\phi_m}(x')\underset{(x')}{*} h(t,x')  
         \big\|_{L^p(\re^{n-1}_{x'})}^p
         d\eta
        \bigg)^{1/p}\\
 \le &\; C\sum_{m\in \Z} 2^{sm}
       \bigg( \int_{\re_+}
        \big\|
         \sum_{k\ge 2m} 
         \Psi_D(t,x',\eta)\underset{(t,x')}{*}        
         \psi_k(t)\underset{(t)}{*} 
         \phi_m(x')\underset{(x')}{*}
          h_{m}(t,x')  
        \big\|_{L^p_{x'}}^p
        d\eta
      \bigg)^{1/p} \\
   &+ C\sum_{m\in \Z} 2^{sm}
        \bigg(\int_{\re_+}
         \big\| 
          \sum_{k<2m}
          \Psi_D(t,x',\eta)\underset{(t,x')}{*}
          \psi_k(t)\underset{(t)}{*}
          \phi_m(x')\underset{(x')}{*}
          h_{m}(t,x')  
         \big\|_{L^p_{x'}}^p
         d\eta
        \bigg)^{1/p}\\
  \equiv &  L_1 +  L_2 .
  \eqntag\label{eqn;potential-besov-10}
}
}
Setting 
\eq{ \label{eqn;chi}
  \Psi_{D,k,m}(t,x',\eta)
  \equiv 
     \int_{\re}\int_{\re^{n-1}}
     \Psi_D(t-s,x'-y',\eta)\psi_k(s)\phi_m(y')
      dy'ds,
}
we see that  $ L_1$ is the time-dominated region and
applying the Minkowski and  the Hausdorff-Young inequality
with using \eqref{eqn;chi}, we have 
{\allowdisplaybreaks
\begin{align*} 
  L_1
  \le &C\sum_{m\in\Z}  2^{sm}
       \Bigg(\int_{\re_+}
        \bigg\{ \sum_{k\ge 2m} 
        \int_{\re_+} 
        \bigg\|
           \int_{\re^{n-1}}              
              \Psi_{D,k,m}(t-s,x'-y',\eta) \\
    &\hskip7.5cm\times 
             \Big[\psi_k(s)\underset{(s)}{*} h_m(s,y')
             \Big] 
           dy' \bigg\|_{L^p_{x'}}ds\bigg\}^p
           d\eta
       \Bigg)^{1/p}
  \\
    \le &C\sum_{m\in\Z}  2^{sm}
       \biggl(\int_{\re_+}
        \Bigl\{  \sum_{k\ge 2m} 
             \int_{\re_+}
             \Big\|
                  \Psi_{D,k,m}(t-s,\cdot,\eta)
             \Big\|_{L^1_{x'}}
             \Big\|
                  \psi_k(s)\underset{(s)}{*} h_m(s,\cdot)  
             \Big\|_{L^p_{x'}}
             ds
        \Bigr\}^p  
         d\eta
        \biggr)^{1/p}.
     \addtocounter{equation}{1}\tag{\theequation} \label{eqn;L_1}
\end{align*}
}
Then by the almost orthogonal estimate between the 
boundary potential $\Psi_{D}$  and the Littlewood-Paley 
decomposition $\psi_k$ in time, 
namely we invoke Lemma \ref{lem;pt-orthogonal-1} below,
for any $t, \eta\in \re_+$,
\begin{equation}   \label{eqn;Ortho-D0}
   \bigl\|
       \Psi_{D,k,m}(t,\cdot,\eta)
   \bigr\|_{L^1_{x'}}
   \le 
   \left\{
   \begin{aligned}
       & C 2^{k}
          \exp\big(-2^{\frac{k}{2}-1}\eta \big) 
          \frac{2^k}{\<2^kt\>^2},
          \quad  {k\ge 2m}, \\
       & C 2^{k} 
         \exp\big(-2^{m-1}\eta\big) 
         \frac{2^k}{\<2^kt\>^2},
         \quad  {k<2m}.
  \end{aligned}
  \right.
\end{equation}

Noting the restriction  {$k\ge 2m$} on  the time-dominated region, 
we  apply \eqref{eqn;Ortho-D0} to 
\eqref{eqn;L_1} 
and obtain that
{\allowdisplaybreaks
\begin{align*} 
 &\big\| L_1 \big\|_{L^1_t(\re_+)} \\
  \le& C\Bigl\| 
         \sum_{m\in \Z}  2^{sm}
         \Bigl(\int_{\re_+}
         \Bigl\{ \sum_{k\ge 2m} 
               \big(2^{k} \exp(-2^{\frac{k}{2}-1}\eta)\big)
         \int_{\re}
              \frac{2^k}{\<2^k(t-s)\>^2}
              \big\|
                 \psi_k\underset{(s)}{*}h_m(s,\cdot) 
              \big\|_{L^p_{x'}} 
              ds 
          \Bigr\}^p
           d\eta
          \Bigr)^{1/p}
         \Bigr\|_{L^1_t(\re_+)} 
 \\
    =&C\Bigl\| 
         \sum_{m\in \Z}  2^{sm}
         \Bigl\{ \sum_{k\ge 2m}
         2^{k}\int_{\re}
              \frac{2^{k}}{\<2^{k}(t-s)\>^2}
              \big\|
                   \psi_{k}\underset{(s)}{*}h_m(s,\cdot) 
              \big\|_{L^p_{x'}} 
              ds      
          \Big(    
          \int_{\re_+} \exp(-p2^{\frac{k}{2}-1}\eta) d\eta
          \Big)^{1/p}
          \Bigr\}   
         \Bigr\|_{L^1_t(\re_+)} 
 \\
    \le&C\Bigl\| 
         \sum_{m\in \Z}  2^{sm}
         \sum_{k\in\Z}
         2^{(1-\frac{1}{2p})k}
         \int_{\re}
              \frac{2^{k}}{\<2^{k}(t-s)\>^2}
              \big\|
                   \psi_{k}\underset{(s)}{*}h_m(s,\cdot) 
              \big\|_{L^p_{x'}} 
              ds 
         \Bigr\|_{L^1_t(\re_+)} 
   \\
       \le&C  \sum_{k\in\Z}
               2^{(1-\frac{1}{2p})k}
         \sum_{m\in \Z}  2^{sm}
         \Bigl\|
         \int_{\re}
              \frac{2^{k}}{\<2^{k}(t-s)\>^2}
              \big\|
                   \psi_{k}\underset{(s)}{*}h_m(s,\cdot) 
              \big\|_{L^p_{x'}} 
              ds 
         \Bigr\|_{L^1_t(\re_+)}
     \\
       \le&C  \sum_{k\in\Z}
               2^{(1-\frac{1}{2p})k}
              \sum_{m\in \Z}  2^{sm}
         \Bigl\|
              \big\|
                   \psi_{k}\underset{(s)}{*}h_m(s,\cdot) 
              \big\|_{L^p_{x'}}  
         \Bigr\|_{L^1_t(\re_+)} 
 \\
    =&C \Bigl\|\sum_{k\in\Z}
               2^{(1-\frac{1}{2p})k}
              \sum_{m\in \Z}   2^{sm} 
              \big\|
                   \psi_{k}\underset{(s)}{*}h_m(s,\cdot) 
              \big\|_{L^p_{x'}}  
        \Bigr\|_{L^1_t(\re_+)} 
   \\
  =& C\big\| h 
      \big\|_{\dF^{1-1/2p}_{1,1}(\re_+;\dB^s_{p,1}(\re^{n-1}_{x'}))}.
          \addtocounter{equation}{1}\tag{\theequation}\label{eqn;L_1-2}
\end{align*}
}
 Meanwhile the second term  $L_2$ is 
the  space-dominated region  and letting 
$h_m(t,x')\equiv \widetilde{\phi_m}*h(t,x')$,
we apply again the Minkowski inequality, the Hausdorff-Young inequality 
and the estimate  \eqref{eqn;chi}, 
\begin{align*} 
&\big\| L_2 \big\|_{L^1_t(\re_+)}\nonumber\\ 
  \le& C\bigg\| 
         \sum_{m\in \Z}  2^{sm}
         \Bigl(\int_{\eta\in \re_+}
          \Bigl\{ \sum_{k<2m}
               \big(2^{k}e^{-2^{m-1}\eta}\big)
             \int_{\re}
              \frac{2^k}{\<2^k(t-s)\>^2}
              \big\| h_m(s,\cdot)  \big\|_{L^p_{x'}}
              ds 
          \Bigr\}^p  
           d\eta
          \Bigr)^{1/p}
        \bigg\|_{L^1_t(\re_+)} 
 \\
   =& C\bigg\| 
         \sum_{m\in \Z}  2^{sm}
         \sum_{k<2m}
               2^{k}
             \int_{\re}
              \frac{2^k}{\<2^k(t-s)\>^2}
              \big\| h_m(s,\cdot)  \big\|_{L^p_{x'}}
              ds  
          \Big(
           \int_{\eta\in \re_+}\exp(-p2^{m-1}\eta)d\eta
          \Big)^{1/p}
        \bigg\|_{L^1_t(\re_+)} 
 \\
    \le& C\bigg\| 
         \sum_{m\in \Z}  2^{sm+2m-\frac{m}{p}}
         \sum_{k<2m}
               2^{k-2m}
             \int_{\re}
              \frac{2^k}{\<2^k(t-s)\>^2}
              \big\| h_m(s,\cdot)  \big\|_{L^p_{x'}}
              ds             
        \bigg\|_{L^1_t(\re_+)} 
 \\
     \le& C \sum_{m\in \Z} 2^{sm+2m-\frac{m}{p}}
            \sum_{k<2m}
               2^{k-2m}
        \bigg\| 
             \int_{\re}
              \frac{2^k}{\<2^k(t-s)\>^2}
              \big\| h_m(s,\cdot)  \big\|_{L^p_{x'}}
              ds             
        \bigg\|_{L^1_t(\re_+)} 
 \\
    \le & C 
         \Bigl\|
         \sum_{m\in \Z}
          2^{(s+2-\frac{1}{p})m}               
         \sum_{k<2m}               
               \big(2^{k-2m}\big)
               \big\| h_m(t,\cdot)  \big\|_{L^p_{x'}}
         \Bigr\|_{L^1_t(\re_+)} 
 \\
  \le&  C\big\|h\big\|_{ L^1 (\re_+;\dB^{s+2-1/p}_{p,1}(\re^{n-1}_{x'}) )}.
    \eqntag\label{eqn;L_2-2}
\end{align*}
From \eqref{eqn;potential-besov-10}, \eqref{eqn;L_1-2} and 
\eqref{eqn;L_2-2}, the estimate \eqref{eqn;P1-est} is shown. 
\footnote{Up to this level there is no restriction on $p$ nor $s$.}
\par

We then prove \eqref{eqn;P2-est}.  By \eqref{eqn;space-time-splitting}, 
we split $P_2^D$ into the time-like region and the space-like region;
{\allowdisplaybreaks
\begin{align*} 
 P_2^D
 =&C\sum_{m\in \Z} 2^{sm} \Big(\int_{\re_+}
               \Big\|
                   \phi_{m}(\eta)\underset{(\eta)}{*}
                   \Psi_D(t-s,x',\eta) 
                   \underset{(t,x')}{*} \\
         &\hskip3cm \times 
                   \sum_{k\in \Z} \sum_{j\le m+1}
                   \psi_k(t)\underset{(t)}{*}\phi_j(y')
                            \underset{(y')}{*} h(t,y')
                 \Big\|_{L^p(\re^{n-1}_{x'})}^p  d\eta \Big)^{1/p}
 \\
 \le& C\sum_{m\in \Z} 2^{sm}
      \Big\| 
        \Big\|
         \sum_{k> 2m}\sum_{j\le m+1}
         \phi_m(\eta)\underset{(\eta)}{*} 
         \Psi_D(t,x',\eta)
         \underset{(t,x')}{*}        
         \psi_k(t)\underset{(t)}{*}\phi_j(x')
         \underset{(x')}{*}h(t,x')  
 \\
   &\hskip15mm+
         \sum_{k\le 2m}\sum_{2j\le k}
         \phi_m(\eta)\underset{(\eta)}{*} 
         \Psi_D(t,x',\eta)
         \underset{(t,x')}{*}        
         \psi_k(t)\underset{(t)}{*}\phi_j(x')
         \underset{(x')}{*}h(t,x')  
        \Big\|_{L^p_{x'}}
      \Big\|_{L^p(\re^+_{\eta})}  \\
   &+C\sum_{m\in \Z} 2^{sm}
       \Big\|
         \Big\| 
         \sum_{k\le 2m}\sum_{k<2j\le 2m+2}
          \phi_m(\eta)\underset{(\eta)}{*} 
          \Psi_D(t,x',\eta)
          \underset{(t,x')}{*} \\
      &\hskip7cm     
          \times\psi_k(t)\underset{(t)}{*}\phi_j(x')
           \underset{(x')}{*}h(t,x')  
         \Big\|_{L^p_{x'}}
       \Big\|_{L^p(\re^+_{\eta})}
\\
 \le& C\sum_{m\in \Z} 2^{sm}
      \Big\|         
         \sum_{k\in\Z}\sum_{2j\le \min(k,2m+2)}
         \Big\|\phi_m(\eta)\underset{(\eta)}{*} 
         \Psi_D(t,x',\eta)
         \underset{(t,x')}{*}      \\
      &\hskip7cm        
         \times  \psi_k(t)\underset{(t)}{*}\phi_j(x')
         \underset{(x')}{*}h(t,x')  
        \Big\|_{L^p_{x'}}
      \Big\|_{L^p(\re^+_{\eta})} 
 \\
    &+ C\sum_{m\in \Z} 2^{sm}
       \Big\|        
          \sum_{k\in\Z}\sum_{2j\le 2m+2}
          \Big\| \phi_m(\eta)\underset{(\eta)}{*} 
          \Psi_D(t,x',\eta)
          \underset{(t,x')}{*}  \\
      &\hskip7cm     
          \times \psi_k(t)\underset{(t)}{*}\phi_j(x')
           \underset{(x')}{*}h(t,x')  
         \Big\|_{L^p_{x'}}
       \Big\|_{L^p(\re^+_{\eta})}
  \\
  \equiv & M_1 +  M_2 .
 \eqntag 
 \label{eqn;potential-besov-11}
\end{align*}
}
The first term  $ M_1 $ of the right hand side 
is the time-dominated part, letting 
$h_j\equiv \widetilde{\phi_j}*h$, 
we apply the Minkowski inequality and the Hausdorff-Young inequality with 
\eqref{eqn;chi} to see
We use the  almost orthogonal estimate \eqref{eqn;Ortho-D0}
between $\psi_m$ and $\Psi_{D,k,m}$ (Lemma \ref{lem;pt-orthogonal-D2}):
For any $N\in \Nt$,
\eqn{
 \big\|\phi_m(\eta)\underset{(\eta)}{*}\Psi_{D,k,j}(t,\cdot,\eta)
 \big\|_{L^1(\re^{n-1}_{x'})} 
 \le
 \left\{
 \spl{
     &\frac{ C_N2^k 2^{-|\frac{k}{2}-m|} } 
           { \<  {2^{ \min(\frac{k}{2},m)}} \eta\>^N}
           \frac{2^k}{\<2^k t\>^2}, 
     \quad { k\ge 2j, }\\
     &\frac{ C_N2^k2^{-|j-m|} }
           {\<{2^{j}}\eta\>^N}
           \frac{2^k}{\<2^k t\>^2}, 
     \quad {k<2j }
     }
\right.
}
for some $C_N>0$.
Then setting $2^m\eta=\tilde{\eta}$ and shifting 
$(m,k,j)\to (m',k,j)$ by $m-\frac{k}{2}=m'$, 
the first term of the right hand side of 
\eqref{eqn;potential-besov-11} can be estimated as follows:

\begin{align*}  
 \big\| &  M_1\big\|_{L^1_t(\re_+)}  
 \\
  \le&\bigg\|\sum_{m\in\Z} 2^{sm}
       \Bigl(\int_{\re_+} 
        \Bigl\{  \sum_{k\in\Z}\sum_{2j\le \min(k,2m+2)} 
             \int_{\re_+}
               \Big\|\phi_m(\eta)\underset{(\eta)}{*} 
                     \Psi_{D,k,j}(t-s,x',\zeta)
               \Big\|_{L^1_{x'}}\\
    &\hskip 8cm \times             
              \Big\|\psi_k\underset{(s)}{*} h_j(s,x')  \Big\|_{L^p_{x'}}
             ds
        \Bigr\}^p d\eta
     \Bigr)^{1/p}
     \bigg\|_{L^1_t(\re_+)}  
\\
  \le&C\bigg\|\sum_{m\in\Z}  2^{sm} \bigg(\int_{\re_+}            
              \Bigl\{ 
              \sum_{k\in\Z}               
                \Big(
                \frac{2^{-|\frac{k}{2}-m|}}
                     {\< {2^{\min(\frac{k}{2},m)}} |\eta|\>^N}
                \Big) 
                2^{k}
              \int_{\re}
              \frac{2^k}{\<2^k(t-s)\>^2}
              \sum_{2j\le k}   
              \Big\|
                 \psi_k \underset{(s)}{*} h_j(s,\cdot) 
              \Big\|_{L^p_{x'}}
             ds 
             \Bigr\}^pd\eta \bigg)^{1/p} 
     \bigg\|_{L^1_t(\re_+)} 
\\
  \le&C\bigg\|\sum_{m\in\Z}  2^{sm} 
              \sum_{k\in\Z}               
              2^{-|\frac{k}{2}-m|} 
              2^{k}
             \int_{\re}
             \frac{2^k}{\<2^k(t-s)\>^2}
             \sum_{2j\le k}  
             \Big\|
             \psi_k \underset{(s)}{*} h_j(s,\cdot) 
             \Big\|_{L^p_{x'}}
             ds 
             \left(\int_{\re_+}
             \Big(\frac{1}{\< {2^{\min(\frac{k}2,m)}}  |\eta|\>^N}
             \Big)^p
             d\eta
             \right)^{1/p} 
     \bigg\|_{L^1_t(\re_+)} 
 \\
  \le&C\bigg\|\sum_{m\in\Z}    2^{sm}
              \sum_{k\in\Z}               
             2^{-|\frac{k}{2}-m|} 
             2^{k} {2^{-\frac{1}{p}\min(\frac{k}2,m)}}
             \int_{\re}
             \frac{2^k}{\<2^k(t-s)\>^2}
            \sum_{2j\le k}  
             \Big\|
             \psi_k \underset{(s)}{*} h_j(s,\cdot) 
             \Big\|_{L^p_{x'}}
             ds  
      \bigg\|_{L^1_t(\re_+)} 
 \\
 \le&C \sum_{k\in\Z}  
              2^{(1-\frac{1}{2p})k}
          \Big( \sum_{m'\in\Z}
              2^{-|m'|+sm'+{\max(0,-\frac{m'}{p})}}
          \Big)  
              2^{\frac{s}{2}k}                               
       \bigg\|\int_{\re}
              \frac{2^k}{\<2^k(t-s)\>^2}
              \sum_{2j\le k}  
              \Big\|
              \psi_k \underset{(s)}{*} h_j(s,\cdot)
              \Big\|_{L^p_{x'}}
             ds 
       \bigg\|_{L^1_t(\re_+)} 
 \\
 \le&C \sum_{k\in\Z}  
              2^{(1-\frac{1}{2p})k}                                
       \bigg\|
       \sum_{2j\le k}  
               2^{sj} 
              \Big\|
              \psi_k \underset{(s)}{*} h_j(s,\cdot)
              \Big\|_{L^p_{x'}}
       \bigg\|_{L^1_t(\re_+)} 
 \\
  \le&C\bigg\| \sum_{k\in\Z}   
          2^{(1-\frac{1}{2p})k}  
              \Big\|
              \psi_k \underset{(s)}{*} h(s,\cdot)
              \Big\|_{\dB^s_{p,1}(\re^{n-1}_{x'})}
         \bigg\|_{L^1_t(\re_+)}
   =C\big\|h\big\|_{\dF^{1-1/2p}_{1,1}(\re_+;\dB_{p,1}^s(\re^{n-1}_{x'}))}, 
   \addtocounter{equation}{1}\tag{\theequation}\label{eqn;M_1-2}
\end{align*}
where we used the fact that
$|s-1/p|<1$,  i.e., 
 $-1+1/p<s< 1+1/p$  and $s\le 0$
at the second line from the bottom.

On the other hand for the estimate $ M_2$, we proceed 
a similar way to treat $ M_1$.  Exchanging the order 
of the summation of $j$ and $k$ and setting 
$
 h_j(t)\equiv \phi_j\underset{(x')}{*}h(t),
$  
it follows by  changing $(m,k,j)\to (m,k,j)$ 
with $m-j=m'$ and \eqref{eqn;potential-besov-11} that
{\allowdisplaybreaks
\begin{align*} 
  &\big\| M_2 \big\|_{L^1} \\
  \le&\bigg\|\sum_{m\in\Z} 2^{sm}
       \bigg(\int_{\re_+} 
        \Bigl\{  \sum_{j\in\Z}\sum_{k< 2j}
          \int_{\re_+}
            \Big\| \phi_m(\eta)\underset{(\eta)}{*} 
                   \Psi_{D,k,j}(t-s,x',\eta)
            \Big\|_{L^1_{x'}}
            \big\| h_j(s,\cdot)  \big\|_{L^p_{x'}}
             ds
        \Bigr\}^p d\eta
     \bigg)^{1/p}
     \bigg\|_{L^1_t(\re_+)}  
\\
  \le&C\bigg\|\sum_{m\in\Z}  2^{sm}
             \bigg(\int_{\re_+}            
             \bigg\{ 
               \sum_{j\in\Z}\sum_{k< 2j}   
               2^k           
              \Bigl(
               \frac{ C_N 2^{-|j-m|} }
                    {\<{2^j} |\eta|\>^N}
              \int_{\re}
              \frac{2^k}{\<2^k(t-s)\>^2}
              \big\| h_j(s,\cdot)\big\|_{L^p_{x'}}
              ds 
             \bigg\}^pd\eta \Bigr)^{1/p} 
      \bigg\|_{L^1_t(\re_+)} 
\\
  \le&C\bigg\|\sum_{m\in\Z}  2^{sm}
               \int_{\re_+}            
               \sum_{j\in\Z}\sum_{k< 2j}   
               2^k           
               2^{-|j-m|} 
              \int_{\re}
              \frac{2^k}{\<2^k(t-s)\>^2}
              \big\| h_j(s,\cdot)\big\|_{L^p_{x'}}
                  ds 
              \left(
              \int_{\re_+}                     
               \frac{1}{\<{2^j} |\eta|\>^{pN}}
               d\eta
              \right)^{1/p}
       \bigg\|_{L^1_t(\re_+)} 
\\
  &\text{(changing the variable $2^j\eta=\tilde{\eta}$)} 
\\
   \le &C \bigg\|\sum_{m\in\Z}  2^{sm}
                \sum_{j\in\Z}
                2^{-|j-m|}{2^{-\frac{j}{p}}}
                2^{2j}
             \sum_{k<2j} 2^{k-2j}
             \int_{\re}
              \frac{2^{k}}{\<2^{k}(t-s)\>^2}
              \Big\| h_j(s,\cdot)
              \Big\|_{L^p_{x'}} 
             ds              
          \left(\int_{\re_+} 
              \frac{1}{\<|\tilde{\eta}|\>^{pN}}  
              d\tilde{\eta}
          \right)^{1/p}
        \bigg\|_{L^1_t(\re_+)}  
\\
   \le &C\bigg\| \sum_{j\in\Z}
         \Big(\sum_{m'\in\Z}  
              2^{sm'} 
              2^{-|m'|}\Big)
              2^{sj}2^{-\frac{j}{p}} 2^{2j}           
             \sum_{k<2j} 2^{k-2j}
              \int_{\re}
              \frac{2^{k}}{\<2^{k}(t-s)\>^2}
              \Big\| h_j(s,\cdot)
              \Big\|_{L^p_{x'}} 
             ds              
        \bigg\|_{L^1_t(\re_+)}
\\
   \le &C \sum_{j\in\Z}
              2^{(s+2-\frac{1}{p})j}            
          \sum_{k<2j} 2^{k-2j}
          \bigg\| 
              \int_{\re}
              \frac{2^{k}}{\<2^{k}(t-s)\>^2}
              \big\| h_j(s,\cdot)
              \big\|_{L^p_{x'}} 
             ds              
        \bigg\|_{L^1_t(\re_+)}
\\
   \le &C\bigg\| 
          \sum_{j\in\Z}
              2^{(s+2-\frac{1}{p})j}            
         \Big(\sum_{k<2j} 2^{k-2j} \Big)        
              \big\| h_j(s,\cdot)
              \big\|_{L^p_{x'}}             
         \bigg\|_{L^1_t(\re_+)}  
  \\  
     \le &C\Bigl\| 
           \sum_{j \in\Z} 2^{(s+2-\frac{1}{p})j}
               \big\|h_j(t,\cdot)\big\|_{L^p_{x'}}             
          \Bigr\|_{L^1_t(\re_+)} 
  = \big\|h\big\|_{L^1(\re_+;\dB^{s+2-1/p}_{p,1}(\re^{n-1}_{x'}))}.
   \addtocounter{equation}{1}\tag{\theequation}\label{eqn;M_2-3}
\end{align*} 
}
Here we used the fact that {$|s|<1$} for convergence of 
the summation on $m'$
In the last estimate, we exchange the order of the integration 
in time and the summation of $\ell'$ and use the Hausdorff-Young 
inequality to remove the convolution with the time potential term
and then recovers the time integration out side.
From \eqref{eqn;potential-besov-11}, \eqref{eqn;M_1-2} and \eqref{eqn;M_2-3} 
the estimate \eqref{eqn;P2-est} is shown. This completes the proof of Lemma 
\ref{lem;P1-P2-estimate}. 
\end{prf}

%
%
\subsection{Proof of Theorem \ref{thm;Lp-boundary-D}}\par
For the proof of Theorem \ref{thm;Lp-boundary-D}, it is simpler than
the proof of  Theorem \ref{thm;boundary-trace-D} since it does not 
involve the Littlewood-Paley decomposition in $(x',\eta)$-variable.   
We recall $I_{-\ell}=[2^{-\ell}, 2^{-\ell+1})$ for all $\ell\in \Z$.
We again split the case for the time-dominated part and space-dominated part 
and proceed the estimate:
{\allowdisplaybreaks
\algn{ 
   &\int_0^{\infty} \|\Del u(t)\|_{L^p(\re^{n}_+)}dt  \\
      =\, &\Bigl\|\Bigl\|\Bigl(\int_{\re_+}
                  \Big|                   
                  \int_0^{\infty}\int_{\re^{n-1}_{x'}}
                   \Psi_D(t-s,x'-y',\eta) h(s,y')ds dy'
                   \Big|^pd\eta
                  \Bigr)^{1/p} 
           \Bigr\|_{L^p(\re^{n-1}_{x'})}
           \Bigr\|_{L^1_t(\re_+)} 
\\
      =\, &\Bigl\|\Bigl\|\Bigl(\sum_{\ell\in\Z}
                  \int_{\ell\in I_{-\ell}}
                  \Big|                   
                  \int_0^{\infty}\int_{\re^{n-1}_{x'}}
                   \Psi_D(t-s,x'-y',\eta) h(s,y')ds dy'
                   \Big|^p 
                   d\eta
                  \Bigr)^{1/p} 
            \Bigr\|_{L^p(\re^{n-1}_{x'})}
           \Bigr\|_{L^1_t(\re_+)}
\\
   \le&\, C   \Bigl\| \Bigl(\sum_{\ell\in \Z}  
             {2^{-\ell} }
                 \big\| 
                  \Psi_D(t,x',\eta)\Big|_{\eta\simeq 2^{-\ell}}
                  \underset{(t,x')}{*} h(t,x')  
                 \big\|_{L^p(\re^{n-1}_{x'})}^p 
                \Bigr)^{1/p}
      \Bigr\|_{L^1_t(\re_+)},
   \eqntag\label{eqn;boudary-Lp-1}
}
}
where we  split the last term  of \eqref{eqn;boudary-Lp-1} 
into two parts: 
\algn{ \label{eqn;boudary-Lp-2} 
 \Big\|& \Psi_{D}(t,x',\eta)  \Big|_{\eta\simeq 2^{-\ell}}
        \underset{(t,x')}{*} h(t,x')  
 \Big\|_{L^p(\re^{n-1}_{x'})}\nonumber\\
  \le & \big\| \Psi_{D}(t,x',\eta)
         \underset{(t,x')}{*} 
         \sum_{k\in\Z}\sum_{\widehat{0}\le j\le k/2}
         \psi_k(t)\underset{(t)}{*}\phi_j(x')
         \underset{(x')}{*}h(t,x')  
       \big\|_{L^p_{x'}} \nonumber\\
    &+\big\|  \Psi_{D}(t,x',\eta)
        \underset{(t,x')}{*} 
        \sum_{j\ge\widehat{0}}\sum_{k<2j}
        \psi_k(t)\underset{(t)}{*}\phi_j(x')\underset{(x')}{*}h(t,x')  
      \big\|_{L^p_{x'}} 
\\
  \le & \sum_{k\in\Z}
        \Bigl\|
           \int_{\re^{n-1}}
            \int_{\re}
            \sum_{\widehat{0}\le j\le k/2}
             \Psi_{D,k,j}(t-r,x'-y')
              \Big[\psi_k\underset{(t)}{*}\phi_j(y')
              \underset{(y')}{*}h(r,y')
              \Big] 
              drdy'
        \Bigr\|_{L^p_{x'}}    \\
     & +\sum_{j\ge\tilde{0}}
        \Bigl\|
           \int_{\re^{n-1}}\int_{\re}\sum_{k<2j}
               \Psi_{D,k,j}(t-s,x'-z')
              \Big[\phi_j\underset{(x')}{*}
                   \psi_k\underset{(s)}{*}h(s,z')
              \Big] 
              dsdz'
        \Bigr\|_{L^p_{x'}}
 \nonumber\\
    \le& \sum_{k\in\Z}
          \int_{\re} 
              \sup_{j\in \Z} 
             \Big\|
               \Psi_{D,k,j}(t-r,\cdot)
             \Big\|_{L^1_{x'}}
             \Big\|
                 \sum_{\widehat{0}\le j\le k/2}\phi_j\underset{(x')}{*}
                 \big( \psi_k\underset{(r)}{*}h(r,\cdot)\big)
             \Big\|_{L^p_{x'}} 
         dr\\
     & +\sum_{j\ge\tilde{0}}
          \int_{\re} 
             \Big\|\sum_{k<2j}
               \Psi_{D,k,j}(t-s,\cdot)
             \Big\|_{L^1_{x'}}
             \Big\|
                \phi_j\underset{(\cdot)}{*}
                   h(s,\cdot) 
             \Big\|_{L^p_{x'}} ds 
 \\
 \equiv &D^1_{\ell}+ D^2_{\ell}.   \eqntag 
}
Here the notation $j=\widehat{0}$ for the Littlewood-Paley decomposition is
given in \eqref{eqn;LP-decomp}. 
Thus from \eqref{eqn;boudary-Lp-1}, \eqref{eqn;boudary-Lp-2},  we estimate 
the following two terms:
\begin{align*}
   \int_0^{\infty} \|\Del u(t)\|_{L^p(\re^{n})}dt 
   \le  & C \Bigl\| \Bigl(\sum_{\ell\in \Z} 
               {2^{-\ell} }
                 \big\| 
                  \Psi_{D}(t,x',\eta)
                   \underset{(t,x')}{*}  h(t,x')  
                 \big\|_{L^p(\re^{n-1})}^p 
                 \Big|_{\eta\in I_{-\ell}}
                \Bigr)^{1/p}
           \Bigr\|_{L^1_t(\re_+)}\nonumber\\
\le  & C\Big\|\| \{ D^1_{\ell}\}_{\ell\in\Z} \|_{\ell^{-1/p}_p}
         \Big\|_{L^1_t(\re_+)}
       +C\Big\|\|\{ D^2_{\ell}\}_{\ell\in\Z} \|_{\ell^{-1/p}_p}
         \Big\|_{L^1_t(\re_+)}.
\end{align*}
We now see that the following proposition is enough to ensure that  
Theorem \ref{thm;Lp-boundary-D} hold.
\vskip2mm
\begin{prop}\label{prop;L1'-L2'-estimate} Let $1\le p\le \infty$. 
 Let  $ D^1_{\ell}$ and $D^2_{\ell}$  be defined in 
\eqref{eqn;boudary-Lp-2}. 
Then there exists  constants $C>0$ such that the following estimate hold:
\begin{align}
 &\Big\|\| \{ D^1_{\ell}\}_{\ell\in \Z} \|_{\ell^{-1/p}_p}
  \Big\|_{L^1_t(\re_+)}
  \le C\big\|h \big\|_{ \dot{F}^{1-1/2p}_{1,1}(\re_+;L^p(\re^{n-1}))}, 
       \label{eqn;L'1-est-star}\\
 &\Big\|\| \{D^2_{\ell}\}_{\ell\in\Z} \|_{\ell^{-1/p}_p}
  \Big\|_{L^1_t(\re_+)}
  \le C\big\|h \big\|_{ L^1 (\re_+;B^{2-1/p}_{p,1}(\re^{n-1})) }. 
       \label{eqn;L'2-est-star}
\end{align}
\end{prop}

\begin{prf}{Proposition \ref{prop;L1'-L2'-estimate}}
The estimates \eqref{eqn;L'1-est-star} and \eqref{eqn;L'2-est-star} 
are proved in the similar manner as 
\eqref{eqn;P1-est} and \eqref{eqn;P2-est} in the case $1\le p<\infty$.  
Setting $k-2\ell=-2\ell'$ to change $(\ell,k,j)\to (\ell',k,j)$,
\algn{ \eqntag\label{eqn;boudary-Lp-M1-*}
       & \left\| \left(
                   \sum_{\ell\in \Z} |{2^{- \frac{\ell}{p}}} D^1_{\ell}|^p 
                   \right)^{\frac{1}{p}}
          \right\|_{L^1_t(\re_+)}
\\
  \le & 
        \left\| 
        \left\{
        \sum_{\ell\in \Z} {   2^{(-\frac{1}{p})p\ell}  }
        \left(\sum_{k\in\Z}
          \int_{\re} 
              \sup_{j\in \Z} 
             \Big\|
               \Psi_{D,k,j}(t-r,\cdot)  
             \Big\|_{L^1_{x'}}
             \Big\|
                 \sum_{\widehat{0}\le j\le k/2}\psi_k\underset{(t)}{*}
                    \phi_j\underset{(x')}{*}h(r,\cdot)
             \Big\|_{L^p_{x'}} 
         dr    
        \right)^p
        \right\}^{\frac{1}{p}}
          \right\|_{L^1_t(\re_+)} 
\\
  \le& C\left\| \left\{
           \sum_{\ell\in \Z} 
             \Big( \sum_{k\in\Z}
               2^{(2-\frac{1}{p})\ell}
               \big(2^{k-2\ell}e^{-2^{\frac12(k-2\ell)}}\big)
             \int_{\re}
             \frac{2^k}{\<2^k(t-r)\>^2}
             \Big\|
               \psi_k*h(r,\cdot) 
             \Big\|_{L^p_{x'}}
             dr 
             \Big)^p\right\}^{1/p}
         \right\|_{L^1_t(\re_+)} 
\\
  \le& C\left\| \left\{
           \sum_{\ell\in \Z} 
             \Big( \sum_{k\in\Z}
               \big(2^{k-\frac{\ell}{p}}e^{-2^{\frac12(k-2\ell)}}\big)
             \int_{\re}
             \frac{2^k}{\<2^k(t-r)\>^2}
             \Big\|
               \psi_k\underset{(r)}{*}h(r,\cdot)             
             \Big\|_{L^p_{x'}}
             dr 
             \Big)^p\right\}^{1/p}
         \right\|_{L^1_t(\re_+)}
\\
  =& C\left\| \left\{
           \sum_{\ell'\in \Z/2} 
             \Big( \sum_{k\in\Z}
               \big(2^{(1-\frac1{2p})k-\frac{1}{p}\ell'}e^{-2^{-\ell'}}\big)
             \int_{\re}
             \frac{2^k}{\<2^k(t-r)\>^2}
             \Big\|
               \psi_k\underset{(r)}{*}h(r,\cdot)            
             \Big\|_{L^p_{x'}}
             dr 
             \Big)^p\right\}^{1/p}
         \right\|_{L^1_t(\re_+)}
 \\
  =& C\left\| \left\{
             \sum_{\ell'\in \Z/2}  
             \big(2^{-\frac{1}{p}\ell'}e^{-2^{-\ell'}}\big)^p
             \Big( \sum_{k\in\Z}
              2^{(1-\frac1{2p})k}
             \int_{\re}
             \frac{2^k}{\<2^k(t-r)\>^2}
             \Big\|
               \psi_k\underset{(r)}{*}h(r,\cdot)             
             \Big\|_{L^p_{x'}}
             dr 
             \Big)^p\right\}^{1/p}
         \right\|_{L^1_t(\re_+)}
 \\
 \le& C\sum_{k\in\Z}
             2^{(1-\frac1{2p})k}
            \left\| \int_{\re}
             \frac{2^k}{\<2^k(t-r)\>^2}
             \Big\|
               \psi_k\underset{(r)}{*}h(r,\cdot)             
             \Big\|_{L^p_{x'}}
             dr 
         \right\|_{L^1_t(\re_+)}
 \\
 \le& C\sum_{k\in\Z}
             2^{(1-\frac1{2p})k}
            \left\| 
             \Big\|
               \psi_k\underset{(t)}{*}h(t,\cdot)            
             \Big\|_{L^p_{x'}}
         \right\|_{L^1_t(\re_+)} 
   \le C\left\| 
           \sum_{k\in \Z}  
              2^{(1-\frac{1}{2p})k}
              \Big\|\psi_{k}\underset{(t)}{*}h(t,\cdot)\Big\|_{L^p_{x'}}
         \right\|_{L^1_t(\re_+)} 
\\
   \le &C
       \big\|h
       \big\|_{\dot{F}^{1-1/(2p)}_{1,1}(\re_+;L^p_{x'})}.
}

%
%
For the case $p=\infty$  in \eqref{eqn;L'1-est-star}, 
setting $k-2\ell= -2\ell'$ and changing $\ell\to \ell'$, we have 
{\allowdisplaybreaks
\begin{align*}
    & \Bigl\| \sup_{\ell\in \Z}                    | D^1_{\ell}|
      \Bigr\|_{L^1_t(\re_+)}\\
  \le& C\Bigl\| 
             \sup_{\ell\in \Z} 
             \sum_{k\in\Z}
               2^{2\ell}
               \big(2^{k-2\ell}e^{-2^{\frac12(k-2\ell)}}\big)
             \int_{\re}
             \frac{2^k}{\<2^k(t-r)\>^2}
             \Big\|
               \sum_{\widehat{0}\le j\le k/2}\phi_j\underset{(\cdot)}{*}h_k(r,\cdot)
             \Big\|_{L^{\infty}_{x'}}
             dr 
         \Bigr\|_{L^1_t(\re_+)} 
 \\
      =&C\Bigl\| 
           \sup_{\ell'\in \Z} 
            \sum_{k\in\Z}
              2^{k+2\ell'}
               \big(2^{-2\ell'}e^{-2^{-\ell'}}\big)
             \int_{\re}
             \frac{2^{k}}{\<2^{k}(t-r)\>^2}
             \Big\|
                \sum_{\widehat{0}\le j\le k/2}
                \phi_j\underset{(\cdot)}{*}h_{k}(r,\cdot)
             \Big\|_{L^{\infty}_{x'}} 
             dr
         \Bigr\|_{L^1_t(\re_+)} 
  \\
      =&C\Bigl\| 
           {\sup_{\ell'\in \Z}  \big(e^{-2^{-\ell'}}\big)}
            \sum_{k\in\Z}
              2^{k}              
             \int_{\re}
             \frac{2^{k}}{\<2^{k}(t-r)\>^2}
             \Big\|
                \sum_{\widehat{0}\le j\le k/2}
                \phi_j\underset{(\cdot)}{*}h_{k}(r,\cdot)
             \Big\|_{L^{\infty}_{x'}} 
             dr
         \Bigr\|_{L^1_t(\re_+)}
  \\
      =&C\Bigl\| 
            \sum_{k\in\Z}
              2^{k}              
             \int_{\re}
             \frac{2^{k}}{\<2^{k}(t-r)\>^2}
             \Big\|
                h_{k}(r,\cdot) 
             \Big\|_{L^{\infty}_{x'}} 
             dr
         \Bigr\|_{L^1_t(\re_+)} 
   \le C\big\|h\big\|_{\dot{F}^{1}_{1,1}(\re_+;L^{\infty}_{x'})}.
\end{align*}
}
To show the estimate \eqref{eqn;L'2-est-star} for $1\le p<\infty$, 
we derive it 
by setting $g_j(s,y')\equiv \phi_j*g(s,y')$ and then
$j-\ell=-\ell'$  and $(\ell,k,j)\to (\ell',k,j)$ that 
\begin{align*}
 &\left\| \left(
           \sum_{\ell\in \Z} {2^{(-\frac1p)p\ell}}
           | D^2_{\ell}|^p 
          \right)^{\frac{1}{p}}
  \right\|_{L^1_t(\re_+)}
\\
  \le & \left\| \left\{
           \sum_{\ell\in \Z} {2^{(-\frac1p)p\ell} }
             \Big( 
             \int_{\re} 
             \sum_{j\ge \widehat{0}}
             \sum_{k<2j}
             \Big\|
              {
               \Psi_{D,k,j}(t-r,\cdot)
              }
             \Big\|_{L^1_{x'}}
             \Big\|
                  h_j(r,\cdot) 
             \Big\|_{L^p_{x'}} 
              dr 
             \Big)^p\right\}^{1/p}
         \right\|_{L^1_t(\re_+)} 
\\
    \le &C\left\| \left\{
           \sum_{\ell\in \Z} 2^{(2-\frac1p)p\ell}
             \Big( \sum_{j\ge \widehat{0}} 
              \big(2^{ 2(j-\ell)} e^{-2^{(j-\ell)}}\big)
              \sum_{k<2j}2^{k- 2j}
              \int_{\re}
              \frac{2^k}{\<2^k(t-r)\>}
              \big\|h_j(r,\cdot)\big\|_{L^p_{x'}} 
              dr
             \Big)^p\right\}^{1/p}
         \right\|_{L^1_t(\re_+)} 
 \\
     = &C\left\| \left\{
              \sum_{\ell'\in \Z}
              \big(2^{-\frac{\ell'}p} e^{-2^{-\ell'} }\big)^p
              \bigg(
              \sum_{j\ge  \widehat{0}}              
                 2^{(2-\frac1p)j}               
              \sum_{k<2j}2^{k- 2j}
              \int_{\re}
              \frac{2^k}{\<2^k(t-r)\>}
              \big\|h_j(r,\cdot)\big\|_{L^p_{x'}} 
              dr
             \bigg)^p\right\}^{1/p}
         \right\|_{L^1_t(\re_+)} 
   \\
     \le &C\left\| 
             \sum_{j\ge  \widehat{0}}              
              2^{(2-\frac1p)j}              
              \sum_{k<2j}2^{k- 2j}
              \int_{\re}
              \frac{2^k}{\<2^k(t-r)\>}
              \big\|h_j(r,\cdot)\big\|_{L^p_{x'}} 
              dr
         \right\|_{L^1_t(\re_+)} 
   \\
     \le &C\sum_{j\ge  \widehat{0}}              
              2^{(2-\frac1p)j}              
              \sum_{k<2j}2^{k- 2j}
           \left\|
              \int_{\re}
              \frac{2^k}{\<2^k(t-r)\>}
              \big\|h_j(r,\cdot)\big\|_{L^p_{x'}} 
              dr
         \right\|_{L^1_t(\re_+)} 
    \\
     \le &C \Big\|\sum_{j\ge  \widehat{0}}              
              2^{(2-\frac1p)j}                        
              \big\|h_j(t,\cdot)\big\|_{L^p_{x'}} 
            \Big\|_{L^1_t(\re_+)} 
    \\
    =&C\big\|h\big\|_{L^1(\re_+;B^{2-1/p}_{p,1}(\re^{n-1}_{x'}))}.
\end{align*}
For the case $p=\infty$ in 
\eqref{eqn;L'2-est-star}, by setting $j-\ell=-\ell'$ to change 
$(\ell,j)\to (\ell',j)$, we see that 
{\allowdisplaybreaks 
\begin{align*}
  &\Bigl\| 
           \sup_{\ell\in \Z} 
           |D^2_{\ell}|
  \Bigr\|_{L^1_t(\re_+)}\\
  \le & C
         \Bigl\| 
           \sup_{\ell\in \Z} 2^{2\ell}
             \Bigl( 
             \int_{\re} 
             \sum_{j\ge \widehat{0}}
             \sum_{k<2j}
               2^{k-2\ell}
               e^{-2^{j-\ell}}
              \frac{2^k}{\<2^k(t-r)\>^2}
             \big\| h_j(r,\cdot)\big\|_{L^{\infty}_{x'}} 
              dr 
             \Bigr)
         \Bigr\|_{L^1_t(\re_+)} 
 \\
    \le &C\Bigl\| 
           \sup_{\ell\in \Z} 
             \Bigl( \sum_{\bk j\ge \widehat{0}} 
              2^{2\ell}
              \big(2^{ 2(j-\ell)} e^{-2^{(j-\ell)}}\big)
              \sum_{k<2j}2^{k-2j}
              \int_{\re}
              \frac{2^k}{\<2^k(t-r)\>^2}
              \big\|h_j(r,\cdot)\big\|_{L^{\infty}_{x'}} 
              dr
             \Bigr)
         \Bigr\|_{L^1_t(\re_+)} 
 \\
  = &C\Bigl\|
             \sup_{\ell' \in \Z} 
             \Bigl( \sum_{\bk j\ge  \widehat{0}}
             2^{2\ell'+2j}
                    2^{-2\ell'} e^{-2^{-\ell'}}
             \sum_{k<2j}2^{k-2j}
             \int_{\re}
              \frac{2^k}{\<2^k(t-r)\>^2}
              \big\|h_{j}(r,\cdot)\big\|_{L^{\infty}_{x'}}
              dr
             \Bigr)
         \Bigr\|_{L^1_t(\re_+)} 
  \\
    = &C\Bigl\|
             \sup_{\ell' \in \Z} e^{-2^{-\ell'}}
             \Bigl( \sum_{\bk j\ge  \widehat{0}}
                 2^{2j}  
             \sum_{k<2j}2^{k-2j}
             \int_{\re}
              \frac{2^k}{\<2^k(t-r)\>^2}
              \big\|h_{j}(r,\cdot)\big\|_{L^{\infty}_{x'}}
              dr
             \Bigr)
         \Bigr\|_{L^1_t(\re_+)}
  \\
    \le &C \sum_{\bk j\ge  \tilde{0}} 2^{2j}  
           \sum_{k<2j}2^{k-2j}
            \biggl\|
             \int_{\re}
              \frac{2^k}{\<2^k(t-r)\>^2}
              \big\|h_{j}(r,\cdot)\big\|_{L^{\infty}_{x'}}
              dr
          \Bigr\|_{L^1_t(\re_+)}
   \\
    \le &C \sum_{\bk j\ge  \widehat{0}} 2^{2j}  
            \Bigl\|
               \big\|h_{j}(t,\cdot)\big\|_{L^{\infty}_{x'}}
            \Bigr\|_{L^1_t(\re_+)}
    = C \big\|h\big\|_{L^1(\re_+;B^{2}_{\infty,1}(\re^{n-1}_{x'}))}.
\end{align*}
}
This concludes the proof of Proposition \ref{prop;L1'-L2'-estimate}.

\end{prf}

%
%
\sect{The Neumann boundary condition}
\subsection{The Neumann boundary condition and the boundary trace}\par
For the Neumann boundary condition, the boundary potential associated to 
the Green's function is slightly different from the one for the Dirichlet case.
However as we see the explicit form of the potential function, it simply 
adjust the order of the derivative and the estimate is very similar to the 
Dirichlet case. 
From \eqref{eqn;green-fn-N}, \eqref{eqn;potential-N} and \eqref{eqn;potential-besov-1},
{\allowdisplaybreaks
\begin{align}
  \overline{\Phi_m}&\underset{(x',\eta)}{*}\big(\eta^{-1} \Psi_N(t,x',\eta)\big) \nonumber\\
    =&  -c_n^2 \zeta_{m-1}(\eta)\underset{(\eta)}{*}\int_{\re^{n-1}}\int_{\re} 
                e^{it\t +ix'\cdot \xi'}  \widehat{\phi_m}(|\xi'|)
                \frac{ \t }{\sqrt{i\t+|\xi'|^2}}
                e^{-\sqrt{i\t+|\xi'|^2}\eta}
                d\xi'd\t  
        \nonumber\\
    &  - c_n^2 \phi_{m}(\eta)\underset{(\eta)}{*}   
             \int_{\re}\int_{\re^{n-1}}\int_{\re} 
                e^{it\t +ix'\cdot \xi'} \widehat{\zeta_{m}}(|\xi'|)
                \frac{ \t }{\sqrt{i\t+|\xi'|^2}}
                e^{-\sqrt{i\t+|\xi'|^2}\eta}
                d\xi'd\t. \label{eqn;potential-besov-21}
\end{align}
}
To estimate the Besov-norm of the solution, we use the Littlewood-Paley 
decomposition for direct sum type \eqref{eqn;direct-sum-L-P}. 
The estimates for the right hand side of 
\eqref{eqn;potential-besov-21} follows very similar manner 
to the case of the Dirichlet boundary condition. 
\algn{
   \int_0^{\infty} &\|\Del u(t)\|_{\dB^s_{p,1}(\re^{n})}dt  
 \\
     \le &C\Bigl\|\sum_{m\in \Z}2^{sm}\Bigl(\int_{\re_+} 
         \Bigl\| \Psi_N(t,x',\eta)
                \underset{(t,x')}{*} 
               \phi_m(x')\underset{(x')}{*} h(t,x')  
         \Bigr\|_{L^p(\re^{n-1}_{x'})}^p 
         d\eta \Bigr)^{1/p}\Bigr\|_{L^1_t(\re_+)}
 \\
        &+ C\Bigl\|
              \sum_{m\in \Z} 2^{sm}\Bigl(\int_{\re_+}
               \Bigl\|
                   \big(
                   \phi_{m}(\eta)\underset{(\eta)}{*}
                   \Psi_N(t,x',\eta)\big)
                   \underset{(t,x')}{*}  \\
         &\hskip4cm \times 
                   \sum_{k\in \Z} \sum_{j\le m}
                   \psi_k(t)\underset{(t)}{*}\phi_j(x')
                   \underset{(y')}{*} h(t,x') 
                  \Bigr\|_{L^p(\re^{n-1}_{x'})}^p  d\eta \Bigr)^{1/p}
               \Bigr\|_{L^1_t(\re_+)}
 \\
\equiv &\|P_1^N\|_{L^1_t(\re^+)}+\|P_2^N\|_{L^1_t(\re^+)}.
 \eqntag\label{eqn;potential-besov-19}
}

\subsection{Proof of Theorem \ref{thm;boundary-trace-N}}
We prove the sufficiency part of Theorem \ref{thm;boundary-trace-N}, which is reduced to prove Lemma \ref{lem;P1-P2-estimate-N}. We assume that 
$h$ satisfies \eqref{eqn;boundary_data_cond_besov-N}. 
For the Neumann case, we use the time-space splitting argument 
same as the Dirichlet case. 
We split each of two terms in the right hand side of  \eqref{eqn;potential-besov-19}
into time-dominated area and  space-dominated area. 
%
\begin{lem}\label{lem;P1-P2-estimate-N}
Let $1<p<\infty$ and $-1+\frac{1}{p}<s\le 0$.
There exists a constant $C>0$ such that $P_1^N$ defined in 
\eqref{eqn;potential-besov-19} satisfies 
\eq{\label{eqn;P1-est-N}
 \|P_1^N\|_{L^1_t(\re_+)}
 \le C\big(
       \big\|h \big\|_{ \dF^{1/2-1/2p}_{1,1}(\re_+;\dB^s_{p,1}(\re^{n-1}))}
      +\big\|h \big\|_{ L^1 (\re_+;\dB^{s+1-1/p}_{p,1}(\re^{n-1})) }
      \big).
}
Similarly $P_2^N$ defined in \eqref{eqn;potential-besov-19} satisfies 
the following estimate:
\eq{ \label{eqn;P2-est-N}
 \|P_2^N\|_{L^1_t(\re_+)}
 \le C\big(
       \big\|h \big\|_{ \dF^{1/2-1/2p}_{1,1}(\re_+;\dB^s_{p,1}(\re^{n-1}))}
      +\big\|h \big\|_{ L^1 (\re_+;\dB^{s+1-1/p}_{p,1}(\re^{n-1})) }
      \big).
}
\end{lem}

\begin{prf}{Lemma \ref{lem;P1-P2-estimate-N}}
Introducing 
\eqn{ 
  \Psi_{N,k,j}(t,x',\eta)
  \equiv 
     \int_{\re}\int_{\re^{n-1}}
     \Psi_N(t-s,x'-y',\eta)\psi_k(s)\phi_j(y')
      dy'ds,
}
we apply \eqref{eqn;space-time-splitting}
to $P_1^N$ and set $h_m\equiv \widetilde{\phi_m}\underset{(x')}{*}h$.
Similar to the estimates \eqref{eqn;potential-besov-10} and 
\eqref{eqn;L_1} in the Dirichlet boundary case, we obtain that 
\algn{
 & P_1^N \\
 \le& C\sum_{m\in \Z} 2^{sm}
       \Bigl( \int_{\re_+}
        \Big\|
         \sum_{k\in\Z}\sum_{2j\le k}
         \Psi_N(t,x',\eta)
         \underset{(t,x')}{*}        
         \psi_k(t)\underset{(t)}{*}
         \phi_m(x')\underset{(x')}{*} 
          h_j(t,x')  
        \Big\|_{L^p(\re^{n-1}_{x'})}^p
        d\eta
      \Bigr)^{1/p}  \\
  & +  C\sum_{m\in \Z} 2^{sm}
        \Bigl(\int_{\re_+}
         \Big\| 
          \sum_{k\in\Z}\sum_{2j>k}
          \Psi_N(t,x',\eta)\underset{(t,x')}{*}
          \psi_k(t)\underset{(t)}{*}
          \phi_m(x')\underset{(x')}{*} 
          h_j(t,x')  
         \Big\|_{L^p(\re^{n-1}_{x'})}^p
         d\eta
        \Bigr)^{1/p}
 \\
   \le &  C\sum_{m\in \Z} 2^{sm}
       \Bigl( \int_{\re_+}
        \Big\|
         \sum_{k\in\Z}
         \Psi_N(t,x',\eta)
         \underset{(t,x')}{*}        
         \psi_k(t)\underset{(t)}{*}
         \phi_m(x')\underset{(x')}{*} 
         h_m(t,x')  
        \Big\|_{L^p_{x'}}^p
         d\eta
      \Bigr)^{1/p}  \\
  &  +C\sum_{m\in \Z} 2^{sm}
       \Bigl( \int_{\re_+}
        \Big\|
         \sum_{k\in\Z} 
         \Psi_N(t,x',\eta)
         \underset{(t,x')}{*}        
         \psi_k(t)\underset{(t)}{*}
         \phi_m(x')\underset{(x')}{*} 
         h_m(t,x')  
        \Big\|_{L^p_{x'}}^p
        d\eta
      \Bigr)^{1/p} 
 \\ 
   \le &C\sum_{m\in\Z}  2^{sm}
       \Bigl(\int_{\re_+}
        \Bigl\{ \sum_{k\ge 2m} 
             \int_{\re_+}
             \Big\|\Psi_{N,k,m}(t-s,x',\eta)\Big\|_{L^1_{x'}}
             \Big\|\psi_k\underset{(s)}{*} h_m(s,x')  \Big\|_{L^p_{x'}}
             ds
        \Bigr\}^p  
         d\eta
     \Bigr)^{1/p}  \\
    &+C\sum_{m\in\Z}  2^{sm}
       \Bigl(\int_{\re_+}
        \Bigl\{ \sum_{k\le 2m} 
             \int_{\re_+}
             \big\|\Psi_{N,k,m}(t-s,x',\eta)
             \big\|_{L^1_{x'}}
             \big\| h_m(s,x')  \big\|_{L^p_{x'}}
             ds
        \Bigr\}^p  
        d\eta
     \Bigr)^{1/p} \\
   \equiv&  L_1+ L_2.
   \eqntag\label{eqn;P1-N}
}
Then similar to the Dirichlet boundary condition, we invoke
Lemma \ref{lem;pt-orthogonal-2} in Section \ref{sec;almostorthogonal}
 and it implies 
\begin{equation}   \label{eqn;Ortho-N}
  \Big\| \Psi_{N,k,m}(t,\cdot)\Big \|_{L^1_{x'}}
  \le
   \left\{
   \begin{aligned}
      & C_n 2^{\frac{k}{2}}e^{-2^{\frac{k}{2}}\eta } 
            \frac{2^k}{\<2^kt\>^{2}}{,} 
         &{ k\ge 2m},\\
      & C_n 2^{\frac{k}{2}}e^{-2^{m}\eta}
           \frac{2^k}{\<2^kt\>^{2}}{,}
         &{k<2m.}
  \end{aligned}
  \right.
\end{equation}
Applying \eqref{eqn;Ortho-N} to \eqref{eqn;P1-N}, by 
the similar manner as the proof of \eqref{eqn;L_1-2} 
in Lemma \ref{lem;P1-P2-estimate}, we have for the 
first term of the right hand side that 
\algn{ 
 \big\|  L_1\|_{L^1_t(\re_+)}  
   \le &C\left\|
         \sum_{m\in \Z} 2^{sm}
         \left(\int_{\re_+}
          \bigg\{ \sum_{k\ge2m}
          \big( 2^{\frac{k}{2}}e^{-2^{\frac{k}{2}}\eta } \big)
          \int_{\re}
              \frac{2^k}{\<2^k(t-s)\>^2}
              \Big\|
                 \psi_k\underset{(t)}{*}h_m(s,\cdot)
              \Big\|_{L^p_{x'}} 
              ds 
          \bigg\}^p         
           d\eta
          \right)^{1/p}
         \right\|_{L^1_t(\re_+)} 
\\
  = &C\left\|
         \sum_{m\in \Z} 2^{sm}
         \sum_{k\ge 2m}
           2^{\frac{k}{2}}
          \int_{\re}
              \frac{2^k}{\<2^k(t-s)\>^2}
              \Big\|
                 \psi_k\underset{(t)}{*}h_m(s,\cdot)
              \Big\|_{L^p_{x'}} 
              ds 
          \left(
          \int_{\re_+}       
            e^{-p2^{\frac{k}{2}}\eta } d\eta
          \right)^{1/p}
         \right\|_{L^1_t(\re_+)} 
\\
  = &C\left\|
         \sum_{m\in \Z} 2^{sm}
         \sum_{k\in\Z}
           2^{\frac{k}{2}-\frac{k}{2p}}
          \int_{\re}
              \frac{2^k}{\<2^k(t-s)\>^2}
              \Big\|
                 \psi_k\underset{(t)}{*}h_m(s,\cdot)
              \Big\|_{L^p_{x'}} 
              ds 
         \right\|_{L^1_t(\re_+)} 
\\
  \le &C \sum_{m\in \Z} 2^{sm}
         \sum_{k\in\Z}
           2^{\frac{k}{2}-\frac{k}{2p}}
          \left\|
              \Big\|
                 \psi_k\underset{(t)}{*}h_m(s,\cdot)
              \Big\|_{L^p_{x'}} 
         \right\|_{L^1_t(\re_+)} 
 \\
  \le&C\big\|
          h 
       \big\|_{\dF^{1/2-1/2p}_{1,1}(\re_+;\dB^s_{p,1}(\re^{n-1}_{x'}))}.
       \eqntag\label{eqn;L_1-2-N}
}

The term $ L_2$ stands for the space-dominated part and we set 
$h_m\equiv \widetilde{\phi_m}*h$. 
Observing $\phi_j*\phi_m\simeq
\phi_{m}$  ($|j-m|\le1$) 
we apply \eqref{eqn;Ortho-N} to \eqref{eqn;P1-N}
to see that
\algn{ 
 \big\| L_2 \|_{L^1_t(\re_+)} 
  \le& C\Bigg\| 
         \sum_{m\in \Z} 2^{sm}
         \bigg(\int_{\re_+}
          \bigg\{ \sum_{k<2m}
               \big(2^{\frac{k}{2}} e^{-2^{m}\eta} \big)
             \int_{\re}
              \frac{2^k}{\<2^k(t-s)\>^2}
              \big\| h_m(s,x')  \big\|_{L^p_{x'}}
              ds 
          \bigg\}^p  
           d\eta
          \bigg)^{1/p}
         \Bigg\|_{L^1_t(\re_+)} 
  \\
    \le& C\Bigg\|  
           \sum_{m\in \Z} 2^{sm}
           \sum_{k<2m}
                2^{\frac{k}{2}}
              \int_{\re}
              \frac{2^k}{\<2^k(t-s)\>^2}
              \big\| h_m(s,x')  \big\|_{L^p_{x'}}
               ds  
          \left(\int_{\re_+}
             e^{-p2^{m}\eta} 
           d\eta
          \right)^{1/p}
         \Bigg\|_{L^1_t(\re_+)} 
  \\
    \le& C\Bigg\| 
           \sum_{m\in \Z} 2^{sm}
           \sum_{k<2m}
                2^{\frac{k}{2}} 2^{-\frac{m}{p}}
              \int_{\re}
              \frac{2^k}{\<2^k(t-s)\>^2}
              \big\| h_m(s,x')  \big\|_{L^p_{x'}}
               ds  
         \Bigg\|_{L^1_t(\re_+)} 
   \\
    \le& C\sum_{m\in \Z} 2^{(s+1-\frac{1}{p})m}
          \sum_{k<2m}    2^{\frac{k}{2}-m}
         \Bigg\|
               \big\| h_m(s,x')  \big\|_{L^p_{x'}}        
         \Bigg\|_{L^1_t(\re_+)} 
  \\
   \le& C\Bigg\| \sum_{m\in \Z}  2^{(s+1-\frac{1}{p})m}
               \big\| h_m(s,x')  \big\|_{L^p_{x'}}        
         \Bigg\|_{L^1_t(\re_+)} 
  \\
    = & C\big\|h\big\|_{ L^1 (\re_+;\dB^{s+1-1/p}_{p,1}(\re^{n-1}_{x'}) )}.
    \eqntag\label{eqn;L_2-2-N}
}

To see \eqref{eqn;P2-est-N},  we split in a similar way to 
\eqref{eqn;potential-besov-11} again and obtain that 
\algn{
 &P_2^N \\
 \le &C\sum_{m\in \Z} 2^{sm}
      \Big\|         
         \sum_{k\in\Z}\sum_{2j\le k}
         \Big\|
             \big(
             \phi_m(\eta)\underset{(\eta)}{*} 
             \Psi_N(t,x',\eta)\big)
         \underset{(t,x')}{*}        
         \psi_k(t)\underset{(t)}{*}\phi_j(x')
         \underset{(x')}{*}h(t,x')  
        \Big\|_{L^p(\re^{n-1}_{x'})}
      \Big\|_{L^p(\re^+_{\eta})} \\
  &+C\sum_{m\in \Z} 2^{sm}
       \Bigl\|
           \sum_{k\in\Z}\sum_{2j>k}
          \Big\|
              \big( 
              \phi_m(\eta)\underset{(\eta)}{*} 
              \Psi_N(t,x',\eta)\big)
                    \underset{(t,x')}{*}  
          \psi_k(t)\underset{(t)}{*}\phi_j(x')
          \underset{(x')}{*}h(t,x')  
         \Big\|_{L^p(\re^{n-1}_{x'})}
       \Bigr\|_{L^p(\re^+_{\eta})}
\\
    \le &C\sum_{m\in\Z} 2^{sm} 
       \left(\int_{\re_+}
        \left\{ \sum_{k\in\Z}\sum_{2j\le k}
             \int_{\re_+}
             \Big\| \phi_m(\eta)\underset{(\eta)}{*}
                    \Psi_{N,k,j}(t-s,x',\eta)
             \Big\|_{L^1_{x'}}
             \Big\|\psi_k\underset{(s)}{*} h_j(s,x')  \Big\|_{L^p_{x'}}
             ds
        \right\}^p  d\eta
      \right)^{1/p} \\
     &+C\sum_{m\in\Z}  2^{sm}
       \left(\int_{\re_+}
        \bigg( \sum_{j\in\Z}\sum_{k<2j}
             \int_{\re_+}
             \Big\| \phi_m(\eta)\underset{(\eta)}{*}
                    \Psi_{N,k,j}(t-s,x',\eta)
             \Big\|_{L^1_{x'}}
             \Big\| h_j(s)  \Big\|_{L^p_{x'}}
             ds
        \bigg)^p  d\eta
       \right)^{1/p} \\
 \equiv & M_1 + M_2 .
 \eqntag \label{eqn;potential-besov-31} 
}

We use the almost orthogonal estimate \eqref{eqn;Ortho-N}
between $\psi_m$ and $\Psi_{N,k,j}$.  
From \eqref{eqn;potential-besov-31}, 
\algn{  
 \big\| & M_1\big\|_{L^1_t(\re_+)}  
\\
  \le&C\Bigg\|\sum_{m\in\Z} 2^{sm} \bigg(\int_{\re_+}            
             \bigg\{ 
              \sum_{k\in\Z}\sum_{2j\le k}               
              \Big(
                \frac{C_N 2^{-|\frac{k}{2}-m|}}
                     {\<{2^{\min(\frac{k}{2},m)}} |\eta|\>^N}
              2^{\frac{k}{2}}           
              \int_{\re}
              \frac{2^k}{\<2^k(t-s)\>^2}
              \Big\|
                \psi_k \underset{(s)}{*} h_j(s,\cdot) 
              \Big\|_{L^p_{x'}}
             ds 
             \bigg\}^pd\eta \bigg)^{1/p} 
     \Bigg\|_{L^1_t(\re_+)} 
\\
  \le&C\Bigg\|\sum_{m\in\Z}  2^{sm}
              \sum_{k\in\Z}\sum_{2j\le k}               
              \Big(
                 2^{-|\frac{k}{2}-m|}  
                 2^{\frac{k}{2}}           
              \int_{\re}
              \frac{2^k}{\<2^k(t-s)\>^2}
              \Big\|
                \psi_k \underset{(s)}{*} h_j(s,\cdot) 
              \Big\|_{L^p_{x'}}
             ds  \\
       &\hskip5cm\times
             \, {2^{-\frac{1}{p}\min(\frac{k}{2},m))}}
             \bigg(\int_{\re_+}
                \frac{1}
                     {\<|\teta|\>^{pN}}
              d\teta \bigg)^{1/p} 
     \Bigg\|_{L^1_t(\re_+)} 
\\
  \le&C\Bigg\|  
          \sum_{k\in \Z}  2^{\frac{k}{2}-\frac{k}{2p}}
          \sum_{2j\le k}  2^{\frac{s}{2}k}
    \big(\sum_{m'\in\Z}2^{-|m'|}  2^{sm'}
                { 2^{\max(0,-\frac{m'}{p})} }\big) 
              \Big\|
                \psi_k \underset{(s)}{*} h_j(s,\cdot) 
              \Big\|_{L^p_{x'}} 
     \Bigg\|_{L^1_t(\re_+)} 
\\
 = &C \|h\|_{\dF^{\frac12-\frac{1}{2p}}_{1,1}
                (\re_+;\dB^s_{p,1}(\re^{n-1}_{x'})) }
    \eqntag\label{eqn;M_1-2-N}
} 
under the condition $-1+1/p<s\le 0$.
Setting 
$
 h_j(t)\equiv \phi_j\underset{(x')}{*}h(t)
$
in \eqref{eqn;potential-besov-11}
and changing the order of summation between $j$ and $k$, 
we have similar to the estimate for $M_1$
by the Minkowski and the Hausdorff-Young inequalities with 
\eqref{eqn;chi} that 
\algn{ 
   \big\| & M_2\big\|_{L^1_t(\re_+)}  
\\
  \le&C\Bigg\|\sum_{m\in\Z} 2^{sm} \bigg(\int_{\re_+}            
             \bigg\{ 
              \sum_{j\in\Z}\sum_{k\le 2j}               
              \Big(
                \frac{C_N 2^{-|j-m|}}
                     {\<{2^j} |\eta|\>^N}
              2^{\frac{k}{2}}           
              \int_{\re}
              \frac{2^k}{\<2^k(t-s)\>^2}
              \big\|
                 h_j(s,\cdot) 
              \big\|_{L^p_{x'}}
             ds 
             \bigg\}^pd\eta \bigg)^{1/p} 
     \Bigg\|_{L^1_t(\re_+)} 
\\
  \le&C\Bigg\|\sum_{m\in\Z} 2^{sm}   
              \sum_{j\in\Z}\sum_{k\le 2j}                
              \Big(
                 2^{-|j-m|}  
                 2^{\frac{k}{2}}           
              \int_{\re}
              \frac{2^k}{\<2^k(t-s)\>^2}
              \big\|
                 h_j(s,\cdot) 
              \big\|_{L^p_{x'}}
             ds 
             \,{2^{-\frac{j}{p}}}
             \bigg(\int_{\re_+}
                \frac{1}
                     {\<|\teta|\>^{pN}}
              d\teta \bigg)^{1/p} 
     \Bigg\|_{L^1_t(\re_+)} 
\\
  \le&C\Bigg\|  
         \sum_{j\in\Z}   2^{sj+j-\frac{j}{p}}
    \big(\sum_{k\le 2j}  2^{\frac{k}{2}-j} \big)
    \big(\sum_{m'\in\Z}2^{-|m'|}   2^{sm'}\big)
              \big\|
                 h_j(s,\cdot) 
              \big\|_{L^p_{x'}} 
     \Bigg\|_{L^1_t(\re_+)} 
\\
 = &C \|h\|_{L^1(\re_+;\dB^{s+1-\frac{1}{p}}_{p,1}(\re^{n-1}_{x'})) }
   \eqntag\label{eqn;M_2-3-N}
}  
under the condition {$|s|< 1$}.
From \eqref{eqn;M_1-2-N} and \eqref{eqn;M_2-3-N}, we conclude the proof of 
\eqref{eqn;P2-est-N}.
\end{prf}

%
%
\subsection{Proof of Theorem \ref{thm;Naumann-BC-Lp-trace}}\par
The proof of Theorem \ref{thm;Naumann-BC-Lp-trace} can be done by 
a similar way to the case of the Dirichlet boundary condition.
Indeed we obtain a similar estimate to \eqref{eqn;P1-est-N}:
{\allowdisplaybreaks
\begin{align}
   \int_0^{\infty}& \|\Del u(t)\|_{L^p(\re^{n}_+)}dt  \nonumber\\
        = &\Bigl\|\Bigl\|\Bigl(\int_{\re_+}
                  \Big|                   
                  \int_0^{\infty}\int_{\re^{n-1}_{x'}}
                   \Psi_N(t-s,x'-y',\eta) h(s,y')ds dy'
                   \Big|^pd\eta
                  \Bigr)^{1/p} 
            \Bigr\|_{L^p(\re^{n-1}_{x'})}
           \Bigr\|_{L^1_t(\re_+)}
\nonumber\\
   =  &  \Bigl\|\Bigl\| \Bigl(\sum_{\ell\in \Z} {2^{-\ell} }
                 \big\| 
                 \Psi_N(t,x',\eta) \big|_{\eta\simeq 2^{-\ell}}
                 \underset{(t,x')}{*}  
                  h(t,x')  
                 \big\|_{L^p(\re^{n-1}_{x'})}^p 
                \Bigr)^{1/p}
       \Bigr\|_{L^p(\re^{n-1}_{x'})}
           \Bigr\|_{L^1_t(\re_+)},\label{eqn;boudary-Lp-1-N}
\end{align}
}
then we split the last term into two terms:
{\allowdisplaybreaks
\begin{align}
 \big\| \Psi_{N}(t,x',\eta)
          \underset{(t,x')}{*}  
          & h(t,x')  
 \big\|_{L^p(\re^{n-1}_{x'})}\nonumber\\
  \le & \big\| \Psi_{N}(t,x',\eta)
               \underset{(t,x')}{*}  
        \sum_{k\in\Z}\sum_{\widehat{0}\le j\le k/2}
       \psi_k(t)*\phi_j(x')*h(t,x')  
      \big\|_{L^p(\re^{n-1})} \nonumber\\
    &+\big\|  \Psi_{N}(t,x',\eta)
              \underset{(t,x')}{*}  
        \sum_{j\ge \widehat{0}}\sum_{k<2j} 
        \psi_k(t)*\phi_j(x')*h(t,x')  
      \big\|_{L^p(\re^{n-1})}\nonumber\\
   \le & \sum_{k\in\Z}
          \int_{\re} 
              \sup_{j\in\Z} 
             \Big\|
               \Psi_{N,k,j}(t-r,\cdot,\eta)
             \Big\|_{L^1_{x'}}
             \Big\|
                   \psi_k\underset{(t)}{*}
                   \sum_{\widehat{0} \le j\le k/2}\phi_j*h(r,\cdot) 
             \Big\|_{L^p_{x'}} \nonumber\\
    & +\sum_{j\ge \widehat{0}}
          \int_{\re} 
              \sup_{k\in\Z}
             \Big\|
                \Psi_{N,k,j}(t-r,\cdot,\eta)
             \Big\|_{L^1_{x'}}
             \Big\|
              \sum_{k<2j}\psi_k*h(s,\cdot)
             \Big\|_{L^p_{x'}} ds    \nonumber\\
  \equiv &\;  N^1_{\ell}+  N^2_{\ell}.
\label{eqn;boudary-Lp-2-N}
\end{align}
}
From \eqref{eqn;boudary-Lp-1-N} and 
\eqref{eqn;boudary-Lp-2-N},  
we reduce the estimate for the following:
{\allowdisplaybreaks
\begin{align*}
   \int_0^{\infty} \|\Del u(t)\|_{L^p(\re^{n}_+)}dt 
   =  &  \Bigl\| \Bigl(\sum_{\ell\in \Z} {2^{-\ell} }
                 \big\| 
                  \Psi_{N}(t,x',\eta)
                  \underset{(t,x')}{*}  
                  h(t,x')  
                 \big\|_{L^p(\re^{n-1})}^p 
                 \Big|_{\eta\in I_{-\ell}}
                \Bigr)^{1/p}
      \Bigr\|_{L^1_t(\re_+)}
\\
  \le & \Big\|
           \| \{N^1_{\ell}\}_{\ell} \|_{\ell^{-1/p}_p}
        \Big\|_{L^1_t(\re_+)}
       +\Big\|
           \| \{N^2_{\ell}\}_{\ell} \|_{\ell^{-1/p}_p}
        \Big\|_{L^1_t(\re_+)}.
  \end{align*}
}
The proof of Theorem \ref{thm;Naumann-BC-Lp-trace} is now reduced 
to show the following proposition. 

\begin{prop}\label{prop;P1-P2-estimate}\par\noindent
Let $1\le p\le \infty$.  For $ N^1_{\ell}$ and $ N^2_{\ell}$ 
given by \eqref{eqn;boudary-Lp-2-N}, the following estimates holds.
\begin{align}
 &\Big\|\| \{N^1_{\ell}\}_{\ell}\|_{\ell^{-1/p}_p}\Big\|_{L^1_t(\re_+)}
  \le C\big\|h \big\|_{ \dot{F}^{1/2-1/2p}_{1,1}(\re_+;L^p(\re^{n-1}))}, 
       \notag \\
 &\Big\|\| \{N^2_{\ell}\}_{\ell} \|_{\ell^{-1/p}_p}\Big\|_{L^1_t(\re_+)}
  \le C\big\|h \big\|_{ L^1 (\re_+;B^{1-1/p}_{p,1}(\re^{n-1})) }. 
       \notag
\end{align}
\end{prop}

The proof of Proposition \ref{prop;P1-P2-estimate} can be shown
in a similar way to Proposition \ref{prop;L1'-L2'-estimate}.

\sect{Almost orthogonality}\label{sec;almostorthogonal}
In this section we prove the almost orthogonality 
\eqref{eqn;crucial-potential-orthogonarity-org} 
(or \eqref{eqn;Ortho-D0} in Section \ref{sec;4})
between the boundary potential $\Psi_D$ for the Dirichlet boundary 
case and $\Psi_N$ for the Neumann boundary case with 
the time and space  Littlewood-Paley decomposition $\{\psi_k\}_{k\in Z}$ 
and  $\{\phi_j\}_{j\in Z}$. 
The difficulty is that the fundamental solution $\Psi_{D}(t,x',\eta)$ 
and $\Psi_{N}(t,x',\eta)$ 
made by heat kernel is time and space convolution. 

\subsection{The Dirichlet potential case}\par

\begin{lem}[Crucial potential orthogonality]
\label{lem;pt-orthogonal-1}
For  $k,j,\ell\in \Z$ let $\{\psi_k(t)\}_{k\in \Z}$ and 
$\{\phi_j(x)\}_{j\in \Z_+}$ 
be the time and the space Littlewood-Paley dyadic decomposition and
let $\Psi_D(t,x',\eta)$ be the boundary potential defined in 
\eqref{eqn;potential-D}. 
Set
\eq{\label{eqn;potential-Dkj}
  \Psi_{D,k,j}(t,x',\eta)
   \equiv \int_{\re}\int_{\re^{n-1}}
        \Psi_D(t-s,x'-y',\eta)\psi_k(s)\phi_j(y')dy'ds,
  }
for $\eta>0$.
Then there exists a constant $C_n>0$ depending only on the dimension $n$ 
satisfying 
\eq{ \label{eqn;crucial-potential-orthogonarity-D}
  \|\Psi_{D,k,j}(t,\cdot,\eta)\|_{L^1_{x'}}
  \le
  \left\{
   \begin{aligned}
      &C_n 2^{k}\big(1+(2^{\frac{k}{2}}\eta)^{n+2}\big)
       e^{-2^{\frac{k}{2}-1}\eta } \frac{2^k}{\<2^kt\>^{2}},
        & {k\ge 2j},\\
      &C_n 2^{k} \big(1+(2^{j}\eta)^{n+2}\big)
       e^{-2^{j-1}\eta } \frac{2^k}{\<2^kt\>^{2}},
       &{ k<2j}.
   \end{aligned}
  \right.
  }
\end{lem}

\begin{prf}{Lemma \ref{lem;pt-orthogonal-1}}
We consider a time-like estimate $k> 2j$. 
Taking $\zeta'$-space cut-off, and using the change of variables 
$\t=2^k\s$, $\xi'=2^j\zeta'$ we have  
{\allowdisplaybreaks
\begin{align*}
   \big\|& \Psi_{D,k,j}(t,\cdot,\eta)\big\|_{L^1_{x'}} \\
  =&\Bigl\|
     c_{n+1}\int_{\re}\int_{\re^{n-1}}
         e^{it\t+ix'\cdot \xi'}
            \t 
                \exp\big(-\sqrt{i\t+|\xi'|^2}\eta\big)
                \widehat{\psi}(2^{-k}\t)\widehat{\phi}(2^{-j}\xi')
     d\xi'd\t   
    \Bigr\|_{L^1_{x'}}
 \\
  =&\Bigl\|
     c_{n+1}\int_{\re}\int_{\re^{n-1}}
         e^{i2^kt\s+i2^jx'\cdot \zeta'}
           \s   2^k  
                \exp\bigl(-\sqrt{2^{k}i\s+2^{2j}|\zeta'|^2}
                    \eta
                    \bigr)
          \widehat{\psi}(\s)\widehat{\phi}(\zeta')
        2^{(n-1)j}d\zeta'\cdot 2^kd\s   
    \Bigr\|_{L^1_{x'}}
 \\
   =&  c_{n+1} 2^k
     \Bigl\| \int_{\re}\int_{\re^{n-1}}
             e^{i2^kt\s+i2^jx'\cdot \zeta'} 
             \s \exp\bigl(-{\bk  2^\frac{k}{2}\eta } 
                          \sqrt{i\s+2^{2j-k}|\zeta'|^2}
                    \bigr)
             \widehat{\psi}(\s)\widehat{\phi}(\zeta')
             2^{(n-1)j}d\zeta'\cdot 2^kd\s   
         \Bigr\|_{L^1_{x'}}
\\
       =& c_{n+1} 2^k
        \Bigl\|\int_{\re}\int_{\re^{n-1}}
            e^{i2^kt\s+iy'\cdot \zeta'}
             \exp\bigl(-{\bk  2^\frac{k}{2}\eta } 
                          \sqrt{i\s+2^{2j-k}|\zeta'|^2}
                 \bigr)
             \s \widehat{\psi}(\s)\widehat{\phi}(\zeta')
             d\zeta'\cdot 2^kd\s   
         \Bigr\|_{L^1_{y'}}
 \eqntag \label{eqn;Dirichlet-ortho-1}
\end{align*}
}
by  setting $x'=2^{-j}y'$ in the last equality.
Using the identity 
\eq{ \label{eqn;exp-intgralbyparts}
   e^{i(2^kt\sg+y'\cdot \zeta')}
   =-\frac{1}{i2^kt}\frac{1}{|y'|^2}\Del_{\zeta'}\pt_{\s}
    e^{i(2^kt\sg+y'\cdot \zeta')}, 
}
and integrating by parts with respect to $\sg$ and $\zeta'$, 
setting $p(\sg,\zeta';2^{2j-k})\equiv\sqrt{i\s+2^{2j-k}|\zeta'|^2}$, we proceed 
{\allowdisplaybreaks
\begin{align*}
  &\|\Psi_{D,k,j}(t,\cdot)
     \|_{ L^1_{x'}(B_{2^{j}}^c) } \\
   =&c_{n+1}2^k
     \Bigl\|
       \int_{\re}\int_{\re^{n-1}}
         \frac{1}{i2^kt}
         \left(-\frac{1}{|y'|^2}\right)e^{i(2^kt\s+ y'\cdot \zeta')}\\
    &\phantom{2^{k-2\ell}\|c_n\int_{\re}\int} 
         \times
         \Del_{\zeta'}\frac{\pt}{\pt \s}
          \bigl[\exp\left(-{\bk  2^\frac{k}{2}\eta } 
                          p(\sg,\zeta';2^{2j-k})
                    \right)
          \widehat{\psi}(\s)
          \widehat{\phi}(\zeta')
        2^{k}\s\bigr] d\zeta' d\s   
    \Bigr\|_{L^1_{y'}(B_1^c)}.
\end{align*}
}
Here
{\allowdisplaybreaks 
\begin{align*}
 &\frac{\pt}{\pt \s}
         \bigl[ \exp\left(
                   -{ 2^\frac{k}{2}\eta }\,
                    p(\sg,\zeta';2^{2j-k})
                    \right)
               \widehat{\psi}(\s)2^{k}\s 
               \widehat{\phi}(\zeta')
          \bigr]\\
     =&\exp\big(-  2^\frac{k}{2}\eta \, p(\sg,\zeta';2^{2j-k})\big) 
      \biggl\{ 
                  \frac{-{  2^\frac{k}{2}\eta }i}
                       {2p(\sg,\zeta';2^{2j-k})}
                 \widehat{\psi}(\s)2^{k}\s 
                   +\widehat{\psi}'(\s)2^{k}\s 
                    + \widehat{\psi}(\s)2^{k}
          \biggr\}\widehat{\phi}(\zeta').
\end{align*}
}
In order to take the second derivative for space,  
setting $r=|\zeta'|$ and using the relation 
$\Del_{\zeta'}=\pt_r^2+{\frac{n-2}{r}}\pt_r$, 
we have by
$\widehat{\phi}(\zeta')=\widehat{\phi}(|\zeta'|)$ that 

\algn{
& \Del_{\zeta'} \frac{\pt}{\pt \s}
         \Bigl[ \exp\Bigl(-{  2^\frac{k}{2}\eta } \,
                    \sqrt{i\s+ 2^{2j-k}|\zeta'|^2}
               \Bigr)
               \widehat{\psi}(\s)2^{k}\s 
               \widehat{\phi}(\zeta')
          \Bigr]\\
=&\Del_{\zeta'}\exp\big(-{  2^\frac{k}{2}\eta } 
                   \, p(\sg,\zeta';2^{2j-k})\big) \\
     &\qquad \times
           \Bigg\{ 
                 \frac{-{  2^\frac{k}{2}\eta }i}
                      {2p(\sg,\zeta';2^{2j-k})}
                  \widehat{\psi}(\s)2^{k}\s 
                   +\widehat{\psi}'(\s)2^{k}\s 
                    + \widehat{\psi}(\s)2^{k}
           \Bigg\}\widehat{\phi}(\zeta')
\\
   =&\Bigl(\pt_r+{\frac{n-2}{r}}\Bigr)
           \exp\big(-{  2^\frac{k}{2}\eta } 
           \, p(\sg,\zeta';2^{2j-k})\big)\\
        &\qquad \times
           \Bigl[
             \Bigl(-{  2^\frac{k}{2}\eta } 
                    \frac{2^{2j-k}|\zeta'|}{p(\sg,\zeta';2^{2j-k})}
                    \widehat{\phi}(\zeta')
                    +\widehat{\phi}'(\zeta')
             \Bigr)
             \Bigr( 
                 \frac{-{  2^\frac{k}{2}\eta }i}
                      {2p(\sg,\zeta';2^{2j-k})}
                 \widehat{\psi}(\s)2^{k}\s 
                + \widehat{\psi}'(\s)2^{k}\s 
                + \widehat{\psi}(\s)2^{k}
             \Bigr)\\
       &\qquad\qquad
             +{  2^\frac{k}{2}\eta } 
                 \frac{2^{2j-k} i |\zeta'|}
                      {2p(\sg,\zeta';2^{2j-k})^{3}}
                 \widehat{\psi}(\s)2^{k}\s
                 \widehat{\phi}(\zeta')
          \Bigr]
   \\
  =&\exp\big(-{ 2^\frac{k}{2}\eta } p(\sg,\zeta';2^{k-2j})\big)\\
      &\quad \times\Bigl[
         \Bigl(-{ 2^\frac{k}{2}\eta } 
          \frac{2^{2j-k}|\zeta'|}{p(\sg,\zeta';2^{2j-k})}
         \Bigr)\Bigl\{
             \Bigl(-{  2^\frac{k}{2}\eta } 
                    \frac{2^{2j-k}|\zeta'|}{p(\sg,\zeta';2^{2j-k})}
                    \widehat{\phi}(\zeta')
                    +\widehat{\phi}'(\zeta')
             \Bigr)\\
        &\qquad \times
             \Bigl( 
                 \frac{-{ 2^\frac{k}{2}\eta }i}
                      {2p(\sg,\zeta';2^{2j-k})}
                 \widehat{\psi}(\s)2^{k}\s 
                + \widehat{\psi}'(\s)2^{k}\s 
                + \widehat{\psi}(\s)2^{k}
             \Bigr)
             +{ 2^\frac{k}{2}\eta } 
                 \frac{2^{2j-k} i |\zeta'|}
                      {2p(\sg,\zeta';2^{2j-k})^{3}}
                 \widehat{\psi}(\s)2^{k}\s
                 \widehat{\phi}(\zeta')
          \Bigr\}\\ 
      &\qquad\quad 
       +\Bigl(-{ 2^\frac{k}{2}\eta } 
                  \frac{2^{2j-k}i\sigma}
                       {p(\sg,\zeta';2^{2j-k})^3}
                    \widehat{\phi}(\zeta')
                   -{  2^\frac{k}{2}\eta } 
                   \frac{2^{2j-k}|\zeta'|}{p(\sg,\zeta';2^{2j-k})}
                    \widehat{\phi}'(\zeta')
                   +\widehat{\phi}''(\zeta')
             \Bigr)\\
      &\qquad \qquad \quad\times
             \Bigl(-\frac{-{ 2^\frac{k}{2}\eta }i}{2p(\sg,\zeta';2^{2j-k})}
                 \widehat{\psi}(\s)2^{k}\s 
                + \widehat{\psi}'(\s)2^{k}\s 
                + \widehat{\psi}(\s)2^{k}
             \Bigr)\\
        &\qquad \qquad +
             \Bigl(-{ 2^\frac{k}{2}\eta } 
             \frac{2^{2j-k}|\zeta'|}{p(\sg,\zeta';2^{2j-k})}
                    \widehat{\phi}(\zeta')
                    +\widehat{\phi}'(\zeta')
             \Bigr)
             \Bigl({ 2^\frac{k}{2}\eta } 
              \frac{2^{2j-k}i|\zeta'|}{2p(\sg,\zeta';2^{2j-k})^3}
                 \widehat{\psi}(\s)2^{k}\s
             \Bigr)\\
        &   +{ 2^\frac{k}{2}\eta } 
                  \frac{-3 (2^{2j-k})^2 i |\zeta'|^2
                        +  2^{2j-k} i p(\sg,\zeta';2^{2j-k})^2}
                       {2p(\sg,\zeta';2^{2j-k})^{5}}
                 \widehat{\psi}(\s)2^{k}\s
                 \widehat{\phi}(\zeta')             
                  +{  2^\frac{k}{2}\eta } 
                  \frac{ 2^{2j-k}i |\zeta'|}
                       {2p(\sg,\zeta';2^{2j-k})^{3}}
                 \widehat{\psi}(\s)2^{k}\s
                 \widehat{\phi}'(\zeta')
          \Bigr]\\
   &+{(n-2)}\exp\big(-{ 2^\frac{k}{2}\eta } p(\sg,\zeta';2^{2j-k})\big)\\
       &\qquad \times
        \Bigl[
             \Bigl(-{ 2^\frac{k}{2}\eta } 
                    \frac{2^{2j-k}}{p(\sg,\zeta';2^{2j-k})}
                    \widehat{\phi}(\zeta')
                    +\frac{\widehat{\phi}'(\zeta')}{|\zeta'|}
             \Bigr)
             \Bigl( 
                 \frac{-{ 2^\frac{k}{2}\eta }i}
                      {2p(\sg,\zeta';2^{2j-k})}
                 \widehat{\psi}(\s)2^{k}\s 
                + \widehat{\psi}'(\s)2^{k}\s 
                + \widehat{\psi}(\s)2^{k}
             \Bigr)\\
        &\qquad\qquad
             +{ 2^\frac{k}{2}\eta } 
                 \frac{ 2^{2j-k} i }
                      {2p(\sg,\zeta';2^{2j-k})^{3}}
                 \widehat{\psi}(\s)2^{k}\s
                 \widehat{\phi}(\zeta')
          \Bigr].
}
Summarizing the results, all terms have
$$
\exp\bigl(-{ 2^\frac{k}{2}\eta } \sqrt{i\s+2^{2j-k}|\zeta'|^2}\bigr)
 =\exp\big(-{ 2^\frac{k}{2}\eta }\,  p(\sg,\zeta';2^{2j-k})\big)
$$
namely $\widehat{\psi}(\sg)$, $\widehat{\phi}(\zeta')$
or their derivatives and the order of derivatives are 
the order of partial derivatives with respect to $\sg$ and $\zeta'$.  
And ${2^{\frac{k}{2}}\eta}$ arises the same order of partial derivatives at most 
and arises one time at least. 
The function in denominator $p(\sg,\zeta';2^{2j-k})=\sqrt{i\s+2^{2j-k}|\zeta'|^2}$ 
is estimated as
\eq{\label{eqn;power-bound-D-T}
 2^{-1/2}\le
  \Big|p(\sg,\zeta',2^{2j-k})\Big|
  = (\sg^2+2^{2(2j-k)}|\zeta'|^4)^{1/4} 
  \le 20^{1/4},
}
thanks to the cut-off functions $\widehat{\psi}(\sg)$ 
and  $\widehat{\phi}(\zeta')$ or its derivative. 
Therefore if we differentiate it $n+2$ times, 
then the terms involving $p(\sg,\zeta';2^{2j-k})$ are estimated from below by $2^{n+2}$. 
 Moreover $p(\sg,\zeta';2^{2j-k})$ which arises by each derivative contains the surplus 
scale parameter $2^{2j-k}$ and it is estimated from above by $1$ because of 
the restriction $k\ge 2j$.  
Thus for the functions
{\allowdisplaybreaks
\begin{align*}
  k_{k,j}(\sg,\zeta',\eta)
   \equiv& \exp\left(-{ 2^\frac{k}{2}\eta }\,
                    p(\sg,\zeta';2^{2j-k})
                    \right)
          \widehat{\psi}(\s)2^{k}\s
          \widehat{\phi}(\zeta'),   \\
  K_{k,j}^{n+2}(\sg,\zeta',\eta,2^{2j-k})
   \equiv&\exp\left({2^\frac{k}{2}\eta }\, p(\sg,\zeta';2^{2j-k})
               \right)
          D^{n+2}_{\sg,\zeta'} k_{k,j}(\sg,\zeta',\eta),
\end{align*}
}
where we set
$$
  D^{n+2}_{\sg,\zeta'}
   \equiv(1-\Del_{\zeta'})^{\frac{n}{2}}
            \bigl(1-\frac{\pt^2}{\pt \sg^2}\bigr),  
$$
under the condition $k\ge 2j$ namely $2^{2j-k}\le 1$, we obtain 
the following estimate: For $\eta>0$ 
\algn{
 &\|\Psi_{D,k,j}(t,\cdot,\eta)
  \|_{ L^1_{x'}} \\
   =&c_{n+1}2^k
     \Bigl\|\Bigl(\frac{1}{\<y'\>^2}\Bigr)^{\frac{n}{2}}
            \frac{2^k}{\<2^kt\>^2}
            \int_{\re}\int_{\re^{n-1}}  
            e^{i(2^kt\sg+y'\cdot \zeta')}
            \exp\Bigl(-{ 2^\frac{k}{2}\eta }
                    p(\sg,\zeta';2^{2j-k})
                    \Bigr)\\
     &\hskip10cm\times 
        K_{k,j}^{n+2}(\sg,\zeta',\eta,2^{2j-k})
        d\zeta'd\sg
     \Bigr\|_{L^1_{y'}}\\
  =&c_{n+1}2^k \exp(-{2^{\frac{k}{2}-1}\eta}) \frac{2^k}{\<2^kt\>^2}\\
     &\quad\times 
     \bigg\|\Bigl(\frac{1}{|y'|^2}\Bigr)^{\frac{n}{2}}           
        \int_{\re}\int_{\re^{n-1}}  
            e^{i(2^kt\sg+y'\cdot \zeta')}
            e^{\big(-{2^\frac{k}{2}\eta}
                      \big(p(\sg,\zeta';2^{2j-k})-1/2\big)
                \big)
              }
             K_{k,j}^{n+2}(\sg,\zeta',\eta, 2^{2j-k})
            d\zeta'd\sg
     \bigg\|_{L^1_{y'}}\\
 \le &C2^k
      \exp(-{2^{\frac{k}{2}-1}\eta}) \frac{2^k}{\<2^kt\>^2}\\
     &\quad\times          
        \int_{\re}\int_{\re^{n-1}}
            e^{\big(-{ 2^\frac{k}{2}\eta }
                      \big(p(\sg,\zeta';2^{2j-k})-1/2\big)
                \big)
              }
           \Big|
                K_{k,j}^{n+2} (\sg,\zeta',\eta,2^{2j-k})
           \Big|
         d\zeta'd\sg
   \\
 \le&C2^k\exp(-2^{\frac{k}{2}-1}\eta)
    \big(1+C( 2^{\frac{k}{2}}\eta)^{n+2}\big)\frac{2^k}{\<2^kt\>^2},
}
arranging the suffixes we obtain \eqref{eqn;crucial-potential-orthogonarity-D}.

%
%

Next we consider the case $k\le 2j$. 
Setting $\eta=2^{-\ell}$ and using the change of variables 
         $\t=2^k\s$, $\xi'=2^j\zeta'$
\begin{align*}
   \|\Psi_{D,k,j}&(t,\cdot)\|_{L^1_{x'}} \\
  =&\Bigl\|
     c_{n+1}\int_{\re}\int_{\re^{n-1}}
         e^{it\t+ix'\cdot \xi'}
         {  \t\exp\big(-\sqrt{i\t+|\xi'|^2}\eta\big)
         }
         \widehat{\psi}(2^{-k}\t)\widehat{\phi}(2^{-j}\xi')
     d\xi'd\t   
    \Bigr\|_{L^1_{x'}}
\\
  =&\Bigl\|
     c_{n+1}\int_{\re}\int_{\re^{n-1}}
         e^{i2^kt\s+i2^jx'\cdot \zeta'}
           \s2^{k}
                \exp\bigl(-\sqrt{2^{k}i\s+2^{2j}|\zeta'|^2}
                    \eta
                    \bigr)
          \widehat{\psi}(\s)\widehat{\phi}(\zeta')
        2^{(n-1)j}d\zeta'\cdot 2^kd\s   
    \Bigr\|_{L^1_{x'}}
\\
     =& 2^{k}\Bigl\|
        c_{n+1}\int_{\re}\int_{\re^{n-1}}
            e^{i2^kt\s+iy'\cdot \zeta'}
             \exp\bigl(-{2^{j}\eta}\sqrt{2^{k-2j}i\s +|\zeta'|^2}
                 \bigr)
             \s \widehat{\psi}(\s)\widehat{\phi}(\zeta')
             d\zeta'\cdot 2^kd\s   
         \Bigr\|_{L^1_{y'}},
\end{align*}
where we set setting $x'=2^{-j}y'$ in the last equality.
Noting the identity \eqref{eqn;exp-intgralbyparts} and 
integrating by parts with respect to  $\sg$ and $\zeta'$, 
\begin{align*}
  &\|\Psi_{D,k,j}(t,\cdot)
     \|_{ L^1_{x'}(B_{2^{j}}^c) } \\
   =&2^k
     \Bigl\|
      c_{n+1}\int_{\re}\int_{\re^{n-1}}
         \frac{1}{i2^kt}
         \left(-\frac{1}{|y'|^2}\right)e^{i(2^kt\s+ y'\cdot \zeta')}\\
    &\phantom{2^{k-2\ell}\|c_n\int_{\re}\int} 
         \times
         \Del_{\zeta'}\frac{\pt}{\pt \s}
          \Bigl[\exp\big(-{2^{j}\eta}
                    \sqrt{\big(2^{k-2j}\s i+|\zeta'|^2\big)}
                    \big)
          \widehat{\psi}(\s)
          \widehat{\phi}(\zeta')
        2^{k}\s\Bigr] d\zeta' d\s   
    \Bigr\|_{L^1_{y'}(B_1^c)}.
\end{align*}
Set $q(\sg,\zeta';2^{k-2j})=\sqrt{2^{k-2j}i\s +|\zeta'|^2}$.
Then we have
\begin{align*}
 \frac{\pt}{\pt \s}&
         \Bigl[ \exp\left(-{2^{j}\eta}\,
                    q(\sg,\zeta';2^{k-2j})
                    \right)
               \widehat{\psi}(\s)2^{k}\s 
               \widehat{\phi}(\zeta')
          \Bigr]
  \\
     =&\exp\big(-{2^{j}\eta}\, q(\sg,\zeta';2^{k-2j})\big)
           \Bigr\{-{2^{j}\eta} \frac{2^{k-2j}i}{2q(\sg,\zeta';2^{k-2j})}
                 \widehat{\psi}(\s)2^{k}\s 
                   +\widehat{\psi}'(\s)2^{k}\s 
                    + \widehat{\psi}(\s)2^{k}
          \Bigr\}\widehat{\phi}(\zeta').
\end{align*}
In order to take second derivative with respect to space,  
setting $r=|\zeta'|$ and using the relation 
$\Del_{\zeta'}=\pt_r^2+\frac{n-2}{r}\pt_r$, we have 
\begin{align*}
 &\Del_{\zeta'} \frac{\pt}{\pt \s}
         \Bigl[ \exp\bigl(-{2^{j}\eta}
                    \sqrt{2^{k-2j}i\s+|\zeta'|^2}
               \bigr)
               \widehat{\psi}(\s)2^{k}\s 
               \widehat{\phi}(\zeta')
          \Bigr]
 \\
    =&\Del_{\zeta'}\exp\big(-{2^{j}\eta} q(\sg,\zeta';2^{k-2j})\big)
           \Bigg\{-{2^{j}\eta} \frac{2^{k-2j}i}{2q(\sg,\zeta';2^{k-2j})}
                  \widehat{\psi}(\s)2^{k}\s 
                   +\widehat{\psi}'(\s)2^{k}\s 
                    + \widehat{\psi}(\s)2^{k}
           \Bigg\}\widehat{\phi}(\zeta')
  \\
     =&\bigl(\pt_r+\frac{n-2}{r}\bigr)
           \exp\big(-{2^{j}\eta}\, q(\sg,\zeta';2^{k-2j})\bigr)\\
     &\quad \times
           \Bigl[
             \Bigl(-\frac{{2^{j}\eta} |\zeta'|}{q(\sg,\zeta';2^{k-2j})}
                    \widehat{\phi}(\zeta')
                    +\widehat{\phi}'(\zeta')
             \Bigr)
             \Bigl(-{2^{j}\eta} \frac{2^{k-2j}i}{2q(\sg,\zeta';2^{k-2j})}
                 \widehat{\psi}(\s)2^{k}\s 
                + \widehat{\psi}'(\s)2^{k}\s 
                + \widehat{\psi}(\s)2^{k}
             \Bigr)\\
      &\quad\qquad
             +{2^{j}\eta} \frac{2^{k-2\ell}i |\zeta'|}
                       {2q(\sg,\zeta';2^{k-2j})^{3}}
                 \widehat{\psi}(\s)2^{k}\s
                 \widehat{\phi}(\zeta')
          \Bigr]
 \\
  =&\exp\big(-{2^{j}\eta}\, q(\sg,\zeta';2^{k-2j})\big)\\
     &\quad \times\Bigl[
         \Bigl(-\frac{{2^{j}\eta}|\zeta'|}{q(\sg,\zeta';2^{k-2j})}
         \Bigr)
         \Bigl\{
            \Bigl(-\frac{{2^{j}\eta}|\zeta'|}{q(\sg,\zeta';2^{k-2j})}
                    \widehat{\phi}(\zeta')
                    +\widehat{\phi}'(\zeta')
            \Bigr) \\
      &\times
             \Bigl(-{2^{j}\eta} \frac{2^{k-2j}i}{2q(\sg,\zeta';2^{k-2j})}
                 \widehat{\psi}(\s)2^{k}\s 
                + \widehat{\psi}'(\s)2^{k}\s 
                + \widehat{\psi}(\s)2^{k}
             \Bigr) 
              +{2^{j}\eta} \frac{2^{k-2\ell}i |\zeta'|}
                       {2q(\sg,\zeta';2^{2j-k})^{3}}
                 \widehat{\psi}(\s)2^{k}\s
                 \widehat{\phi}(\zeta')
          \Bigr\}\\ 
      &\qquad\quad 
       +
             \Bigl(-{2^{j}\eta} \frac{2^{k-2j}i \sg}
                                   {q(\sg,\zeta';2^{k-2j})^3}
                    \widehat{\phi}(\zeta')
                   - \frac{{2^{j}\eta}|\zeta'|}{q(\sg,\zeta';2^{k-2j})}
                    \widehat{\phi}'(\zeta')
                   +\widehat{\phi}''(\zeta')
             \Bigr)\\
      &\qquad \qquad \quad\times
             \Bigl(-{2^{j}\eta} \frac{2^{k-2j}i}{2q(\sg,\zeta';2^{k-2j})}
                 \widehat{\psi}(\s)2^{k}\s 
                + \widehat{\psi}'(\s)2^{k}\s 
                + \widehat{\psi}(\s)2^{k}
             \Bigr)\\
        &\qquad\qquad +
             \Bigl(- \frac{{2^{j}\eta}|\zeta'|}{q(\sg,\zeta';2^{k-2j})}
                    \widehat{\phi}(\zeta')
                    +\widehat{\phi}'(\zeta')
             \Bigr)
             \Bigl(\frac{2^j\eta 2^{k-2j}i|\zeta'|}{2q(\sg,\zeta';2^{2j-k})^3}
                 \widehat{\psi}(\s)2^{k}\s
             \Bigr)\\
      &\qquad\qquad\quad
             +{2^{j}\eta} \frac{-3  |\zeta'|^2
                        +  q(\sg,\zeta';2^{k-2j})^2}
                       {2q(\sg,\zeta';2^{k-2j})^{5}}
                 \widehat{\psi}(\s)2^{k}\s
                 \widehat{\phi}(\zeta')
             +{2^{j}\eta} \frac{2^{k-2j}i |\zeta'|}
                       {2q(\sg,\zeta';2^{k-2j})^{3}}
                 \widehat{\psi}(\s)2^{k}\s
                 \widehat{\phi}'(\zeta')
          \Bigr]\\
   &+(n-2)\exp\big(-{2^{j}\eta}\, q(\sg,\zeta';2^{k-2j})\big)\\
     &\quad \times
           \Bigl[
             \Bigl(- \frac{2^{j}\eta}{q(\sg,\zeta';2^{k-2j})}
                    \widehat{\phi}(\zeta')
                    +\frac{\widehat{\phi}'(\zeta')}{|\zeta'|}
             \Bigr)
             \Bigl(-{2^{j}\eta} \frac{2^{k-2j}i}{2q(\sg,\zeta';2^{k-2j})}
                 \widehat{\psi}(\s)2^{k}\s 
                + \widehat{\psi}'(\s)2^{k}\s 
                + \widehat{\psi}(\s)2^{k}
             \bigg)\\
      &\qquad
             +{2^{j}\eta} \frac{2^{k-2j}i}
                       {2q(\sg,\zeta';2^{k-2j})^{3}}
                 \widehat{\psi}(\s)2^{k}\s
                 \widehat{\phi}(\zeta')
          \Bigr]. 
\end{align*} 
As the case for $k> 2j$ before, we may summarize that all terms are including 
$$
\exp\big(-{2^{j}\eta} \sqrt{2^{k-2j}i\s+|\zeta'|^2} \big)
 =\exp\big(-{2^{j}\eta}\, q(\sg,\zeta';2^{k-2j})\big)
$$
$\widehat{\psi}(\sg)$, $\widehat{\phi}(\zeta')$, or their derivatives, 
and the order of derivatives is the same as the order of partial derivatives of 
$\sg$ and $\zeta'$. 
And ${2^{j}\eta}$ is multiplied by the same order of the partial derivatives at most, 
and one time at least. 
In the denominator the function $q(\sg,\zeta';2^{k-2j})=\sqrt{2^{k-2j}i\s+|\zeta'|^2}$ 
is estimated from below by $2^{-1}$ thanks to the cut-off function $\widehat{\phi}(\zeta')$ 
or its derivative when $j\ge1$. Thus if we derive it $n+2$ times, the terms involving 
$q(\sg,\zeta';2^{k-2j})$ are estimated from below 
by $2^{n+2}$. 
 Moreover $q(\sg,\zeta';2^{k-2j})$ which arises by derivative contains surplus 
scale parameter $2^{k-2j}$ which is estimated from above by $1$ by  $k<2j$. 
Setting 
\algn{
  h_{k,j}(\sg,\zeta',\eta)
   \equiv & \exp\bigl(-{2^{j}\eta} 
                    \sqrt{2^{k-2j}i\s +|\zeta'|^2}
                    \bigr)
          2^{k}\s\widehat{\psi}(\s)
          \widehat{\phi}(\zeta')       
  \\
   H_{\ell,k,j}^{n+2}(\sg,\zeta',\eta, 2^{k-2j})
   \equiv &\exp\bigl(2^j\eta
                    \sqrt{2^{k-2j}i\s +|\zeta'|^2}
               \bigr)
          D^{n+2}_{\sg,\zeta'} h_{\ell,k,j}(\sg,\zeta',\eta).
}
and
$$
  D^{n+2}_{\sg,\zeta'}
   \equiv(1-\Del_{\zeta'})^{\frac{n}{2}}
            \bigl(1-\frac{\pt^2}{\pt \sg^2}\bigr),  
$$
we obtain
\algn{
 &\| \Psi_{D,k,j}(t,\cdot)
  \|_{ L^1_{x'}} \\
   =& c_{n+1} 2^{k}
     \Bigl\|\Bigl(\frac{1}{\<y'\>^2}\Bigr)^{\frac{n}{2}}
            \frac{2^k}{\<2^kt\>^2}
        \int_{\re}\int_{\re^{n-1}}  
            e^{i(2^kt\sg+y'\cdot \zeta')}
            (1-\Del_{\zeta'})^{\frac{n}{2}}
            \Bigl(1-\frac{\pt^2}{\pt \sg^2}\Bigr)
            h_{k,j}(\sg,\zeta',\eta)
            d\zeta'd\sg
     \Bigr\|_{L^1_{y'}}
\\
   =&c_{n+1}2^{k}
     \Bigl\|\Bigl(\frac{1}{|y'|^2}\Bigr)^{\frac{n}{2}}
            \frac{2^k}{\<2^kt\>^2}
         \int_{\re}\int_{\re^{n-1}}  
            e^{i(2^kt\sg+y'\cdot \zeta')}
            \exp\Bigl(-{2^{j}\eta}
                    \sqrt{2^{k-2j}i\s+|\zeta'|^2}
                    \Bigr)\\
     &\hskip10cm\times         
          H_{\ell,k,j}^{n+2}(\sg,\zeta',\eta, 2^{k-2j})
            d\zeta'd\sg
     \Bigr\|_{L^1_{y'}}
\\
    =&c_{n+1}  {2^{k}\eta^2}
       \exp(-{2^{j-1}\eta}) \frac{2^k}{\<2^kt\>^2}\\
     &\quad\times 
     \Bigl\|\Bigl(\frac{1}{|y'|^2}\Bigr)^{\frac{n}{2}}           
        \int_{\re}\int_{\re^{n-1}}  
            e^{i(2^kt\sg+y'\cdot \zeta')}
            e^{\bigl(-{2^{j}\eta}
                      \big(q(\sg,\zeta';2^{k-2j})-1/2\big)
                \bigr)
              }
            H_{\ell,k,j}^{n+2}(\sg,\zeta',\eta, 2^{k-2j})
            d\zeta'd\sg
     \Bigr\|_{L^1_{y'}}
\\
\le & C2^{k}\exp(- 2^{j-1}\eta) \frac{2^k}{\<2^kt\>^2}
        \int_{\re}\int_{\re^{n-1}}
           e^{\left(-{2^{j}\eta}
                     \big(q(\sg,\zeta';2^{k-2j})-1/2\big)
               \right)
             }
          \bigl|
           H_{\ell,k,j}^{n+2}(\sg,\zeta',\eta, 2^{k-2j})
          \bigr|
        d\zeta'd\sg
  \\
 \le&C 2^{k}\exp(- 2^{j-1}\eta)
    \big(1+C(2^{j}\eta))^{n+2}\big)\frac{2^k}{\<2^kt\>^2},
}
which completes the proof of Lemma \ref{lem;pt-orthogonal-1}. 
\end{prf}

\subsection{The second orthogonal estimate}
For the estimating the term $P_2^D$ in \eqref{eqn;potential-besov-11}
it involving an $\eta$-convolution between 
the potential $\Psi_D$ and $\phi_m(\eta)$.
We show the following second orthogonal estimate:
\vskip2mm
\begin{lem}[Potential orthogonarity 2]
\label{lem;pt-orthogonal-D2}
Let $k,j,m\in \Z$ and assume $j\le m+1$. Let $\Psi_D(t,x',\eta)$ be
the potential of the solution for the Dirichlet data 
defined by \eqref{eqn;potential-D} and let 
$\{\psi_k(t)\}_{k\in \Z}$ and $\{\phi_j(x)\}_{j\in \Z_+}$
be a spatial and time Littlewood-Paley decomposition. Let
$\Psi_{D,k,j}(t,x',\eta)$ be defined by \eqref{eqn;potential-Dkj}. 
Then for any $N\in \Nt$, there exists a constant $C_N>0$ such that 
  \eq{ \label{eqn;crucial-potential-orthogonarity-D2}
  \| (\phi_m\underset{(\eta)}{*}\Psi_{D,k,j} )  (t,\cdot,\eta)\|_{L^1_{x'}}
  \le
  \left\{
   \begin{aligned}
      &C_N 2^{k}\frac{ 2^{-|\frac{k}{2}-m|}}
                     {\<2^{\min(\frac{k}2, m)}\eta\>^N}\frac{2^k}{\<2^kt\>^{2}},
        &{k\ge 2j},\\
      &C_N 2^{k}\frac{2^{-|j-m|}}
                     {\<{2^j}\eta\>^N}\frac{2^k}{\<2^kt\>^{2}},
        &{k<2j}.
   \end{aligned}
  \right.
  }
\end{lem}
\vskip2mm

\begin{prf}{Lemma \ref{lem;pt-orthogonal-D2}}
 Assuming   $2j< k$, we show the time-dominated estimate. 
Changing $\t=2^k\s$, $\xi'=2^j\zeta'$, as we have seen in
\eqref{eqn;Dirichlet-ortho-1}, it also follows  that 
\eqn{
 \spl{
   \|&\phi_m\underset{(\eta)}{*}\Psi_{D,k,j}(t,\cdot,\eta)\|_{L^1_{x'}} \\
  =&\left\|
     c_{n+1}\int_{\re}\int_{\re^{n-1}}
         e^{it\t+ix'\cdot \xi'}
            i\t 
            \Big(\phi_m\underset{(\eta)}{*}
                 \exp\big(-\sqrt{i\t+|\xi'|^2}\eta\big)
            \Big)\widehat{\psi}(2^{-k}\t)\widehat{\phi}(2^{-j}\xi')
     d\xi'd\t   
    \right\|_{L^1_{x'}}\\
  =& c_{n+1}\left\|
    \int_{\re}\int_{\re^{n-1}}
         e^{i2^kt\s+i2^jx'\cdot \zeta'}
          i 2^k\s  \right.\\
   &\hskip 30mm \times\left.
         \Big(\phi_m\underset{(\eta)}{*}
              \exp\left(-\sqrt{2^{k}i\s+2^{2j}|\zeta'|^2}\eta
                  \right)
         \Big)
         \widehat{\psi}(\s)\widehat{\phi}(\zeta')
        2^{(n-1)j}d\zeta'\cdot 2^kd\s   
    \right\|_{L^1_{x'}}\\
   =& c_{n+1}\left\|
     \int_{\re}\int_{\re^{n-1}}
         e^{i2^kt\s+i2^jx'\cdot \zeta'} 
         i\s  2^k \right.\\
     &\hskip 30mm \times\left.
            \Big(\phi_m\underset{(\eta)}{*}
                  \exp\big(-  2^\frac{k}{2}\eta  
                      \sqrt{i\s+2^{2j-k}|\zeta'|^2}\big)
            \Big)
            \widehat{\psi}(\s)\widehat{\phi}(\zeta')
            2^{(n-1)j}d\zeta'\cdot 2^kd\s   
         \right\|_{L^1_{x'}}\\
     =&  c_{n+1} 2^k 
        \left\|\int_{\re}\int_{\re^{n-1}}
            e^{i2^kt\s+iy'\cdot \zeta'}
             \Big(\phi_m\underset{(\eta)}{*}
              \exp\big(- 2^{\frac{k}{2}}\eta  
                          \sqrt{i\s+2^{2j-k}|\zeta'|^2}
                 \big)
             \Big)
             \s \widehat{\psi}(\s)\widehat{\phi}(\zeta')
             d\zeta'\cdot 2^kd\s   
         \right\|_{L^1_{y'}}.
   }
}
Then using 
\eqn{
   e^{i(2^kt\sg+y'\cdot \zeta')}
   =-\frac{1}{i2^kt}\frac{1}{|y'|^2}\Del_{\zeta'}\pt_{\s}
    e^{i(2^kt\sg+y'\cdot \zeta')}
}
we apply the integration by parts by 
$\sg$, $\zeta'$ as before.  Then the following lemma is 
a crucial step.

\vskip4mm
\noindent
\begin{lem}\label{lem;almost-orthogonal-4} 
Let $k,j,m\in \Z$ and assume {$j< m+1$} and
$2^{-1}\le \sg, |\zeta'|\le 2$. 
For $\al=\al_1+\al_2$, $0\le \al_1\le 1$, $0\le \al_2\le n$,
the following estimates hold:
\algn{\label{eqn;almost-orthogonal-4} \eqntag
 &\bigg| \phi_m\underset{(\eta)}{*} 
       \Big(\sum_{|\al|\le n+2} C_{n,\al}
            \pt^{\al_1}_{\sg}
            \pt^{\al_2}_{\zeta'}   
            \exp\big(- 2^{\frac{k}{2}\eta}
            \sqrt{i\sg+2^{2j-k}|\zeta'|^2} \big)
       \Big)
  \bigg| 
  \le  \frac{ C_N 2^{-|\frac{k}{2}-m|}}
            { \<2^{\min(\frac{k}2,m)}\eta\>^{N} },  
           \quad  
           {
           k\ge 2j,
           }\\
  \label{eqn;almost-orthogonal-5} \eqntag
  &\bigg| \phi_m\underset{(\eta)}{*} 
       \Big(\sum_{|\al|\le n+2} C_{n,\al}
            \pt^{\al_1}_{\sg}
            \pt^{\al_2}_{\zeta'}   
           \exp\big(- 2^{j}\eta
                   \sqrt{i2^{k-2j}\sg+|\zeta'|^2} \big)
       \Big)
  \bigg| 
  \le  \frac{ C_N 2^{-|j-m|}}
            {\<2^{j}\eta\>^{N} },  
           \quad  {k<2j. }  
 }
\end{lem}
\vskip3mm

\begin{prf}{Lemma \ref{lem;almost-orthogonal-4}}
We first show the lemma for the case of {$k\ge 2j $}.
For $p(\sg,\zeta',2^{2j-k})\equiv\sqrt{i\s+2^{2j-k}|\zeta'|^2}$, 
let
\algn{
 P_{\al}(\sg,\zeta', 2^{k/2}\eta)
 =&\exp\big( 2^{\frac{k}{2}}\eta \sqrt{i\s+2^{2j-k}|\zeta'|^2}\big)
   \sum_{|\al|\le n+2}
    C_{n,\al}
     \pt^{\al_1}_{\sg}
     \pt^{\al_2}_{\zeta'}
     \exp\big(- 2^{\frac{k}{2}}\eta \sqrt{i\s+2^{2j-k}|\zeta'|^2}\big). 
      \eqntag
      \label{eqn;almost-orthogonal-6}
}
Then $P_{\al}$ is a polynomial of $ 2^{\frac{k}{2}}\eta$ and 
$|p(\sg,\zeta',2^{2j-k})|^{-1}=
|\sqrt{i\sg+2^{2j-k}|\zeta'|^2}|^{-1}$
and by \eqref{eqn;power-bound-D-T},
it holds that
\eq{\label{eqn;polynomial-bound-D-T}
 | P_{\al}(\sg,\zeta',  2^{k/2}\eta)|
   \le C\big(1+( 2^{\frac{k}{2}} \eta)^{n+2}\big).
 }
We also set
\eq{\label{eqn;change-of-parameters}
 \left\{
 \spl{
  &\bar\th =  2^{\frac{k}{2}}\sqrt{i\s+2^{2j-k}|\zeta'|^2} \th
           =  2^{\frac{k}{2}}p \th, 
            \\
  &\bar\eta = 2^{\frac{k}{2}}\sqrt{i\s+2^{2j-k}|\zeta'|^2}\eta
            = 2^{\frac{k}{2}}p  \eta.
}\right.
}

\noindent
$\<$Step 1$\>$: The case {$k\ge 2m\ge 2j$}.
Using \eqref{eqn;power-bound-D-T},
\algn{
\bigg| \int_{\re_+}&\phi_m(\eta-\th) 
             \Big(
               P_{\al}(\sg,\zeta', 2^{k/2}\th)
              \exp\big(-  2^{\frac{k}{2}}\th
                       \sqrt{i\s+2^{2j-k}|\zeta'|^2} \big)
             \Big)
             d\th
 \bigg| \\
 =&\bigg|\int_{\re_+} 
             2^m\phi
             \big(2^m  2^{-\frac{k}{2}}\sqrt{i\s+2^{2j-k}|\zeta'|^2}^{-1}
                  (\bar\eta-\bar\th)
             \big) 
          \bar{P}_{\al}(\bar\th)
          \exp\big(-\bar\th \big)
          2^{-\frac{k}{2}}\sqrt{i\s+2^{2j-k}|\zeta'|^2}^{-1}
          d\bar\th
   \bigg|
 \\
 \le&\, C2^{m-\frac{k}{2}}
        \int_{\re} 
          \big|\phi
             \big(2^m  2^{-\frac{k}{2}}
                  \Big|\sqrt{i\s+2^{2j-k}|\zeta'|^2}\Big|^{-1}
                  (\bar\eta-\bar\th)\big)
          \big| 
          \bar{P}_{\al}(\bar\th)
          \exp\big(-|\bar\th| \big)
          d\bar\th 
 \\
\le&\, C2^{m-\frac{k}{2}}
     \int_{|\bar\th|> \frac12|\bar\eta|} 
          \frac{C_N}
               {\big\<2^m  2^{-\frac{k}{2}}
                  \big|\sqrt{i\s+2^{2j-k}|\zeta'|^2}\big|^{-1}
                  (\bar\eta-\bar\th)\big\>^N}
          \bar{P}_{\al}(\bar\th)
          \exp\big(-|\bar\th| \big)
          d\bar\th  \\
   &+C2^{m-\frac{k}{2}}
     \int_{|\bar\th|\le \frac12|\bar\eta|} 
          \frac{C_N}
          {\big\<2^m  2^{-\frac{k}{2} }
                  \big|\sqrt{i\s+2^{2j-k}|\zeta'|^2}\big|^{-1}
                  (\bar\eta-\bar\th)\big\>^N}
          \tilde{P}_{\al}(\bar\th)
          \exp\big(-|\bar\th| \big)
          d\bar\th \\
  \equiv& I+II,
  \label{eqn;rev-almost-ortho-1} \eqntag
}
where we set $ \bar{P}_{\al}(\bar\th)$ as a polynomial of 
$2^{\frac{k}{2}}\bar\th$ and has the following estimate
form \eqref{eqn;polynomial-bound-D-T};
\eq{   \label{eqn;rev-almost-ortho-1.5}
 \big|\bar{P}_{\al}(\bar\th)\big|
 \le C_{\al}\big(1+ |2^{\frac{k}{2}}\bar\th|\big)^{n+2}
}
under the restrictions $2^{-1}\le \sg<2$, 
$2^{-1}<|\zeta'|\le 2$ and {$k\ge 2j$}.

Then the first term of the right hand side of 
\eqref{eqn;rev-almost-ortho-1} is estimated by using 
 \eqref{eqn;rev-almost-ortho-1.5} 
as follows: 
\algn{  
 I \le&\,C2^{m-\frac{k}{2}}
      \exp\big(-\frac14 |\bar\eta| \big)  
      \int_{ |\bar\th|> \frac12|\bar\eta| } 
        \frac{C_N}
             {\big\<2^m  2^{-\frac{k}{2}}
              \big|\sqrt{i\s+2^{2j-k}|\zeta'|^2}\big|^{-1}
              (\bar\eta-\bar\th)\big\>^N}     
         \bar{P}_{\al}(\bar\th) \exp\big(-\frac12 |\bar\th| \big)  
      d\bar\th  
  \\
 \le&\,C2^{m-\frac{k}{2}} 
      \frac{C_N}{\<{ 2^{\frac{k}{2}}
                   \big|\sqrt{i\s+2^{2j-k}|\zeta'|^2}\big| 
                   }
                  \eta\>^N}  
      \int_{|\bar\th|> \frac12|\bar\eta|}     
      \bar{P}_{\al}(\bar\th) \exp\big(-\frac12 |\bar\th| \big)  
      d\bar\th  
  \\
  \le& \frac{C_N 2^{m-\frac{k}{2}} }
            { \big\< {2^{\frac{k}{2}}}\eta\big\>^N }  
  \le  \frac{ C_N  2^{ m-{\frac{k}{2}} } } 
            { \big\< 2^m \eta\big\>^N } 
 \label{eqn;rev-almost-ortho-2}\eqntag  
}
by {$k\ge 2m$}.
For the second term in \eqref{eqn;rev-almost-ortho-1},
we note that $|\bar\th|< \frac12|\bar\eta|$ implies
$|\bar\eta-\bar\th|\ge |\bar\eta|-|\bar\th|
  \ge |\bar\eta|-\frac12|\bar\eta|=\frac12|\bar\eta|$,
and it follows that 
\algn{
 II\le&\,C2^{m-\frac{k}{2}}
      \int_{|\bar\th|\le \frac12|\bar\eta|} 
       \frac{C_N}
            {\big\<2^m 2^{-\frac{k}{2}}
       \big|\sqrt{i\s+2^{2j-k}|\zeta'|^2}\big|^{-1}
       (\bar\eta-\bar\th)\big\>^N} 
        \bar{P}_{\al}(\bar\th) \exp\big(-|\bar\th| \big)   
        d\bar\th 
\\
 \le& \frac{ C_N 2^{m-\frac{k}{2}} }
           { \big\<2^{m-1} 2^{-\frac{k}{2}}
                  \big|\sqrt{i\s+2^{2j-k}|\zeta'|^2}\big|^{-1}
                 |\bar\eta|\big\>^N 
           }
      \int_{\re}       
          \bar{P}_{\al}(\bar\th)  \exp\big(-|\bar\th| \big)  
          d\bar\th
 \\
 \le& \frac{ C_N 2^{m-\frac{k}{2}} }
           { \big\<2^{m-1} \eta\big\>^N }.
  \label{eqn;rev-almost-ortho-3}\eqntag
}
Hence by \eqref{eqn;rev-almost-ortho-1}, 
\eqref{eqn;rev-almost-ortho-2} and \eqref{eqn;rev-almost-ortho-3},
 we obtain 
\eq{ \label{eqn;rev-almost-ortho-4}
 \bigg|  \int_{\re_+}\phi_m(\eta-\th) 
             \Big(
               P_{\al}(\t,\xi', 2^{k/2}\th)
              \exp\big(- 2^{\frac{k}{2}}\th
                       \sqrt{i\s+2^{2j-k}|\zeta'|^2} \big)
             \Big)
             d\th
 \bigg|
 \le 
 \frac{ C_N 2^{m-\frac{k}{2}} }
      { \big\<2^{m-1} \eta\big\>^N }.
 }

\noindent
$\<$ Step 2$\>$: When {$2j\le  k\le 2m$},
we note that 
$$
 \int_{\re}\phi_m(\th)d\th=0
$$
and 
\algn{
\bigg| \int_{\re_+}&\phi_m(\th) 
             \Big(
               P_{\al}(\t,\xi', 2^{k/2}(\eta-\th))
              \exp\big(- 2^{\frac{k}{2}}(\eta-\th)
                       \sqrt{i\s+2^{2j-k}|\zeta'|^2} \big)
             \Big)
             d\th
 \bigg| 
\\
 =&\bigg| \int_{\re_+}\phi_m(\th) 
             \Big(
               P_{\al}(\t,\xi', 2^{k/2}(\eta-\th))
              \exp\big(- 2^{\frac{k}{2}}(\eta-\th)
                       \sqrt{i\s+2^{2j-k}|\zeta'|^2} \big) \\
     &\hskip4cm
     - P_{\al}(\t,\xi', 2^{k/2}\eta)
              \exp\big(-2^{\frac{k}{2}}\eta
                       \sqrt{i\s+2^{2j-k}|\zeta'|^2} \big)
             \Big)
             d\th
 \bigg| 
\\
 =&\bigg| \int_{\re_+}\phi_m(\th) 
              \Big(\int_0^1\frac{d}{d\nu}
               P_{\al}(\t,\xi', 2^{k/2}(\eta-\nu\th))
              \exp\big(- 2^{\frac{k}{2}}(\eta-\nu\th)
                       \sqrt{i\s+2^{2j-k}|\zeta'|^2} \big) 
              d\nu\Big)
              d\th
 \bigg| 
 \\
 =&\bigg| \int_0^1 \int_{\re_+}\phi_m(\th)             
                  2^{\frac{k}{2}}\th 
                \sqrt{i\s+2^{2j-k}|\zeta'|^2} 
                 \Big(
                P_{\al}(\t,\xi', 2^{k/2}(\eta-\nu\th))
              -\pt_{\mu}P_{\al}(\t,\xi',\mu)
                \Big|_{\mu=2^{\frac{k}{2}}(\eta-\nu\th)}
              \Big) \\
      &\hskip3cm
           \times \exp\big(- 2^{\frac{k}{2}}(\eta-\nu\th)
                       \sqrt{i\s+2^{2j-k}|\zeta'|^2} \big) 
             d\th
             d\nu
 \bigg| 
\\
 = &\bigg|\int_0^1 
            \int_{\re} 
              2^m  \phi\big(2^m  2^{-\frac{k}{2}} 
                \sqrt{i\s+2^{2j-k}|\zeta'|^2}^{-1}\bar\th\big)
             \,\bar\th\\
     &\hskip25mm\times 
             \widetilde{P_{\al}}(\t,\zeta',(\bar\eta-\nu \bar\th))  
             \exp\big(-(\bar{\eta}-\nu\bar{\th}) 
                 \big)                          
              2^{-\frac{k}{2}}
             \sqrt{i\s+2^{2j-k}|\zeta'|^2}^{-1}  d\bar\th
             d\nu
     \bigg|
\\
   \le& \int_0^1 
        \int_{|\bar\th|> \frac12|\bar\eta|} 
          \big|\phi\big(2^m  2^{ -\frac{k}{2} }   
               \sqrt{i\s+2^{2j-k}|\zeta'|^2}^{-1}\bar\th\big)
          \big|
              2^{ m - \frac{k}{2} }  
             \sqrt{i\s+2^{2j-k}|\zeta'|^2}^{-1}|\bar\th| \\
        &\hskip30mm\times
             \widetilde{P_{\al}}(\t,\zeta',(\bar\eta-\nu \bar\th))
             \exp\big(-|\bar{\eta}-\nu\bar{\th}| \big)                          
             d\bar\th
             d\nu \\
    &+\int_0^1 
          \int_{|\bar\th|\le \frac12|\bar\eta|} 
          \big|\phi\big(2^m  2^{-\frac{k}{2}} 
               \sqrt{i\s+2^{2j-k}|\zeta'|^2}^{-1}\bar\th\big)\big|
             2^{m - \frac{k}{2} }  
             \sqrt{i\s+2^{2j-k}|\zeta'|^2}^{-1}|\bar\th|
             \\
        &\hskip30mm\times
             \widetilde{P_{\al}}(\t,\zeta',(\bar\eta-\nu \bar\th))
             \exp\big(-|\bar{\eta}-\nu\bar{\th}|\big)                          
             d\bar\th
             d\nu
     \\
    \equiv &III+IV,
    \label{eqn;rev-almost-ortho-10} \eqntag
 }
where we use the variables defined in \eqref{eqn;change-of-parameters} and 
set 
\eqn{
 \spl{
  \widetilde{P_{\al}}(\t,\zeta',(\bar\eta-\nu \bar\th))
  \equiv &
     P_{\al}(\t,\xi', 2^{k/2}(\eta-\nu\th))   
      - \pt_{\mu}P_{\al}(\t,\xi',\mu)
                 \big|_{\mu= 2^{\frac{k}{2}}(\eta-\nu\th)}\\
  =  &P_{\al}(\t,\xi',|p|^{-1}(\bar\eta-\nu\bar\th))   
       -\pt_{\mu}
                 P_{\al}(\t,\xi',\mu)
                 \big|_{\mu=|p|^{-1}(\bar\eta-\nu\bar\th)}. 
}
}
Since $\phi$ is rapidly decreasing function, for any $N\in \Nt$,
there exists a constant $C_N>0$ and the first term  of the 
right hand side of \eqref{eqn;rev-almost-ortho-10} is 
\algn{
   III\le &\int_0^1 \int_{|\bar\th|> \frac12|\bar\eta|} 
                \frac{C_N|2^m 2^{-\frac{k}{2}} 
                          \sqrt{i\s+2^{2j-k}|\zeta'|^2}^{-1} \bar\th|}
                     {\< 2^m 2^{-\frac{k}{2}} 
                         \sqrt{i\s+2^{2j-k}|\zeta'|^2}^{-1} \bar\th\>^{2N}}  
              \widetilde{P_{\al}}(\t,\zeta',(\bar\eta-\nu \bar\th))
              \exp\big(-|\bar{\eta}-\nu\bar{\th}| 
                 \big)                       
             d\bar\th
             d\nu
 \\
  \le & \frac{C_N 2^{-m} 2^{ \frac{k}{2}}}  
             {\< 2^m     2^{-\frac{k}{2}} 
                  \sqrt{i\s+2^{2j-k}|\zeta'|^2}^{-1} \bar\eta\>^N
             }
        \int_0^1 
          \int_{|\bar\th|> \frac12|\bar\eta|} 
              \frac{(2^m  2^{-\frac{k}{2}} 
                  \sqrt{i\s+2^{2j-k}|\zeta'|^2}^{-1})^2||\bar\th|}
                   {\<2^m 2^{-\frac{k}{2}} 
                    \sqrt{i\s+2^{2j-k}|\zeta'|^2}^{-1}\bar\th\>^N}  \\
    &\hskip50mm\times
          \widetilde{P_{\al}}(\t,\zeta',(\bar\eta-\nu \bar\th))
          \exp\big(-|\bar{\eta}-\nu\bar{\th}|\big)                        
             d\bar\th
             d\nu
\\
 \le& \frac{C_N 2^{-m+\frac{k}{2}}}
               {\< 2^m\eta\>^N}
      \int_0^1 
        \int_{\re} \frac{x}{\<x\>^N} dx
      d\nu
 \le  \frac{  C_N  2^{-m+\frac{k}{2} } }
              { \<2^m \eta\>^N}.
 \label{eqn;rev-almost-ortho-11} \eqntag
}

\vskip4mm
For the second term in 
\eqref{eqn;rev-almost-ortho-10}, we prepare the following {simple} estimate:
%
%
\vskip2mm
\begin{lem}\label{lem;b-bound} For {$N=2,3,4\cdots$} and
 $a>0$, 
\eq{
\int_{|x|\le a}\frac{dx}{(1+|x|^2)^{N/2}}
 \le \frac{C_N\; a}{(1+|a|^2)^{1/2}},
 \label{eqn;almost-orthogonal-A}
 }
 {where $C_N>0$ is a constant depending on $N$.}
 \end{lem}

\vskip3mm
\vskip3mm
{
Using  \eqref{eqn;polynomial-bound-D-T},
\eqref{eqn;change-of-parameters} and \eqref{eqn;almost-orthogonal-A},
$2j\le  k\le 2m$ and \eqref{eqn;power-bound-D-T}
with changing
\eqn{
  \bar{\bar\th}\equiv 2^m 2^{-\frac{k}{2}}|p|^{-1}\bar\th,
  \qquad   
  \bar{\bar\eta}\equiv 2^m 2^{-\frac{k}{2}}|p|^{-1}\bar\eta, 
   }
we notice that 
$|\bar{\bar\eta}-\nu\bar{\bar\th}|\ge |\bar{\bar\eta}|-|\bar{\bar\th}|
\ge \frac12|\bar{\bar\eta}|$ 
under the restriction $|\bar{\bar\th}|\le \frac12|\bar{\bar\eta}| $ 
and it follows 
}
{
\algn{
 IV
 \le &{\bk
      \int_0^1
      \int_{|\bar\th|\le \frac12|\bar\eta|} 
           \frac{C_N|2^m 2^{-\frac{k}{2}} 
                          |\sqrt{i\s+2^{2j-k}|\zeta'|^2}|^{-1}  \bar\th|}
                {\< 2^m  2^{-\frac{k}{2}} 
                          |\sqrt{i\s+2^{2j-k}|\zeta'|^2}|^{-1}  \bar\th\>^{3}}  
           \big|\widetilde{P_{\al}}(\t,\zeta',(\bar\eta-\nu \bar\th))\big|
           \exp\big(-|\bar{\eta}-\nu\bar{\th}| \big)                       
             d\bar\th
             d\nu
     }
  \\
\le &C_N2^{-m+\frac{k}{2}}\int_0^1
      \int_{|\bar\th|\le \frac12|\bar\eta|} 
           \frac{(2^m 2^{-\frac{k}{2}} 
                          |\sqrt{i\s+2^{2j-k}|\zeta'|^2}|^{-1})^2  |\bar\th|}
                {\< 2^m  2^{-\frac{k}{2}} 
                          |\sqrt{i\s+2^{2j-k}|\zeta'|^2}|^{-1}  \bar\th\>^{3}}   \\
    &\hskip5cm \times
           \big|\widetilde{P_{\al}}(\t,\zeta',(\bar\eta-\nu \bar\th))\big|
           \exp\big(-|\bar{\eta}-\nu\bar{\th}| \big)                       
             d\bar\th
             d\nu
  \\
 \le& C_N 2^{-m+\frac{k}{2}}
          \int_0^{1}
             \int_{|\bar{\bar\th}|< \frac12|\bar{\bar\eta}|} 
               \frac{|\bar{\bar\th}| }
                    {\< \bar{\bar\th}\>^{3}} 
               {\exp\big(-2^{-1}2^{-m+\frac{k}{2}}|p|
               |\bar{\bar\eta}-\nu\bar{\bar\th}| \big) }
              d\bar{\bar\th}
             d\nu
   \\
 \le& C_N 2^{-m+\frac{k}{2}}
          \int_0^{1}
             \int_{|\bar{\bar\th}|< \frac12|\bar{\bar\eta}|} 
               \frac{1}
                    {\< \bar{\bar\th}\>^{2}} 
               {
               \exp\big(
                -2^{-1}2^{-m+\frac{k}{2}}|p|
                |\bar{\bar\eta}-\nu\bar{\bar\th}| \big) }
              d\bar{\bar\th}
             d\nu
    \\
 \le& C_N 2^{-m+\frac{k}{2}}
         { \exp\big(
                -2^{-2}2^{-m+\frac{k}{2}}|p|
                |\bar{\bar\eta}| \big) 
         }
             \int_{|\bar{\bar\th}|< \frac12|\bar{\bar\eta}|} 
               \frac{1}{\< \bar{\bar\th}\>^{2}} 
              d\tilde\th
    \\
 \le& C_N 2^{-m+\frac{k}{2}}
         { \exp\big(
                -2^{-2}2^{-m+\frac{k}{2}}|p|
                |\bar{\bar\eta}| \big) 
         }
          \frac{|\bar{\bar\eta}|}{(1+ |\bar{\bar\eta}|^2)^{1/2}} 
\\
\le& C_N2^{-m+\frac{k}{2}}{ \exp\big(
                -2^{-2}
                |\bar{\eta}| \big) 
         }
\\
 \le & \frac{ C_N 2^{-m+\frac{k}{2}} }
            { \<2^{\frac{k}{2}}\eta\>^{N} }. 
 \label{eqn;rev-almost-ortho-12} \eqntag
}
}
Hence we obtain from \eqref{eqn;rev-almost-ortho-11}  and
 \eqref{eqn;rev-almost-ortho-12} that 
\eq{ \label{eqn;rev-almost-ortho-15} \eqntag
 \bigg| \int_{\re_+}\phi_m(\eta-\th) 
             \Big(
               P_{\al}(\t,\xi', 2^{k/2} \th)
              \exp\big(- 2^{\frac{k}{2}}\th
                       \sqrt{i\s+2^{2j-k}|\zeta'|^2} \big)
             \Big)
              d\th
 \bigg|
 \le 
 \frac{ C_N 2^{-m+\frac{k}{2}} }
      {\<2^{\frac{k}{2}}\eta\>^{N} }.
 }
The estimates 
\eqref{eqn;rev-almost-ortho-4} and \eqref{eqn;rev-almost-ortho-15}
yield \eqref{eqn;almost-orthogonal-4}. 

\noindent
$\<$ Step 3$\>$ The case  {$k<2j$}:\par
To show \eqref{eqn;almost-orthogonal-5}, we use
$q(\sg,\zeta',2^{2j-k})\equiv\sqrt{i2^{k-2j}\s+|\zeta'|^2}$
instead of \eqref{eqn;almost-orthogonal-6} and let 
\algn{
 Q_{\al}(\sg,\zeta',\eta,2^j)
 =&\exp\big( 2^j\eta\sqrt{i 2^{k-2j}\sg+|\zeta'|^2}\big)
   \sum_{|\al|\le n+2}
    C_{n,\al}
     \pt^{\al_1}_{\sg}
     \pt^{\al_2}_{\zeta'}
    \exp\big(- 2^j\eta\sqrt{i 2^{k-2j}\sg+|\zeta'|^2}\big). 
}
Since $|Q_{\al}|$ is a polynomial of $\eta$
and $|q(\sg,\zeta',2^{2j-k})|^{-1}$, it follows from the 
assumption that 
\eqn{ \label{eqn;power-bound-D-S}
 2^{-1}\le
  \Big|q(\sg,\zeta',2^{k-2j})\Big|
  = \Big( 2^{2(k-2j)} \sg^2+|\zeta'|^4\Big)^{1/4} 
  \le 20^{1/4}{.}
}
Hence as in the previous step, 
\eqn{  \label{eqn;polynomial-bound-D-S}
 | Q_{\al}(\t,\xi', 2^{j}\eta)|
 \le C\big(1+(2^{j}\eta)^{n+2}\big).
 }
All the other estimate is very similar to the case of 
$\<$Step 2$\>$ and all the terms involving 
$ 2^{\frac{k}{2}}$ are arranged into  $ 2^j$.
The estimate corresponding to the case  $\<$Step 1$\>$ is 
redundant by the assumption $j<m$.
\end{prf}

\noindent
{\bf Proof of Lemma \ref{lem;pt-orthogonal-D2}, continued.}
The proof of Lemma \ref{lem;pt-orthogonal-D2} goes in a similar 
way to the case of Lemma \ref{lem;pt-orthogonal-1}.
After integrating by parts, the integrable factors  
$\<t\>^{-2}$ and $\<y'\>^{-n}$ appear and then estimate the 
integrant to obtain the desired estimate in 
Lemma \ref{lem;almost-orthogonal-4} for the case $2j< k$.
The other case $2j\ge k$ is also obtained from \eqref{eqn;almost-orthogonal-5}
in Lemma \ref{lem;almost-orthogonal-4}.
This complete the proof of Lemma \ref{lem;pt-orthogonal-D2}.
\end{prf}

\subsection{The Neumann potential case}
The almost orthogonal estimate for the Neumann boundary potential 
is very similar to the case of Dirichlet potential case except the
order of the derivative. The following lemma shows the estimates 
\eqref{eqn;Ortho-N} hold valid.

\begin{lem}[A crucial potential orthogonality]
\label{lem;pt-orthogonal-2}
For  $k,j,\ell\in \Z$ let $\{\psi_k(t)\}_{k\in \Z}$ and 
$\{\phi_j(x)\}_{j\in \Z_+}$ 
be the time and the space Littlewood-Paley dyadic decomposition and
let $\Psi_N(t,x',\eta)$ be the boundary potential defined in \eqref{eqn;potential-N}. 
Set 
\eqn{
 \Psi_{N,k,j}(t,x',\eta)
   \equiv \int_{\re}\int_{\re^{n-1}}
        \Psi_N(t-s,x'-y',\eta)\psi_k(s)\phi_j(y')dy'ds
 }
for $\eta=x_n\in I_{\ell}=[2^{-\ell},2^{-\ell+1})$.
Then there exists a constant $C_n>0$ depending only on the dimension $n$ 
satisfying 
\eq{ \label{eqn;crucial-potential-orthogonarity-N}
  \| \Psi_{N,k,j}(t,\cdot,\eta) \|_{L^1_{x'}}
  \le
   \left\{
   \begin{aligned}
      &C_n 2^{\frac{k}{2}}\big(1+(2^{\frac{k}{2}}\eta)^{n+2}\big)
         e^{-2^{\frac{k}{2}-1}\eta} \frac{2^k}{\<2^kt\>^{2}} ,
      &{k\ge 2j},\\
      &C_n2^{\frac{k}{2}}\big(1+(2^{j}\eta)^{n+2}\big)
         e^{-2^{j-1}\eta}\frac{2^k}{\<2^kt\>^{2}},
      &{k<2j },
   \end{aligned}
   \right.
  }
and
\eq{ \label{eqn;crucial-potential-orthogonarity-N2}
  \|\phi_m\underset{(\eta)}{*}\Psi_{N,k,j}(t,\cdot,\eta) \|_{L^1_{x'}}
  \le
   \left\{
   \begin{aligned}
      &C_n 2^{\frac{k}{2}}
          \frac{2^{-|\frac{k}{2}-m|}}{\<{2^{\min(\frac{k}2,m)} }\eta\>^N} 
          \frac{2^k}{\<2^kt\>^{2}} ,
      &{k\ge 2j},\\
      &C_n 2^{\frac{k}{2}}
          \frac{2^{-|j-m|}}{\<{2^j} \eta\>^N}
          \frac{2^k}{\<2^kt\>^{2}},
      &{k<2j}.
   \end{aligned}
   \right.
  }
\end{lem}

The proof for Lemma \ref{lem;pt-orthogonal-2} is shown in a parallel way 
to the proof of Lemma \ref{lem;pt-orthogonal-1}.  The only difference stems
from the difference of two potentials $\Psi_D$ and $\Psi_N$ of 
boundary data and the difference from \eqref{eqn;Dirichlet-ortho-1} 
reflects the order of spacial derivatives 
appearing the Fourier image of the following expression:
\eqn{
 \spl{
   &\Psi_{N,k,j}(t,x',\eta)\\
  =& -c_{n+1}\int_{\re}\int_{\re^{n-1}}
         e^{it\t+ix'\cdot \xi'}
        \frac{i\t}{\sqrt{i\t+|\xi'|^2}}
        \exp\big(-\sqrt{i\t+|\xi'|^2}\eta\big)
        \widehat{\psi}(2^{-k}\t)\widehat{\phi}(2^{-j}\xi')
     d\xi'd\t.
}
}

\vskip1mm
\begin{prf}{Lemma \ref{lem;pt-orthogonal-2}}
We only show  the out-lined proof of \eqref{eqn;crucial-potential-orthogonarity-N}.
The other estimate \eqref{eqn;crucial-potential-orthogonarity-N2} follows very much in a similar way to the proof of the Dirichlet case
in Lemma \ref{lem;pt-orthogonal-D2}.

We consider a time-like estimate {$k\ge  2j$} as is in the Dirichlet case.
Taking $\zeta'$-space cut-off, we have by 
using the change of variables 
         $\t=2^k\s$, $\xi'=2^j\zeta'$ that
\eqn{
 \spl{
   \| &\Psi_{N,k,j}(t,\cdot,\eta)\|_{L^1_{x'}} \\
  =&\left\|
     c_{n+1}\int_{\re}\int_{\re^{n-1}}
         e^{it\t+ix'\cdot \xi'}
        \frac{i\t }{\sqrt{i\t+|\xi'|^2}}
               \exp\big(-\sqrt{i\t+|\xi'|^2}\eta\big)
                \widehat{\psi}(2^{-k}\t)\widehat{\phi}(2^{-j}\xi')
     d\xi'd\t   
    \right\|_{L^1_{x'}}
\\
     =& 2^{\frac{k}{2}}
         \left\|
         c_{n+1}\int_{\re}\int_{\re^{n-1}}
            e^{i2^kt\s+iy'\cdot \zeta'}
            \frac{\s}
                 {\sqrt{i\s+2^{2j-k}|\zeta'|^2}}
            \exp\left(-2^{ \frac{k}{2} }\eta\,
                          \sqrt{i\s+2^{2j-k}|\zeta'|^2}
                    \right)
             \widehat{\psi}(\s)\widehat{\phi}(\zeta')
             d\zeta'\cdot 2^kd\s   
         \right\|_{L^1_{y'}},
   }
}
where we set $x'=2^{-j}y'$ in the last line.
Using \eqref{eqn;exp-intgralbyparts} and integrating by parts, 
we see by setting 
$p(\sg,\zeta';2^{2j-k})\equiv\sqrt{i\s+2^{2j-k}|\zeta'|^2}$
that
\algn{
  &\|\Psi_{N,k,j}(t,\cdot,\eta)
     \|_{ L^1_{x'}(B_{2^{j}}^c) } \\
   =&2^{\frac{k}{2}}
     \bigg\|
      c_{n+1}\int_{\re}\int_{\re^{n-1}}
         \frac{1}{i2^kt}
         \left(-\frac{1}{|y'|^2}\right)e^{i(2^kt\s+ y'\cdot \zeta')}\\
    &\phantom{2^{k-2\ell}\|c_n\int_{\re}\int} 
         \times
         \Del_{\zeta'}\frac{\pt}{\pt \s}
          \left[\frac{1}{p(\sg,\zeta';2^{2j-k})}
                \exp\left(-2^{ \frac12 k} \eta\,
                          p(\sg,\zeta';2^{2j-k})
                    \right)
          \widehat{\psi}(\s)
          \widehat{\phi}(\zeta')
        2^{k}\s\right] d\zeta' d\s   
    \bigg\|_{L^1_{y'}(B_1^c)}.
 }
Since 
\algn{
 \frac{\pt}{\pt \s}&
         \left[ \frac{1}{p(\sg,\zeta';2^{2j-k})}
                \exp\left(
                   -2^{\frac{k}2}\eta\,
                    p(\sg,\zeta';2^{2j-k})
                    \right)
               \widehat{\psi}(\s)2^{k}\s 
               \widehat{\phi}(\zeta')
          \right]\\
   =& \frac{1}{p(\sg,\zeta';2^{2j-k})}
      \exp\Big(-2^{\frac{k}2}\eta\,p(\sg,\zeta';2^{2j-k})\Big) \\
      &\hskip1cm
       \times \Bigg\{-\frac{2^{\frac{k}2}\eta\,i}{2p(\sg,\zeta';2^{2j-k})}
                       \widehat{\psi}(\s)2^{k}\s
                     -\frac{i}{2p(\sg,\zeta';2^{2j-k})^2}
                      \widehat{\psi}(\s)2^{k}\s 
                     +\widehat{\psi}'(\s)2^{k}\s 
                     + \widehat{\psi}(\s)2^{k}
             \Bigg\}\widehat{\phi}(\zeta').
 }
we take second derivative by setting 
$r=|\zeta'|$ and $\Del_{\zeta'}=\pt_r^2+\frac{n-2}{r}\pt_r$ 
to have 
\begin{align*}
 &\Del_{\zeta'} \frac{\pt}{\pt \s}
         \left[\frac{1}{p(\sg,\zeta';2^{2j-k})}
               \exp\left(-2^{\frac{k}2}\eta\,
                         p(\sg,\zeta';2^{2j-k})
               \right)
               \widehat{\psi}(\s)2^{k}\s 
               \widehat{\phi}(\zeta')
          \right]
 \\
 =&\left(\pt_r+\frac{n-2}{r}\right)
    \frac{\exp\big(-2^{\frac{k}2}\eta\,p(\sg,\zeta';2^{2j-k})\big)}
           {p(\sg,\zeta';2^{2j-k})}\\
     &\qquad \times
           \Bigg[
             \bigg(-\frac{ { 2^{\frac{k}{2}}\eta} 2^{2j-k}\,|\zeta'|}
                         {p(\sg,\zeta';2^{2j-k})}
                    \widehat{\phi}(\zeta')
                   -\frac{ 2^{2j-k}\,|\zeta'|}
                         { p(\sg,\zeta';2^{2j-k})^2}
                    \widehat{\phi}(\zeta')
                    +\widehat{\phi}'(\zeta')
             \bigg)\\
     &\qquad \qquad \times
             \bigg(-\frac{2^{\frac{k}2}\eta\,i}{2p(\sg,\zeta';2^{2j-k})}
                       \widehat{\psi}(\s)2^{k}\s
                     -\frac{i}{2p(\sg,\zeta';2^{2j-k})^2}
                      \widehat{\psi}(\s)2^{k}\s 
                     +\widehat{\psi}'(\s)2^{k}\s 
                     + \widehat{\psi}(\s)2^{k}
             \bigg)\\
      &\qquad\qquad
           +\bigg(\frac{ 2^{\frac{k}2}\eta\,2^{2j-k} |\zeta'|i}
                        {2p(\sg,\zeta';2^{2j-k})^3}
                  +\frac{ 2^{2j-k} |\zeta'|i}{p(\sg,\zeta';2^{2j-k})^4}
            \bigg)
                 \widehat{\psi}(\s)2^{k}\s
                 \widehat{\phi}(\zeta')
          \Bigg].
 \eqntag\label{eqn;almost-orth-N-1}
\end{align*} 

As the Dirichlet boundary case,  all terms again have
$
\exp\big(- 2^\frac{k}{2}\eta  p(\sg,\zeta';2^{2j-k})\big),
$
with
$\widehat{\psi}(\sg)$, $\widehat{\phi}(\zeta')$
or their derivatives.  
The function in denominator $p(\sg,\zeta';2^{2j-k})=\sqrt{i\s+2^{2j-k}|\zeta'|^2}$ 
is estimated from below by $2^{-1}$ thanks to the cut-off function 
$\widehat{\phi}(\zeta')$ or its derivative. Therefore if we derivate $n+2$ times, 
then it is estimated from below by $2^{n+2}$. 
 Moreover $p(\sg,\zeta';2^{2j-k})$ which arises by derivative contains surplus 
scale parameter $2^{2j-k}$ and it is estimated from above by $1$ by $k\ge 2j$.  
Thus for the functions
\eqn{
 \spl{
  \widetilde{k_{k,j}}(\sg,\zeta',\eta)
   \equiv& \frac{\exp\big(-2^{\frac{k}2}\eta\,
                        p(\sg,\zeta';2^{2j-k})\big)}
                {p(\sg,\zeta';2^{2j-k})}
          \widehat{\psi}(\s)
          \widehat{\phi}(\zeta')
        2^{k}\s {,}\\
   \widetilde{K_{k,j}^{n+2} }(\sg,\zeta',\eta,2^{2j-k})
   \equiv&\exp\big(2^{\frac{k}2}\eta\, p(\sg,\zeta';2^{2j-k})\big)
          D^{n+2}_{\sg,\zeta'}\widetilde{k_{k,j}}(\sg,\zeta',\eta),
 }
}
where we put
$$
  D^{n+2}_{\sg,\zeta'}
   \equiv(1-\Del_{\zeta'})^{\frac{n}{2}}
            \left(1-\frac{\pt^2}{\pt \sg^2}\right). 
$$
Under the condition  $k\ge 2j$, i.e., the time-like 
condition   $2^{2j-k}\le 1$, we obtain the following estimate:
\begin{align*}
 &\|\Psi_{N,k,j}(t,\cdot,\eta)
  \|_{ L^1_{x'}} \\
   =& 2^{\frac{k}{2}}
     \left\|\left(\frac{1}{\<y'\>^2}\right)^{\frac{n}{2}}
            \frac{2^k}{\<2^kt\>^2}
        c_{n+1}\int_{\re}\int_{\re^{n-1}}  
            e^{i(2^kt\sg+y'\cdot \zeta')}
            (1-\Del_{\zeta'})^{\frac{n}{2}}
            \left(1-\frac{\pt^2}{\pt \sg^2}\right)
            \widetilde{k_{k,j}}(\sg,\zeta',\eta)
            d\zeta'd\sg
     \right\|_{L^1_{y'}}\\
    =&\, C\,2^{\frac{k}{2}} \exp(-2^{\frac{k}{2}-1}\eta) \frac{2^k}{\<2^kt\>^2}\\
     &\quad\times 
     \Bigg\|\left(\frac{1}{|y'|^2}\right)^{\frac{n}{2}}           
         \int_{\re}\int_{\re^{n-1}}  
            e^{i(2^kt\sg+y'\cdot \zeta')}
            e^{\big(-2^{\frac{k}{2}}\eta\,
                      \big(p(\sg,\zeta';2^{2j-k})-\frac12\big)
                \big)
              }
            \widetilde{K_{k,j}^{n+2} }(\sg,\zeta',\eta,2^{2j-k})
            d\zeta'd\sg
     \Bigg\|_{L^1_{y'}}\\
 \le &\,C\,2^{\frac{k}{2}} \exp(-2^{\frac{k}{2}-1}\eta) 
        \frac{2^k}{\<2^kt\>^2}        
         \int_{\re}\int_{\re^{n-1}}
            e^{\big(-2^{\frac{k}{2}}\eta\,
                      \big(\sqrt{i\s+2^{2j-k}|\zeta'|^2}-\frac12\big)
               \big)
              }
           \Big|
            \widetilde{K_{k,j}^{n+2} }(\sg,\zeta',\eta,2^{2j-k})
           \Big|
         d\zeta'd\sg
   \\
 \le&C2^{\frac{k}{2}}\, \exp(-2^{\frac{k}{2}-1}\eta)
    \big(1+(2^{\frac{k}{2}}\eta)^{n+2}\big)\frac{2^k}{\<2^kt\>^2}.
\end{align*}

 For the case of the space-like region {$k< 2j$}, 
we proceed similar way. By changing  
$\t=2^k\s$, $\xi'=2^j\zeta'$, 
\begin{align*}
   \|\Psi_{N,k,j}&(t,\cdot,\eta)\|_{L^1_{x'}} \\
  =&\left\|
     c_{n+1}\int_{\re}\int_{\re^{n-1}}
         e^{it\t+ix'\cdot \xi'}
          \frac{\t }{\sqrt{i\t+|\xi'|^2}}
                \exp\big(-\sqrt{i\t+|\xi'|^2}\eta\big)
                \widehat{\psi}(2^{-k}\t)\widehat{\phi}(2^{-j}\xi')
     d\xi'd\t   
    \right\|_{L^1_{x'}}
\\
     =& 2^{\frac{k}{2}}
      \bigg\|
      c_{n+1}\int_{\re}\int_{\re^{n-1}}
            e^{i2^kt\s+iy'\cdot \zeta'}
             \frac{\sg}
             {\sqrt{i\s+2^{2j-k}|\zeta'|^2}}
             \exp\Big(-2^{j}\eta\sqrt{2^{k-2j}i\s +|\zeta'|^2}
                 \Big)
      \\
      &\hskip8cm  \times 
              \widehat{\psi}(\s)\widehat{\phi}(\zeta')
             d\zeta'\cdot 2^kd\s   
         \bigg\|_{L^1_{y'}}.
\end{align*}
Here we set $x'=2^{-j}y'$  in the last line.
Using \eqref{eqn;exp-intgralbyparts} and integrating by parts, 
\eqn{
 \spl{
  &\|\Psi_{N,k,j}(t,\cdot,\eta)
     \|_{ L^1_{x'}(B_{2^{\frac{k}{2}}}^c) } \\
   =& 2^{\frac{k}{2}}
     \bigg\|
      c_{n+1}\int_{\re}\int_{\re^{n-1}}
         \frac{1}{i2^kt}
         \left(-\frac{1}{|y'|^2}\right)e^{i(2^kt\s+ y'\cdot \zeta')}\\
    &\phantom{2^{k} } 
         \times
         \Del_{\zeta'}\frac{\pt}{\pt \s}
          \left[
          \frac{\sg }
               {\sqrt{i\s+2^{2j-k}|\zeta'|^2}}        
             \exp\Big(-2^{j}\eta
                    \sqrt{\big(2^{k-2j}\s i+|\zeta'|^2\big)}
                 \Big)
          \widehat{\psi}(\s)
          \widehat{\phi}(\zeta')
        2^{k}\right] d\zeta' d\s   
    \bigg\|_{L^1_{y'}(B_1^c)}.
   }
}
Setting
\eqn{\label{eqn;p-q}
 \begin{cases}
  &p(\sg,\zeta';2^{2j-k})=\sqrt{i\s +2^{2j-k}|\zeta'|^2}, \\
  &q(\sg,\zeta';2^{k-2j})=\sqrt{2^{k-2j}i\s +|\zeta'|^2},
 \end{cases}
}
we see that 
\begin{align*}
 \frac{\pt}{\pt \s}&
         \left[ \frac{\exp\big(-2^{j}\eta\,
                    q(\sg,\zeta';2^{k-2j})
                    \big)}
                  {p(\sg, \zeta';2^{2j-k})}
               \widehat{\psi}(\s)2^{k}\s 
               \widehat{\phi}(\zeta')
          \right]\\
   =& \frac{\exp\Big(-2^{j}\eta\,
            q(\sg,\zeta';2^{k-2j})\Big)}
          {p(\sg,\zeta';2^{2j-k})}\\
    & \qquad \times
           \Bigg[
             \bigg(-\frac{2^{j}\eta\, 2^{k-2j} i}
                       {2q(\sg,\zeta';2^{k-2j})}
                  -\frac{i}
                       {2p(\sg,\zeta';2^{2j-k})^2}
             \bigg)
                 \widehat{\psi}(\s)2^{k}\s 
                + \widehat{\psi}'(\s)2^{k}\s 
                      + \widehat{\psi}(\s)2^{k}
          \Bigg]\widehat{\phi}(\zeta').
\end{align*}
Taking the second derivative in $r=|\zeta'|$
\begin{align*}
 &\Del_{\zeta'} \frac{\pt}{\pt \s}
         \left[\frac{\sg}
               {\sqrt{i\s+2^{2j-k}|\zeta'|^2}}  
               \exp\left(-2^{j}\eta
                    \sqrt{2^{k-2j}i\s+|\zeta'|^2}
               \right)
               \widehat{\psi}(\s)2^{k}\s 
               \widehat{\phi}(\zeta')
          \right]
\\
    =&\Del_{\zeta'}\frac{\exp\Big(-2^{j}\eta\,
                         q(\sg,\zeta';2^{k-2j})\Big)}
                        {p(\sg,\zeta';2^{2j-k})} \\
    & \qquad \times
           \Bigg[
             \bigg(-\frac{2^{j}\eta\,2^{k-2j} i}
                         {2q(\sg,\zeta';2^{k-2j})}
                  -\frac{i}
                        {2p(\sg,\zeta';2^{2j-k})^2}
             \bigg)
                 \widehat{\psi}(\s)2^{k}\s 
                + \widehat{\psi}'(\s)2^{k}\s 
                      + \widehat{\psi}(\s)2^{k}
          \Bigg]\widehat{\phi}(\zeta').
\end{align*}

Similarly before all terms again have
$
\exp\big(-  2^\frac{k}{2}\eta  q(\sg,\zeta';2^{k-2j})\big),
$
 $\widehat{\psi}(\sg)$, $\widehat{\phi}(\zeta')$
or their derivatives. 
The denominator $q(\sg,\zeta';2^{k-2j})$ and $p(\sg,\zeta',2^{2j-k})$
are estimated from below by $2^{-1}$ thanks to the cut-off functions 
$\widehat{\psi}(\sg)$ and $\widehat{\phi}(\zeta')$ or its derivative. 
Moreover $p(\sg,\zeta';2^{2j-k})$ is estimated from above by $1$ by
$k\le 2j$.  
Thus introducing  
\eqn{
 \spl{
  \widetilde{ h_{k,j}}(\sg,\zeta',\eta)
   \equiv& \frac{\exp\big(-2^{j}\eta
                    \sqrt{2^{k-2j}\s i+|\zeta'|^2}\big)
                }
                { \sqrt{\s i+2^{2j-k}|\zeta'|^2}}
          2^{k}\s\widehat{\psi}(\s)
          \widehat{\phi}(\zeta'),       
  \\
  \widetilde{H_{k,j}^{n+2} }(\sg,\zeta',\eta,2^{k-2j})
   \equiv&\frac{\exp\big(-2^{j}\eta
                    \sqrt{2^{k-2j}\s i+|\zeta'|^2}\big)
                }
                { \sqrt{\s i+2^{2j-k}|\zeta'|^2}}
          D^{n+2}_{\sg,\zeta'} \widetilde{ h_{k,j}}(\sg,\zeta',\eta),
 }
}
we obtain
\begin{align*}
 &\|\Psi_{N,k,j}(t,\cdot,\eta)
  \|_{ L^1_{x'}} \\
   =&2^{\frac{k}{2}}
     \left\|\left(\frac{1}{\<y'\>^2}\right)^{\frac{n}{2}}
            \frac{2^k}{\<2^kt\>^2}
     c_{n+1}\int_{\re}\int_{\re^{n-1}}  
            e^{i(2^kt\sg+y'\cdot \zeta')}
            (1-\Del_{\zeta'})^{\frac{n}{2}}
            \left(1-\frac{\pt^2}{\pt \sg^2}\right)
            \widetilde{ h_{\ell,k,j}}(\sg,\zeta',\eta)
            d\zeta'd\sg
     \right\|_{L^1_{y'}}
\\
 \le &C\,2^{\frac{k}{2}} \exp(-2^{j-1}\eta) \frac{2^k}{\<2^kt\>^2}
           \int_{\re}\int_{\re^{n-1}}
            e^{\left(-2^{j}\eta
                      \big(q(\sg,\zeta';2^{k-2j})-1\big)
                \right)
              }
           \Big|
            \widetilde{H_{k,j}^{n+2} }(\sg,\zeta',\eta,2^{k-2j})
           \Big|
         d\zeta'd\sg
   \\
 \le&C\,2^{\frac{k}{2}}\exp(-2^{j-1}\eta)
    \big(1+(2^{j}\eta)^{n+2}\big)\frac{2^k}{\<2^kt\>^2}.
\end{align*}

This completes the proof of 
\eqref{eqn;crucial-potential-orthogonarity-N}  in 
Lemma \ref{lem;pt-orthogonal-2}.
\end{prf}

%
%

\sect{Optimal boundary trace estimates}
In this section we show the optimality for the boundary 
trace estimate shown in Theorem \ref{thm;boundary-trace-D} and 
\ref{thm;boundary-trace-N}. 
This shows that the condition on the boundary data in 
those theorems are not only the sufficient condition but also 
the necessary condition (see for more detailed estimates 
for the boundary trace \cite{JS}, \cite{M-V14})

\noindent
\begin{thm}[The Dirichlet boundary trace]
\label{thm;sharp-boundary-trace-D} 
Let $1\le p\le \infty$ and $s\in \re$.  
There exists a constant $C>0$ such that 
for all function 
$ u=u(t,x',\eta)\in \dot{W}^{1,1}(\re_+; \dB^s_{p,1}(\re^{n}_+))$, ${\Del u(t,x',\eta)\in } L^1(\re_+; \dB^{s}_{p,1}(\re^{n}_+))$ with satisfying 
$u(0,x',\eta)=0$ almost everywhere in $(x',\eta)\in \re^{n}_+$, 
the following estimate holds:
\begin{align}
   \sup_{\eta\in\re_+}
   \Big(\|u(\cdot,\cdot, \eta)
        \|_{\dot{F}^{1-1/2p}_{1,1}(\re_+; \dB^s_{p,1}(\re^{n-1}))}
  &+\|u(\cdot,\cdot, \eta)\|_{L^1(\re_+; \dB^{s+2-1/p}_{p,1}(\re^{n-1}))}
  \Big)  \nonumber\\
  \le & C\Big(
         \|\pt_t u\|_{L^1(\re_+;\dB^s_{p,1}(\re^{n}_+))}
       + \|\N^2  u\|_{L^1(\re_+;\dB^s_{p,1}(\re^{n}_+))}
       \Big).\label{eqn;sharp-trace-D}
\end{align}
\end{thm}

From Theorem \ref{thm;sharp-boundary-trace-D} the following 
corollary immediately follows:
\vskip2mm
\noindent
\begin{cor}[Sharp boundary trace for Dirichlet data]
\label{cor;heat-boundary-trace-D} 
For $1\le p\le \infty$, there exists a constant $C>0$ such that
the solution $u$ to the initial-boundary value problem 
\eqref{eqn;heat-D-0}
in
$$
  {u\in} \dot{W}^{1,1}\big(\re_+;\dB^s_{p,1}(\re^n_+)\big),\quad
  {\Del u\in} L^1\big(\re_+;\dB^{s+2}_{p,1}(\re^n_+)\big)
$$
gives the following estimate on the boundary condition $h$:
\eqn{
   \|h\|_{\dot{F}^{1-1/2p}_{1,1}(\re_+; \dB^s_{p,1}(\re^{n-1}))}
  +\|h\|_{L^1(\re_+; \dB^{s+2-1/p}_{p,1}(\re^{n-1}))} 
  \le 
  \left\{
  \begin{aligned}
   & C\|\N^2  u\|_{L^1(\re_+;\dB^s_{p,1}(\re^{n}_+))}, \\
   & C\|\pt_t  u\|_{L^1(\re_+;\dB^s_{p,1}(\re^{n}_+))}.
   \end{aligned} 
  \right.
}
\end{cor}

\begin{prf}{Theorem \ref{thm;sharp-boundary-trace-D}}
For $1\le p\le \infty$, we assume that 
$u\in \dot{W}^{1,1}(\re_+;\dB^s_{p,1}{(\re^n_+)})$,
${\Del u\in  L^1(\re_+;\dB^{s}_{p,1} (\re^n_+)})$. 
{\allowdisplaybreaks
\begin{align*}
 &\big\|u(\cdot,\cdot, \eta) 
 \big\|_{\dF^{1-1/2p}_{1,1}\big(\re_+;\dB^{s}_{p,1}(\re^{n-1}_{x'})\big)} 
  \\
  \le & 
  \Big\|{ \sum_{j\in\Z}
              \sum_{k\ge 2j}  }
              2^{(1-1/2p)k}2^{sj}
              \big\|
                 \psi_k\underset{(t)}{*}
                 \phi_j\underset{(x')}{*}
                 u(t,\cdot,\eta)
              \big\|_{L^p(\re^{n-1}_{x'})}
       \Big\|_{L^1_t(\re_+)} \\
  &+\Big\|{ \sum_{j\in\Z}
                \sum_{k\le 2j}  }
            2^{(1-1/2p)k}2^{sj}
              \big\|
                 \psi_k\underset{(t)}{*}
                 \phi_j\underset{(x')}{*}
                 u(t,\cdot,\eta)
              \big\|_{L^p(\re^{n-1}_{x'})}
       \Big\|_{L^1_t(\re_+)} 
 \equiv  I+II. \addtocounter{equation}{1}
 \tag{\theequation} \label{eqn;time-like-0}
 \end{align*}
 }
Using the assumption $u(0,x',\eta)=0$ almost everywhere,
\begin{align}
 \psi_k(t)\underset{(t)}{*} u(t,x',\eta)
 =&-\int_{\re_+}\pt_s\big(\int_s^{\infty}\psi_k(t-r)dr\big)  u(s,x',\eta)ds  \nonumber\\
 =&-\Bigl[\big(\int_s^{\infty} \psi_k(t-r)dr\big)  u(s,x',\eta)
    \Bigr]_{s=0}^{\infty}
   +\int_{\re_+}\int_s^{\infty}\psi_k(t-r)dr  \pt_s u(s,x',\eta)ds  \nonumber\\
 =&\pt_t^{-1}\psi_k(t)\underset{(t)}{*}  \pt_t u(t,x',\eta),
 \label{eqn;time-like-int}
\end{align}
where we set
\eq{\label{eqn;int-psi-k}
  \pt_t^{-1}\psi_k(t-s)
  \equiv\int_{-\infty}^{t-s}\psi_k(r)dr
  =\int_s^{\infty}\psi_k(t-r)dr.
}
Here we recall the Littlewood-Paley decomposition of direct sum type 
defined in \eqref{eqn;direct-sum-L-P}.
Since $\pt_t^{-1} \psi_k=2^{-k}(\pt_t^{-1} \psi)_k$ 
is also rapidly decreasing 
smooth function,  $\zeta_{j}(\eta)\underset{(\eta)}{*}
\zeta_{j-1}(\eta)=\zeta_{j-1}(\eta)$ by definition \eqref{eqn;zeta} 
and using \eqref{eqn;Sum_Phi_is_1},   
the Hausdorff-Young inequality with \eqref{eqn;time-like-int},  
there exists a constant $C>0$ such that 
{\allowdisplaybreaks
\begin{align*}
  I = & \bigg\|{\sum_{j\in\Z}
                    \sum_{k\ge 2j}  }
                    2^{(1-1/2p)k}2^{sj}
              \big\|
                    \pt_t^{-1}\psi_k\underset{(t)}{*}
                    \phi_j\underset{(x')}{*}
                   \pt_t u(t,\cdot,\eta)
             \big\|_{L^p(\re^{n-1}_{x'})}
        \bigg\|_{L^1_t(\re_+)}
  \\
   \le &C\bigg\|{\sum_{j\in\Z}2^{sj}
                     \sum_{k\ge 2j}  } 
                  2^{(1-/2p)k}
                 |2^{-k}(\pt_t^{-1}\psi)_k(t)| 
                \underset{(t)}{*}
                \big\|
                    \phi_j\underset{(x')}{*}
                    \pt_t u(t,\cdot,\eta)
                \big\|_{L^p(\re^{n-1}_{x'})}
         \bigg\|_{L^1_t(\re_+)}
 \\
   \le &C{\sum_{j\in \Z} 2^{sj} \sum_{k\ge 2j} 2^{-\frac1{2p}k} }
         \bigg\|\int_{\re_+}
                 \frac{2^k}{\<2^k(t-s)\>^2}
                 \big\|
                 \phi_j\underset{(x')}{*}
                 \pt_t u(t,\cdot,\eta) 
                \big\|_{L^p(\re^{n-1}_{x'})}               
                ds
        \bigg\|_{L^1_t(\re_+)}
 \\
   \le &C{\sum_{j\in \Z} 2^{sj} \sum_{k\ge 2j}  2^{-\frac1{2p}k} }
          \Big\|\frac{2^k}{\<2^k t\>^2}\Big\|_{L^1_t(\re_+)}
          \bigg\|
                 \Big\|
                    \sum_{m\in \Z} \overline{\Phi_m}\underset{(x',\eta)}{*} 
                 \phi_j\underset{(x')}{*}
                 \pt_t u(t,\cdot,\eta)  
                \Big\|_{L^p(\re^{n-1}_{x'})}
          \bigg\|_{L^1_t(\re_+)}
  \\
   \le &C\sum_{j\in \Z}2^{sj}
          \bigg\| {\sum_{k\ge 2j}  2^{-\frac1{2p}k} }
             \Big\|
              \sum_{|m-j|\le 1} \overline{\Phi_m}\underset{(x',\eta)}{*} 
              \phi_j\underset{(x')}{*}
              \zeta_{j}(\eta)\underset{(\eta)}{*}
              \zeta_{j-1}(\eta)\underset{(\eta)}{*} 
              \pt_t u(t,\cdot,\eta)  
            \Big\|_{L^p(\re^{n-1}_{x'})}                
          \bigg\|_{L^1_t(\re_+)}
  \\
   \le &C\sum_{j\in \Z}2^{sj}
         \bigg\| {\sum_{k\ge 2j}2^{-\frac1{2p}k}}
                   \|\zeta_{j}(\eta)
                   \|_{L^{p'}(\re_{+,\eta})}
              \Big\|\big\|
                 \overline{\Phi_j}\underset{(x',\eta)}{*} 
                 \pt_t u(t,\cdot,\eta)  
                 \big\|_{L^p(\re^{n-1}_{x'})}
              \Big\|_{L^p(\re^{+}_{\eta})}                
         \bigg\|_{L^1_t(\re_+)}
  \\
      \le &C\bigg\| \sum_{j\in \Z}2^{sj}
                 \Big\|
                  \overline{\Phi_j}\underset{(x',\eta)}{*} 
                  \pt_t u(t,\cdot,\eta)  
                \Big\|_{L^p(\re^{n}_{+,(x',\eta)})}              
          \bigg\|_{L^1_t(\re_+)} 
    \le C \|\pt_t u\|_{ L^1(\re_+;\dB^s_{p,1}(\re^n_{+}) }.
     \addtocounter{equation}{1}
 \tag{\theequation} \label{eqn;time-like-1}
\end{align*}
}

On the other hand, the second term of the right hand side of 
\eqref{eqn;time-like-0} can be treated by using the Minkowski inequality 
that
\algn{
  II
  = & \Big\|   \sum_{j\in\Z}  \sum_{k\le 2j} 
                2^{(1-1/2p)k}2^{sj}
              \Big\|
                    \psi_k\underset{(t)}{*}
                    \phi_j\underset{(x')}{*}
                     u(t,\cdot,\eta)
             \Big\|_{L^p(\re^{n-1}_{x'})}
        \Big\|_{L^1_t(\re_+)}
    \\
   \le & C   \sum_{j\in\Z} 2^{sj} \sum_{k\le 2j} 
                2^{(1-1/2p)k}
        \Big\| {\bk |\psi_k|\underset{(t)}{*}}
             \Big\|
                   \phi_j\underset{(x')}{*}
                    u(t,\cdot,\eta)
             \Big\|_{L^p(\re^{n-1}_{x'})}
        \Big\|_{L^1_t(\re_+)}
\\
   \le & C  \sum_{j\in\Z} 2^{sj} \sum_{k\le 2j} 2^{(1-1/2p)k}              
        \Big\|
             \Big\|
              \sum_{m\in \Z} \overline{\Phi_m}\underset{(x',\eta)}{*}
              \zeta_{j}(\eta)\underset{(\eta)}{*}  
              \zeta_{j-1}(\eta)\underset{(\eta)}{*} 
              \phi_j\underset{(x')}{*}
               u(t,\cdot,\eta)
             \Big\|_{L^p(\re^{n-1}_{x'})}
        \Big\|_{L^1_t(\re_+)}
\\
     \le & C\bigg\| \sum_{j\in\Z} 2^{sj} 2^{(2-1/p)j}               
             \Big\|
               \sum_{|m-j|\le1} \overline{\Phi_m}\underset{(x',\eta)}{*} 
               \zeta_{j}(\eta)\underset{(\eta)}{*} 
               \phi_j\underset{(x')}{*}
               u(t,\cdot,\eta)
             \Big\|_{L^p(\re^{n-1}_{x'})}
        \Big\|_{L^1_t(\re_+)}
\\
       \le & C\Big\|   \sum_{j\in\Z}  2^{(s+2-1/p)j}  
              \|\zeta_{j}\|_{L^{p'}(\re_{+,\eta})}            
             \Big\|
                 \overline{\Phi_j}\underset{(x',\eta)}{*} 
                 u(t,\cdot,\eta)
             \Big\|_{L^p(\re^{n}_{+})}
        \Big\|_{L^1_t(\re_+)}
\\
       \le & C\Big\| \sum_{j\in\Z}    2^{(s+2)j} 
             \Big\|
               (-\Del)^{-1}\overline{\Phi_j}\underset{(x',\eta)}{*} 
               \Del u(t,\cdot,\eta)
             \Big\|_{L^p(\re^{n}_{+})}
        \Big\|_{L^1_t(\re_+)}
\\
       \le & C\Big\|  \sum_{j\in\Z}  2^{sj}
             \Big\|
               \overline{\Phi_j}\underset{(x',\eta)}{*} 
               \Del u(t,\cdot,\eta)
             \Big\|_{L^p(\re^{n}_{+})}
        \Big\|_{L^1_t(\re_+)}
     = C\big\|\Del u \big\|_{ L^1(\re_+;\dB^s_{p,1}(\re^{n}_{+})) }.
    \eqntag\label{eqn;time-like-2}
}

One can apply the similar treatment for the spatial direction, 
for $1\le p\le \infty$. Let 
{$\Del u\in L^1(\re_+;\dB^{s}_{p,1}(\re^n_+))$} and
$\eta\in I_{-\ell}=[2^{-\ell},2^{-\ell+1})$ with  $\ell\in \Z$,
we insert the unity 
$\dsp \sum_{m\in\Z}\overline{\Phi_m}(x',\eta)\underset{(x',\eta)}{*}$ 
 to the estimate to see
{\allowdisplaybreaks
\begin{align*}
 \big\|u&(\cdot, \cdot,\eta) \big\|_{ L^1(\re_+;\dB^{s+2-1/p}_{p,1}(\re^{n-1}_{x'})) } 
 \\
   \le &C\Big\|\sum_{j\in \Z} 2^{sj}
                 2^{(2-1/p)j}          
                \Big\| \phi_j\underset{(x')}{*}
                       \zeta_{j-1}(\eta)\underset{(\eta)}{*}
                       \sum_{|m-j|\le 1}
                       \overline{\Phi_m}(x',\eta)\underset{(x',\eta)}{*}
                        u(t,\cdot,\eta)
                \Big\|_{ L^p(\re^{n-1}_{x'}) } 
        \Big\|_{L^1_t(\re_+)}
\\
   \le &C\Big\|\sum_{j\in \Z} 2^{(s+2)j} 
                \Big(   2^{-j}        
                \Big\|
                   \phi_j\underset{(x')}{*}
                   \zeta_{j-1}(\eta)\underset{(\eta)}{*}
                   \zeta_{j}(\eta)\underset{(\eta)}{*}
                   \overline{\Phi_j}(x',\eta)\underset{(x',\eta)}{*}
                   u(t,\cdot,\eta)
                \Big\|_{ L^p(\re^{n-1}_{x'}) }^p               
                \Big)^{1/p}
        \Big\|_{L^1_t(\re_+)}
\\
   \le &C\Big\|\sum_{j\in \Z}  2^{(s+2)j} 
                \Big( 2^{-j} 
                \Big(
                \big\|\zeta_{j}\big\|_{L^{p'}(\re_{+,\eta})} \\
       &\hskip25mm \times
                \Big\|\big\|\phi_j\underset{(x')}{*}
                      \zeta_{j-1}(\eta)\underset{(\eta)}{*}
                      \overline{\Phi_j}(x',\eta)\underset{(x',\eta)}{*}
                      u(t,\cdot,\eta)
                \big\|_{L^p(\re^{n-1}_{x'})}
                \Big\|_{L^p(\re_{+,\eta})}\Big)^p              
                \Big)^{1/p}
        \Big\|_{L^1_t(\re_+)}
  \\
   \le &C\Big\|\sum_{j\in \Z}  2^{(s+2)j} 
                \Big(   2^{-j} 2^j
                \Big\| \phi_j\underset{(x')}{*}
                       \overline{\Phi_j}(x',\eta)\underset{(x',\eta)}{*}
                       u(t,\cdot,\eta)
                \Big\|_{L^p(\re^{n}_+)}^p              
                \Big)^{1/p}
        \Big\|_{L^1_t(\re_+)} 
  \\
     \le &C\Big\|\sum_{j\in \Z} 2^{(s+2)j} 
            \Big\|\overline{\Phi_j}(x',\eta)\underset{(x',\eta)}{*}
                   u(s,\cdot,\eta)
            \Big\|_{L^p(\re^{n-1}_{x'}\times \re_{+,\eta})}
           \Big\|_{L^1_t(\re_+)}
 \\
     \le &C\Big\|\sum_{j\in \Z} 2^{sj}
            \Big\|\overline{\Phi_j}(\cdot, \cdot)\underset{(x',\eta)}{*}
                  \Del u(s,\cdot,\cdot)
            \Big\|_{L^p(\re^{n-1}_{x'}\times \re_{+,\eta})}
         \Big\|_{L^1_t(\re_+)}
     \le C\int_{\re_+}\|\Del u(t,x)\|_{\dB^s_{p,1}}dt.
 \addtocounter{equation}{1}
 \tag{\theequation} \label{eqn;space-like}  
\end{align*}
}
Combining the estimates \eqref{eqn;time-like-0}, \eqref{eqn;time-like-1}, 
\eqref{eqn;time-like-2} and \eqref{eqn;space-like}, 
 we conclude the estimate \eqref{eqn;sharp-trace-D}. 
\end{prf}

\begin{thm}[Sharp boundary derivative trace]
\label{thm;sharp-boundary-trace-N} 
For $1\le p\le \infty$ and $s\in \re$, there exists a constant $C>0$ such that for all
function
$ u=u(t,x',\eta)\in \dot{W}^{1,1}(\re_+; \dB^s_{p,1}(\re^{n}_+))$, 
$
{ \Del u\in  L^1(\re_+; \dB^{s}_{p,1}(\re^{n}_+))}
$ 
with 
$\pt_{\eta}u(0,x',\eta)=0$, it holds
\begin{align}
   \sup_{\eta\in\re_+}
   \Bigl(\big\|\pt_{\eta}u(\cdot,\cdot, \eta)
        &\big\|_{\dot{F}^{1/2-1/2p}_{1,1}(\re_+; \dB^s_{p,1}(\re^{n-1}))}
  +\big\|\pt_{\eta}u(\cdot,\cdot, \eta)
    \big\|_{L^1(\re_+; \dB^{s+1-1/p}_{p,1}(\re^{n-1}))}
  \Bigr)  \nonumber\\
   & \qquad\qquad  \le C\Bigl(
         \|\pt_t u\|_{L^1(\re_+;\dB^s_{p,1}(\re^{n}_+))}
       + \|\N^2  u\|_{L^1(\re_+;\dB^s_{p,1}(\re^{n}_+))}
       \Bigr) .\label{eqn;sharp-trace-N}
\end{align}
\end{thm}

Similar to the boundary trace estimate,
 Theorem \ref{thm;sharp-boundary-trace-N} immediately implies the 
optimality of the boundary derivative trace estimate for 
the Neumann boundary condition.
\par
\noindent
\begin{cor}[Sharp boundary trace for Neumann data]
\label{cor;heat-boundary-trace-N} 
For $1\le p\le \infty$ and $s\in \re$, let $u$ be a solution to the Cauchy problem 
\eqref{eqn;heat-N-0} in 
$$     \dot{W}^{1,1}(\re_+; \dB^s_{p,1}(\re^{n}_+)),
       \quad 
       {
        \Del u\in  L^1(\re_+; \dB^{s}_{p,1}(\re^{n}_+)),
       }
$$ 
then there exists a constant $C>0$ such that the boundary condition 
has to satisfy
\eqn{
   \|h\|_{\dot{F}^{1/2-1/2p}_{1,1}(\re_+; \dB^s_{p,1}(\re^{n-1}))}
  +\|h\|_{L^1(\re_+; \dB^{s+1-1/p}_{p,1}(\re^{n-1}))} 
  \le  \left\{
  \begin{aligned}
   & C\|\N^2  u\|_{L^1(\re_+;\dB^s_{p,1}(\re^{n}_+))}, \\
   & C\|\pt_t  u\|_{L^1(\re_+;\dB^s_{p,1}(\re^{n}_+))}.
  \end{aligned}
  \right. 
  }
\end{cor}

\begin{prf}{Theorem \ref{thm;sharp-boundary-trace-N}}
The proof is very similar to the proof for 
Theorem \ref{thm;sharp-boundary-trace-D}.
For $1\le p\le \infty$ and $s\in \re$, assume 
$u\in \dot{W}^{1,1}(\re_+;\dB^s_{p,1}{(\re^n_+)})$ 
with {$\Del u\in  L^1(\re_+;\dB^{s}_{p,1}(\re^n_+))$.} 
Then
{\allowdisplaybreaks
\begin{align*}
 \big\|\pt_{\eta}u(\cdot,\cdot, \eta) 
  &\big\|_{\dF^{1/2-1/2p}_{1,1}\big(\re_+;\dB^{s}_{p,1}(\re^{n-1}_{x'})\big)} 
 \\
  \le & 
  \Big\|{\bk \sum_{j\in\Z} 2^{sj}
              \sum_{k\ge 2j}  }
              2^{(1/2-1/2p)k}
              \big\|
                 \psi_k\underset{(t)}{*}
                 \phi_j\underset{(x')}{*}
                 \pt_{\eta}u(t,\cdot,\eta)
              \big\|_{L^p(\re^{n-1}_{x'})}
       \Big\|_{L^1_t(\re_+)} \\
  &+\Big\|{ \sum_{j\in\Z} 2^{sj}
                \sum_{k\le 2j}  }
            2^{(1/2-1/2p)k}
              \big\|
                 \psi_k\underset{(t)}{*}
                 \phi_j\underset{(x')}{*}
                 \pt_{\eta}u(t,\cdot,\eta)
              \big\|_{L^p(\re^{n-1}_{x'})}
       \Big\|_{L^1_t(\re_+)} 
 \\
 \equiv & I+II.
  \addtocounter{equation}{1}
 \tag{\theequation} \label{eqn;time-like-N0}
 \end{align*}
 }
For all $j\in \Z$, 
\eq{ \label{eqn;compatibility-condition-N}
  \pt_{\eta}u(0,x',\eta)=0,
  \quad (x',\eta)\in \re^{n-1}\times\re_+,
}
and it follows from \eqref{eqn;int-psi-k} that 
\begin{align*}
 \psi_k(t)\underset{(t)}{*} \pt_{\eta}u(t,x',\eta)
 =&\pt_t^{-1}\psi_k(t)\underset{(t)}{*} \pt_{\eta} \pt_t u(t,x',\eta).
 \addtocounter{equation}{1}
\tag{\theequation}  \label{eqn;time-like-int-N} 
\end{align*}

Hence from \eqref{eqn;compatibility-condition-N} and 
\eqref{eqn;time-like-int-N} and using the Hausdorff-Young 
inequality
{\allowdisplaybreaks
\begin{align*}
  I = & \Big\|{\bk \sum_{j\in\Z} 2^{sj}
                   \sum_{k\ge 2j}  }
                    2^{(1/2-1/2p)k}
              \big\|
                    \pt_t^{-1}\psi_k\underset{(t)}{*}
                    \phi_j\underset{(x')}{*}
                    \pt_{\eta}\pt_t u(t,x,\eta)
             \big\|_{L^p(\re^{n-1}_{x'})}
        \Big\|_{L^1_t(\re_+)}
  \\
 \le &C\Big\| {\bk\sum_{j\in\Z} 2^{sj}
                  \sum_{k\ge 2j}  2^{-(\frac12+\frac1{2p})k}
              }
               \int_{\re_+}              
                 |(\pt_t^{-1}\psi)_k(t-s)|  
                \big\|
                    \phi_j\underset{(x')}{*}
                    \pt_{\eta}\pt_t u(s,x,\eta)
                \big\|_{L^p(\re^{n-1}_{x'})}                 
                ds
          \Big\|_{L^1_t(\re_+)}
 \\
   \le &C\sum_{j\in \Z} 2^{sj}
          \Big\| {\bk \sum_{k\ge 2j}  2^{-(\frac12+\frac1{2p})k} } \\
            &\hskip24mm\times 
                 \big\|
                 \sum_{|m-j|\le 1} \overline{\Phi_m}\underset{(x',\eta)}{*} 
                 \phi_j\underset{(x')}{*}
                 \pt_{\eta}\zeta_{j}(\eta)\underset{(\eta)}{*}
                 \zeta_{j-1}(\eta)\underset{(\eta)}{*} 
                 \pt_t u(t,x,\eta)  
                \big\|_{L^p(\re^{n-1}_{x'})}                
          \big\|_{L^1_t(\re_+)}
  \\
   \le &C\sum_{j\in \Z} 2^{sj}
          \Big\|{\bk \sum_{k\ge 2j} 2^{-(\frac12+\frac1{2p})k}  }
                 \|\pt_{\eta}\zeta_{j}(\eta)
                 \|_{L^{p'}(\re_{+,\eta})}\\
          &\qquad\qquad\qquad\qquad \times       
                 \Big\|\big\|
                     \overline{\Phi_j}\underset{(x',\eta)}{*} 
                 \phi_j\underset{(x')}{*}
                \zeta_{j-1}(\eta)\underset{(\eta)}{*} 
                \pt_t u(t,x,\eta)  
                \big\|_{L^p(\re^{n-1}_{x'})}
                \Big\|_{L^p(\re^{+}_{\eta})}                
          \Big\|_{L^1_t(\re_+)}
   \\
      \le &C\Big\|\sum_{j\in \Z}  2^{sj}
                 \big\| 
                     \overline{\Phi_j}\underset{(x',\eta)}{*} 
                 \pt_t u(t,x,\eta)  
                \big\|_{L^p(\re^{n}_{+,(x',\eta)})}              
          \Big\|_{L^1_t(\re_+)} 
    \le C \|\pt_t u\|_{ L^1(\re_+;\dB^s_{p,1}(\re^n_{+}) }.
     \addtocounter{equation}{1}
 \tag{\theequation} \label{eqn;time-like-N1}
\end{align*}
}

For the second term of  \eqref{eqn;time-like-N0}, 
we use the Minkowski inequality to see
{\allowdisplaybreaks
\begin{align*}
  II 
   \le & C {\bk  \sum_{j\in\Z}  2^{sj} \sum_{k\le 2j} }
                2^{(1/2-1/2p)k}
        \Big\| {\bk |\psi_k|\underset{(t)}{*}}
             \big\|
                   \phi_j\underset{(x')}{*}
                    \pt_{\eta}u(t,x',\eta)
             \big\|_{L^p(\re^{n-1}_{x'})}
        \Big\|_{L^1_t(\re_+)}
\\
   \le & C {\bk  \sum_{j\in\Z} 2^{sj} \sum_{k\le 2j} 2^{(1/2-1/2p)k}} \\
       &\hskip20mm\times            
        \Big\|
             \big\|
              \sum_{m\in \Z} \overline{\Phi_m}\underset{(x',\eta)}{*} 
              \zeta_{j}(\eta)\underset{(\eta)}{*} 
              \zeta_{j-1}(\eta)\underset{(\eta)}{*} 
              \phi_j\underset{(x')}{*}
              \pt_{\eta}u(t,x',\eta)
             \big\|_{L^p(\re^{n-1}_{x'})}
        \Big\|_{L^1_t(\re_+)}
\\
     \le & C\Big\| {\bk  \sum_{j\in\Z}  2^{(s+1-1/p)j} }              
             \big\|
              \sum_{|m-j|\le1} \overline{\Phi_m}\underset{(x',\eta)}{*} 
              \pt_{\eta}\zeta_{j}(\eta)\underset{(\eta)}{*} 
              \zeta_{j-1}(\eta)\underset{(\eta)}{*} 
              \phi_j\underset{(x')}{*}
              u(t,x',\eta)
             \big\|_{L^p(\re^{n-1}_{x'})}
        \Big\|_{L^1_t(\re_+)}
\\
       \le & C\Big\|   \sum_{j\in\Z}    2^{(s+1-1/p)j}  
              \|\pt_{\eta}\widetilde{\zeta_{j}}\|_{L^{p'}(\re_{+,\eta})}             
             \big\|
               \overline{\Phi_j}\underset{(x',\eta)}{*} 
               u(t,x',\eta)
             \big\|_{L^p(\re^{n}_{+})}
        \Big\|_{L^1_t(\re_+)}
\\
       \le & C\Big\|   \sum_{j\in\Z}    2^{(s+2)j} 
             \big\|
               (-\Del)^{-1}\overline{\Phi_j}\underset{(x',\eta)}{*} 
              \Del u(t,x',\eta)
             \big\|_{L^p(\re^{n}_{+})}
        \Big\|_{L^1_t(\re_+)}
    \le  C\big\|\Del u \big\|_{ L^1(\re_+;\dB^s_{p,1}(\re^n_+)) }.
     \addtocounter{equation}{1}
\tag{\theequation} \label{eqn;time-like-N2}
\end{align*}
}
The estimate for the spatial direction is slightly simpler.
For $1\le p\le \infty$, 
{$\Del u\in L^1(\re_+;\dB^{s}_{p,1}(\re^n_+))$} and 
noting \eqref{eqn;Sum_Phi_is_1},
we estimate very similar way as in \eqref{eqn;space-like} to obtain
{\allowdisplaybreaks
\begin{align*}
 &\big\|\pt_{\eta}u(\cdot, \cdot,\eta) 
  \big\|_{ L^1\big(\re_+;\dB^{s+1-1/p}_{p,1}(\re^{n-1}_{x'})\big) } \\
   \le &C\Big\|\sum_{j\in \Z}
                 2^{(s+1-1/p)j}          
                \big\|\phi_j\underset{(x')}{*}
                \sum_{m\in\Z}
                \overline{\Phi_m}(x',\eta)\underset{(x',\eta)}{*}
                \pt_{\eta}u(t,x,\eta)
                \big\|_{ L^p(\re^{n-1}_{x'}) }  
        \Big\|_{L^1_t(\re_+)}
  \\
   \le & C\Big\|\sum_{j\in \Z}
                 2^{(s+1-1/p)j}          
                \big\| \phi_j\underset{(x')}{*}
                      \pt_{\eta}\zeta_{j}(\eta)\underset{(\eta)}{*}
                      \zeta_{j-1}(\eta)\underset{(\eta)}{*} 
                      \sum_{|m-j|\le 1}
                     \overline{\Phi_m}(x',\eta)\underset{(x',\eta)}{*}
                u(t,x,\eta)
                \big\|_{ L^p(\re^{n-1}_{x'}) }  
        \Big\|_{L^1_t(\re_+)}
\\
   \le &C\Big\|\sum_{j\in \Z} 2^{(s+1)j} 
              \Big(     
              \|\pt_{\eta}\zeta_{j}\|_{ L^{p'}(\re_{+,\eta}) }^p   \\
       &\hskip 25mm \times  
              \big\|
                   \phi_j\underset{(x')}{*}                   
                   \zeta_{j-1}(\eta)\underset{(\eta)}{*} 
               \overline{\Phi_j}(x',\eta)\underset{(x',\eta)}{*}
              u(t,x,\eta)
              \big\|_{ L^p(\re^{n}_{x'}\times\re_{+,\eta}) }^p 
                2^{-j}
                \Big)^{1/p}
        \Big\|_{L^1_t(\re_+)}
  \\
     \le &\;C\Big\|\sum_{j\in \Z} 2^{(s+2)j} 
            \big\| \overline{\Phi_j}(x',\eta)\underset{(x',\eta)}{*}
                   u(s,\cdot,\eta)
            \big\|_{L^p(\re^{n-1}_{x'}\times \re_{+,\eta})}
           \Big\|_{L^1_t(\re_+)} 
  \\
     \le &\;C\int_{\re_+}\|\Del u(t)\|_{\dB^s_{p,1}(\re^n_+)}dt.
     \eqntag\label{eqn;space-like-N} 
\end{align*}
}
Combining the estimates \eqref{eqn;time-like-N0}, \eqref{eqn;time-like-N1}, 
\eqref{eqn;time-like-N2} and \eqref{eqn;space-like-N}, we obtain 
the result \eqref{eqn;sharp-trace-N}.  
This completes the proof.
\end{prf}

%
%
%

\sect{Concluding remarks}
Since we establish maximal $L^1$-regularity for the inhomogeneous
boundary data of Dirichlet and Neumann boundary conditions, the estimate 
can be generalized into other boundary  conditions.
For instance, one can generalize  to the  case of the oblique boundary  condition:
\eq{ \label{eqn;parabolic-ob}
  \left\{
  \begin{aligned}
    &\pt_t u - \Del u=f(t,x), 
     &\qquad &t>0,\ &&x\in \re^n_+,\\
    &\quad \left. b\cdot \N u\right|_{x_n=0}=h(t,x'), 
     &\qquad &t>0,\ &&x'\in \re^{n-1}, \\
    &\quad \left. u\right|_{t=0}=u_0,
      &\qquad &\quad\ &&x\in \re^n_+,  
    \end{aligned}
  \right.
}
where $b=(b',b_n)$ is a given constant vector in $\re^n$ with $b_n\neq 0$.
We also obtain that the integral kernel to the Laplacian of the 
solution $\Psi_{ob}$ is given by 
\eq{
\Psi_{ob}(t,x)=
      c_{n+1}\int_{\re} \int_{\re^{n-1}}
      e^{i\t t+ix'\cdot \xi'}
       \frac{i\t }{ ib'\cdot \xi'- b_n\sqrt{i\tau+|\xi'|^2}}
       e^{-\sqrt{i\t+|\xi'|^2}x_n} d\xi'd\t.
} 
Then the analogous estimate to the boundary potential in 
Theorem \ref{thm;Naumann-BC-Lp-trace} can be obtained once 
we establish the almost orthogonal estimate as in Lemma \ref{lem;pt-orthogonal-2}.  
The detailed estimate is shown in the forthcoming 
paper \cite{OgSs21}.
\par
For a general domain $\Omega\subset\re^n$, for instance bounded domain with smooth boundary, 
we may generalize our results. In such a case, by standard decompositions of unity near the boundary 
$\partial\Omega$, it may be reduced into a problem in the half-space by a smooth diffeomorphism. 
Then if we establish maximal $L^1$-regularity for the parabolic initial-boundary value problems with lower order spatial derivative terms, we may extend the estimate for general domain cases. 
Further application to the nonlinear problem is also available.  We discuss 
such an application to the fluid mechanics in a forthcoming paper  \cite{OgSs20}
(cf. \cite{OgSs-JEPE}).  
%
%

\vskip2mm
{\bf Acknowledgments}
The authors are grateful to the anonymous referee for 
valuable suggestions and comments that improve the presentation of 
this paper largely.
The first author is partially supported by JSPS grant-in-aid for 
Scientific Research (S) \#19H05597, Scientific Research (B) \#18H01131
and Challenging Research (Pioneering) \#20K20284. 
The second author is partially supported by JSPS grant-in-aid 
for Scientific Research (B) \#16H03945 and (B)  \#21H00992 
and Fostering Joint International Research (B) \#18KK0072.

\vskip3mm\noindent
The authors declare that there is no conflict of interest. 
 

\end{document}